\documentclass[11pt,a4paper]{article}

\usepackage{amsmath,amssymb,amsfonts,amsthm,mathtools,mathrsfs}
\usepackage{geometry}
\usepackage{placeins}
\usepackage{booktabs,array,longtable}
\usepackage{microtype,needspace,xurl}
\usepackage{xcolor}
\usepackage{tikz}
\usepackage[numbers,sort&compress]{natbib}
\usepackage[colorlinks=true,linkcolor=blue!45!black,
  citecolor=green!35!black,urlcolor=blue!50!black]{hyperref}

\newtheorem{theorem}{Theorem}
\newtheorem{lemma}{Lemma}[section]
\newtheorem{proposition}[lemma]{Proposition}
\newtheorem{corollary}[lemma]{Corollary}
\theoremstyle{remark}

\allowdisplaybreaks[2]
\hypersetup{pdftitle={Fourth moments of prime exponential sums and exceptional sets for sums and differences of primes},pdfauthor={}}

\title{Fourth moments of prime exponential sums and exceptional sets for sums and differences of primes}
\author{Huixi Li
  \thanks{\raggedright Corresponding author. School of Mathematical Sciences and LPMC,
  Nankai University, Tianjin, China. Email
  \mbox{lihuixi@nankai.edu.cn}.}
\and
Zhixin Liu
  \thanks{\raggedright School of Mathematical Sciences, Tianjin University,
  Tianjin, China. Email \mbox{zhixinliu@tju.edu.cn}.}
}
\date{}

\begin{document}
\maketitle

\begin{abstract}
We prove a square root bound for the exceptional sets in the binary
Goldbach and prime difference problems under quasi Riemann Hypothesis. If every nontrivial zero of every primitive
Dirichlet $L$-function has real part at most a fixed
$1/2\leq\Theta<11/12$, then
$E_+(N)\ll_\Theta N^{1/2}(\log N)^{100}$ and $E_-(N)\ll_\Theta N^{1/2}(\log N)^{100}$.
Here $E_+(N)$ counts the even integers greater than $2$ and at most $N$
that are not sums of two primes, and $E_-(N)$ counts the positive even
integers at most $N$ that are not differences of two primes.
The common analytic estimate is a fourth moment bound of the minor arcs for a smoothed prime exponential
sum. 
\end{abstract}

\noindent\textbf{Mathematics Subject Classification (2020).}
Primary 11P32; Secondary 11L07, 11M26.

\noindent\textbf{Keywords.}
Goldbach Conjecture, quasi Riemann Hypothesis, prime exponential sums, fourth moments.

\section{Introduction}
\label{sec:introduction}

The binary Goldbach conjecture asserts that every even integer greater
than $2$ is the sum of two primes. Several different weakenings of this
statement have led to substantial progress. By allowing one of
the two summands to be almost prime, Chen~\cite{Chen1973} proved
in 1973 that every sufficiently
large even integer is the sum of a prime and an integer with at most
two prime factors. For surveys of
Goldbach's problem and related additive questions, see ~\cite{KT2005, Vaughan2016, Li2023}.

Another interesting question asks how many integers can fail to have such a representation. 
Let
\begin{align*}
 E_+(N)&=\#\{2<n\leq N:2\mid n,\ n\ne p_1+p_2\},\\
 E_-(N)&=\#\{2\leq h\leq N:2\mid h,\ h\ne p_2-p_1\},
\end{align*}
where the variables on the right are primes and need not be consecutive.
Thus the binary Goldbach conjecture is equivalent to $E_+(N)=0$ for
every $N$. The analogous difference problem asks whether every positive
even integer is a difference of two primes. The definition of $E_-(N)$
places no size restriction on the primes.

The first unconditional results for $E_+(N)$ were obtained independently
by van der Corput~\cite{VanDerCorput1937},
Chudakov~\cite{Chudakov1937,Chudakov1938}, and
Estermann~\cite{Estermann1938} in 1937-1938. Their work gives
$E_+(N)\ll_A N(\log N)^{-A}$ for every fixed $A>0$.
Vaughan~\cite{Vaughan1972} obtained
$E_+(N)\ll N\exp(-c\sqrt{\log N})$ for some $c>0$ in 1972.
The passage to a power saving was made by Montgomery and
Vaughan~\cite{MontgomeryVaughan1975} in 1975: there is a fixed
$\delta>0$ such that $E_+(N)\ll N^{1-\delta}$.
Their argument uses estimates for the density of zeros of Dirichlet $L$-functions,
building on Gallagher~\cite{Gallagher1970Density}, and treats both
the presence and the absence of an exceptional real zero.

Subsequent work made the exponent explicit and reduced it.
Chen and Pan~\cite{ChenPan1980} in 1980, Chen~\cite{Chen1983} in 1983, and
Chen and Liu~\cite{ChenLiu1989} in 1989 obtained the exponents $0.99$,
$0.96$, and $0.95$, respectively. Li~\cite{Li1999Goldbach,Li2000Goldbach}
reduced the exponent to $0.921$ in 1999 and $0.914$ in 2000,
and Lu~\cite{Lu2010} obtained $0.879$ in 2010.
Pintz~\cite{Pintz2018II} obtained $E_+(N)\ll N^{0.72}$ in 2018,
using an explicit formula which retains the contribution of zeros
close to $1$. Zhao~\cite{Zhao2026} obtains
$E_+(N)\ll N^{7/10}$ in 2026.
Bhowmik and Grimmelt~\cite{BhowmikGrimmelt2026} explain the development
of the methods based on the density of zeros and explicit formulas and give a smoothed
formula on the major arcs.

The distribution of exceptions within short intervals is a separate
question. For $H=X^\vartheta$ with $1/5<\vartheta\leq1$,
Matom\"aki~\cite{Matomaki2008} proved in 2008 that
\[
 E_+(X+H)-E_+(X)\ll H^{1-\delta},
\]
where $\delta=\delta(\vartheta)>0$.
Grimmelt~\cite{Grimmelt2022} proved in 2022 that all but
$O_{A,\varepsilon}(H(\log X)^{-A})$ even integers in $[X-H,X]$
are Goldbach numbers when $X^{7/120+\varepsilon}<H\leq X$,
for every fixed $A,\varepsilon>0$.

With respect to the asymptotic formula the number of representations of even $n$, Hardy and
Littlewood~\cite{HardyLittlewood1923} conjectured in 1923 that
\[
 \sum_{p_1+p_2=n}1\sim\mathfrak S(n)\frac{n}{(\log n)^2},
\]
where, for an even nonzero integer $r$,
\begin{equation*}\label{eq:intro-singular-series}
 \mathfrak S(r)=
 2\prod_{p>2}\left(1-\frac1{(p-1)^2}\right)
   \prod_{\substack{p\mid r\\p>2}}\frac{p-1}{p-2}.
\end{equation*}
The same singular series occurs in the expected asymptotic for prime
pairs of difference $h$. Pintz~\cite{Pintz2018I} treats the two
problems in a common explicit formula in 2023; the
unconditional difference bound $E_-(N)\ll N^{\vartheta}$ for some
$\vartheta<2/3$ is given by Pintz~\cite[Theorem~17]{Pintz2009}.

Matom\"aki, Radziwi\l\l\ and
Tao~\cite[Theorem~1.3]{MRT2019} proved in 2019 that the
Hardy--Littlewood asymptotics for  correlations between primes hold for
almost all shifts $|h-h_0|\leq H$, where
$0\leq h_0\leq X^{1-\varepsilon}$ and
$X^{8/33+\varepsilon}\leq H\leq X^{1-\varepsilon}$.
Their Goldbach assertion concerns targets in $[X,X+H]$.
For every fixed $A>0$ and $0<\varepsilon<1/2$, the error is
$O_{A,\varepsilon}(X(\log X)^{-A})$ outside
$O_{A,\varepsilon}(H(\log X)^{-A})$ exceptions.
This extends the earlier $X^{1/3+\varepsilon}$ range of
Mikawa~\cite{Mikawa1991} and
Perelli and Pintz~\cite{PP1993}.
See other related works on gaps between primes. Goldston, Pintz and
Y\i ld\i r\i m~\cite{GPY2009} proved that consecutive prime
gaps can be arbitrarily small relative to the average spacing.
Zhang~\cite{Zhang2014}, Maynard~\cite{Maynard2015}, and the
Polymath project~\cite{Polymath2014} obtained bounded gaps and
prime clusters in bounded intervals. Green and
Tao's~\cite{GreenTao2006,GreenTao2010} restriction estimates and
work on linear equations in primes provide a broader setting.

The connection with zeros also has a converse direction.
Languasco and Zaccagnini~\cite{LZ2012} obtain an average
formula for Goldbach representations under the Riemann Hypothesis.
Bhowmik and Ruzsa~\cite{BhowmikRuzsa2018}, completing an argument
of Granville~\cite{Granville2007}, show that a power saving in
the average error implies a fixed half plane free of zeros for
$\zeta(s)$. 
There are also connections
with exceptional real zeros. Friedlander, Goldston, Iwaniec and
Suriajaya~\cite{FGIS2022} show that suitable upper and lower Goldbach
bounds for every sufficiently large even target rule out such
zeros; see also Goldston and Suriajaya~\cite{GoldstonSuriajaya2021}.
Matom\"aki and Merikoski~\cite{MatomakiMerikoski2023} study
Goldbach asymptotics under the existence of a
Siegel zero. With respect to the Elliott--Halberstam conjecture, Huang and Li~\cite[Corollary~1]{HuangLi2020}
prove that every sufficiently large even integer $N$ is the sum of two
primes if the Elliott--Halberstam estimate for $\Lambda(n)$ and its
analogue twisted by the M\"obius function for $\Lambda(n)\mu(N-n)$ hold with fixed
levels of distribution $\vartheta,\vartheta'\in(0,1)$ satisfying
$\vartheta+\vartheta'>1$.

In the present paper we obtain a better exceptional bound with square root power, under the quasi Riemann Hypothesis. 
For $1/2\leq\Theta<1$, let $(Q_\Theta)$ denote the assertion
\begin{equation*}\label{eq:Q-theta}
 \Re\rho\leq\Theta
 \quad\text{for every nontrivial zero of every primitive
 Dirichlet $L$-function}.
\end{equation*}
At $\Theta=1/2$ this is the Generalized Riemann Hypothesis.
Hardy and Littlewood~\cite{HL1924} proved in 1924 that GRH implies
$E_+(N)\ll_\varepsilon N^{1/2+\varepsilon}$.
Goldston~\cite{Goldston1992} sharpened this under GRH in 1992 to
\begin{equation*}
 E_+(N)\ll N^{1/2}(\log N)^3.
\end{equation*}

For $1/2\leq\Theta<3/4$, the classical method gives
$E_+(N)\ll_{\Theta,\varepsilon}N^{2\Theta-1/2+\varepsilon}$
under $(Q_\Theta)$; see
Bhowmik and Grimmelt~\cite[Theorem~3.1]{BhowmikGrimmelt2026}.
In 2026, Pintz and Alth\"ofer~\cite[Theorem~1]{PintzAlthofer2026}
proved, under $(Q_\Theta)$ for any fixed $\Theta<1$, that
\[
 E_+(N)\ll_{\Theta,\varepsilon}N^{7/12+\varepsilon}.
\]
They combine Pintz's formula on the major arcs with the exponential sum over primes
estimate of Maynard, Pandey and Radziwi\l\l~\cite{MPR2026}.
For $1/2\leq\Theta\leq2/3$, their Theorem~3 gives the refined
conditional bound
\begin{equation*}
 E_+(N)\ll_{\Theta,\varepsilon}N^{g(\Theta)+\varepsilon},
 \qquad
 g(\Theta)=\frac{\Theta(5-4\Theta)}{2(2-\Theta)}.
\end{equation*}
Here $g(1/2)=1/2$ and $g(2/3)=7/12$.
We prove that the square root power is valid throughout the fixed
range $1/2\leq\Theta<11/12$.
The same method gives this conclusion for both
$E_+(N)$ and $E_-(N)$.

For the main parameters, let $X\geq2$ and $N\geq2$.
Write $e(t)=\exp(2\pi it)$ and $L=\log X$.
Auxiliary lengths below are denoted by $\mathsf{M}$, with subscripts
when needed, to distinguish them from the arithmetic cutoff $N$. 
Fix a real weight $w\in C_c^\infty((0,1))$ such that
$0\leq w\leq1$ and $w=1$ on $[1/8,7/8]$, chosen from the
logarithmic Gevrey class of order $2$ specified in
Section~\ref{sec:fixedlog-smoothing}. Put
\[
 S_w(\alpha)=\sum_p(\log p)w(p/X)e(p\alpha),\qquad
 S_{\Lambda,w}(\alpha)=\sum_n\Lambda(n)w(n/X)e(n\alpha).
\]
For $P=X^\kappa$, define the major and minor arcs on
$\mathbb{R}/\mathbb{Z}$ by
\begin{align*}
 \mathfrak M(P)&=\bigcup_{q\leq P}\ 
   \bigcup_{\substack{a\bmod q\\(a,q)=1}}
   \left\{\alpha:\left|\alpha-\frac aq\right|
                       \leq\frac{P}{qX}\right\},\\
 \mathfrak m(P)&=(\mathbb{R}/\mathbb{Z}) \setminus\mathfrak M(P).
\end{align*}
For sufficiently large $X$, these arcs are disjoint, since $2P^2<X$.
Our main analytic result is the following.

\begin{theorem}\label{thm:moment-main}
Assume $(Q_\Theta)$ for a fixed $1/2\leq\Theta<11/12$.
There is a fixed $\kappa=\kappa(\Theta)\in(4/9,9/20)$ such that
\begin{equation}\label{eq:moment-main}
 \int_{\mathfrak m(X^\kappa)}|S_w(\alpha)|^4\,d\alpha
 \ll_{w,\Theta}X^{5/2}(\log X)^{100}.
\end{equation}
\end{theorem}

\begin{theorem}\label{thm:applications}
Assume $(Q_\Theta)$ for a fixed $1/2\leq\Theta<11/12$. Then
\begin{equation*}
 E_+(N)\ll_\Theta N^{1/2}(\log N)^{100}, \quad \text{and } E_-(N)\ll_\Theta N^{1/2}(\log N)^{100}.
\end{equation*}
\end{theorem}

For every fixed $1/2<\Theta<11/12$, our Goldbach bound
attains the same $N^{1/2}$ power as Goldston's estimate under
GRH, while requiring only the weaker quasi Riemann Hypothesis. 
Note that OpenAI~\cite{OAI2026September30} proved the quasi Riemann hypothesis with a $7/8$ zero-free boundary, which was subsequently refined by Liu~\cite{Liu2026} and by Xu and Sang~\cite{XuSang2026}.
For $11/12 \leq \Theta <1$, by a similar method we can also improve the exponent $7/12+\epsilon$. We leave this to the interested reader. 

By additive
orthogonality,
\begin{align*}
 R_+(n;X)
 &=\int_0^1S_w(\alpha)^2e(-n\alpha)\,d\alpha\notag\\
 &=\sum_{p_1+p_2=n}(\log p_1)(\log p_2)
                  w(p_1/X)w(p_2/X),
 \\
 R_-(h;X)
 &=\int_0^1|S_w(\alpha)|^2e(-h\alpha)\,d\alpha\notag\\
 &=\sum_{p_2-p_1=h}(\log p_1)(\log p_2)
                  w(p_1/X)w(p_2/X).
\end{align*}
The corresponding coefficients on the minor arcs $R_{+,\mathfrak m}$ and
$R_{-,\mathfrak m}$ satisfy
\begin{equation}\label{eq:intro-bessel}
 \sum_{r\in\mathbb Z}|R_{+,\mathfrak m}(r;X)|^2
 =\sum_{r\in\mathbb Z}|R_{-,\mathfrak m}(r;X)|^2
 =\int_{\mathfrak m(P)}|S_w(\alpha)|^4\,d\alpha.
\end{equation}
The terms from the major arcs are bounded below by a positive constant times
$X$ for even $n\in[X/2,X]$ and even $h\in[2,X/2]$.
Consequently, a missing representation contributes $\gg X^2$ to
the relevant sum in \eqref{eq:intro-bessel}. This proves the passage
from Theorem~\ref{thm:moment-main} to
Theorem~\ref{thm:applications}.

Guth and Maynard~\cite{GuthMaynard2026} developed estimates for large values of Dirichlet polynomials, and Chen, Gupta and
Li~\cite{CGL2026} obtained estimates in the character aspect for their
energy and the density of Dirichlet zeros.

The fourth moment is the point at which we use more than a
pointwise estimate on the minor arcs. The estimate of Maynard, Pandey and
Radziwi\l\l~\cite{MPR2026}, combined with the second moment,
gives
\[
 \int_{\mathfrak m(P)}|S_w(\alpha)|^4\,d\alpha
 \ll_\varepsilon X^{31/12+\varepsilon}
 \qquad(P=X^\kappa,\ \kappa>4/9).
\]
This is the input on the minor arcs leading to the exponent $7/12$ above.
To reach $5/2$ under $(Q_\Theta)$, we keep the energy and the weighted
square mass of the same set of large values throughout the character
expansion and Mellin inversion.

The smoothing estimates are developed in
Section~\ref{sec:fixedlog-smoothing}, and the estimates for the density of zeros
are recalled in Section~\ref{dens:section}.
Section~\ref{sec:major-principal} treats the major arcs and the
complete principal term. The estimates for difference energy and factor moments are proved in
Sections~\ref{sec:difference-energy-fixed-log} and
\ref{sec:original-factor-fourth-fixed-log}, and are combined with
the mass estimates for complete values in
Sections~\ref{sec:ab-complete-bridge}--\ref{sec:ar-subset-deficits}.
The remaining sections treat the terminal range and complete the
proofs of the two main theorems.

We use $\mu$, $\Lambda$ and $\varphi$ for the M\"obius, von Mangoldt
and Euler totient functions. Dirichlet convolution is denoted by $*$,
and $d_j=1^{*j}$. A star on a character sum restricts it to primitive
characters. The Gauss sum is
\[
 \tau_q(\chi)=\sum_{a\bmod q}\chi(a)e(a/q).
\]
Unless otherwise stated, character sums and character relations use
one fixed original modulus and include induced characters.
The implied constants in $\ll$ may depend on fixed weights and on
the parameters indicated. We write $A\preccurlyeq B$ when
$A\ll_\varepsilon X^\varepsilon B$ for every $\varepsilon>0$.
Dyadic length parameters lie in a fixed polynomial range; a factor
supported at $1$ has length scale one.

\section{Smoothing and a fixed logarithmic value bound}
\label{sec:fixedlog-smoothing}
\label{sec:endpoint-prime-cap}

The smoothing functions restrict all sums to integers of size comparable
to the large parameter. A function in $C_c^\infty((0,1))$ is infinitely
differentiable and vanishes outside a closed subinterval of $(0,1)$.
Consequently integrations by parts involving these cutoffs have zero
endpoint terms.
The condition that the weight equals one on a smaller interval ensures
that smoothing retains a substantial part of the Goldbach convolution.
The auxiliary weight will let us insert a smooth cutoff into a complete
Dirichlet polynomial without changing the summands under consideration.

Fix \(1<s\leq2\). Choose \(w\in C_c^\infty((0,1))\),
\(0\leq w\leq1\), with \(w=1\) on \([1/8,7/8]\), and
\(u\in C_c^\infty((0,\infty))\) equal to one near \(\operatorname{supp}w\).
Choose both in logarithmic coordinates from a fixed Gevrey class
of order \(s\), so that, for \(g(v)=u(e^v)\),
\[
 \|g^{(j)}\|_\infty\leq C_0C_1^j(j!)^s\qquad(j\geq0).
\]
We fix the order $s=2$ in the final logarithmic estimate.
We use \(\widehat f(z)=\int_0^\infty f(y)y^{z-1}\,dy\).

The change of variables $y=e^v$, $dy/y=dv$, gives
\begin{equation*}
 \widehat f(\sigma+it)
 =\int_{\mathbb R}f(e^v)e^{\sigma v}e^{itv}\,dv.
\end{equation*}
Thus the Mellin transform is a Fourier transform in the logarithm of
the variable. In particular, Mellin inversion takes the form
\begin{equation*}
 f(y)=\frac1{2\pi}\int_{\mathbb R}
          \widehat f(\sigma+it)y^{-\sigma-it}\,dt.
\end{equation*}
Its usefulness for products comes from
$(y_1\cdots y_k)^{-it}=\prod_i y_i^{-it}$: a smooth function of a
product becomes an integral of products of functions of one variable.
The Gevrey condition controls how the derivative bounds grow with the
number of differentiations. It gives a quantified decay estimate, rather
than a separate unspecified constant for each fixed number of
integrations by parts.

\begin{lemma}\label{lem:log-gevrey}
Uniformly for $0\le\sigma\le1$,
\[
 |\widehat u(\sigma+iy)|\le C\exp(-c|y|^{1/s}).
\]
If $f_\xi(y)=w(y)e(X\xi y)$, with $w(e^v)$ a fixed Gevrey
cutoff of the same order, and $T=\max(1,X|\xi|)$, then for
$|t|\ge C_2T$ one also has
\[
 |\widehat f_\xi(it)|\le C\exp(-c|t|^{1/s}).
\]
Consequently, for every fixed $A>0$, all additive Mellin frequencies
outside $|t|\le C_2T+C_A L^s$ have total absolute kernel mass
$O_A(X^{-A})$.  The height is additive, not $TL^s$.
\end{lemma}

\begin{proof}
The Mellin integral is the Fourier transform of
$e^{\sigma v}g(v)$ on a fixed compact interval.  Its $j$th derivative
is bounded by $C^{j+1}(j!)^s$, uniformly in $\sigma$.
Integrating by parts $j$ times and choosing $j\asymp |y|^{1/s}$,
we obtain the first assertion.

More explicitly, the fixed support and the derivative estimate give
\begin{equation*}
 |\widehat u(\sigma+iy)|
 \le C\left(\frac{C_3}{|y|}\right)^j(j!)^s
 \le C\left(\frac{C_3j^s}{|y|}\right)^j.
\end{equation*}
For large $|y|$, take $j$ to be the integer part of
$c_3|y|^{1/s}$, with $c_3$ small enough that the parenthesis is
at most $1/2$. The resulting bound is $C2^{-j}$, which has the
required exponential decay. Bounded $|y|$ is covered by the
integral of the absolute value of the amplitude.

For the second, write its phase as
$t\phi_\lambda(v)$, where
$\phi_\lambda(v)=v+\lambda e^v$ and
$\lambda=2\pi X\xi/t$.  If $C_2$ is sufficiently large in
terms of the fixed support, then
$1/2\le\phi_\lambda'(v)\le3/2$ on that support.
The inverses of these analytic changes of variables have uniform
analytic derivative bounds.  After changing variables to
$z=\phi_\lambda(v)$, the amplitude, extended by zero, is a
Gevrey function with uniform constants and support in a common
compact interval. Applying the first argument with
Fourier parameter $t$ gives the required decay. Integrating the
exponential tail introduces only a polynomial factor in $L$, which
is absorbed by increasing $C_A$.

The condition on the phase derivative prevents a stationary point,
where the first derivative of the phase vanishes and integration by
parts would divide by zero. Beyond the stated height, the tail is
bounded by a constant times
\begin{equation*}
 \int_{C_A L^s}^{\infty}e^{-ct^{1/s}}\,dt
 =s\int_{C_A^{1/s}L}^{\infty}v^{s-1}e^{-cv}\,dv.
\end{equation*}
This is smaller than any prescribed negative power of the large scale
once the fixed constant is chosen sufficiently large. A logarithmic
amount is therefore added to the height already required by the
additive oscillation.
\end{proof}

A Dirichlet character modulo $q$ is a multiplicative periodic function
which is zero on integers not coprime to $q$. The principal character
has value one at every invertible residue. A primitive character is one
which does not arise from a smaller modulus. These are different
properties: a nonprincipal character can be induced from a smaller
primitive nonprincipal character. If its conductor is $r$, induction
means
\begin{equation*}
 \chi(n)=\chi^*(n)\mathbf1_{(n,q)=1},\qquad r\mid q.
\end{equation*}
The primitive character of conductor one induces the principal
character at every modulus. We first estimate primitive nonprincipal
characters, where the explicit formula has no pole at one, and then
account for the extra integers removed by induction.

\begin{lemma}\label{lem:endpoint-prime-cap}
Let $u\in C_c^\infty((0,\infty))$ be fixed.  Suppose that every
nontrivial zero of every primitive Dirichlet $L$-function has real
part at most a fixed $1/2\leq\Theta<1$.  For $X\geq2$, any nonprincipal character
$\chi$ modulo $q$, and every real $t$, put
\[
 \Psi_u(t,\chi)=\sum_n\Lambda(n)u(n/X)\chi(n)n^{-it}.
\]
Then
\begin{equation*}
 |\Psi_u(t,\chi)|
 \ll_u X^\Theta\log\bigl(2q(2+|t|)\bigr).
\end{equation*}
In particular, for every fixed $C>0$, this is
$O_{u,C,\Theta}(X^\Theta\log X)$ uniformly when
$q(2+|t|)\leq X^C$.  The conclusion applies in particular to
the fixed Gevrey weights used in the endpoint reconstruction estimate.
\end{lemma}

\begin{proof}
For a primitive nonprincipal \(\chi^*\) of conductor \(r>1\),
Mellin inversion gives
\[
 \Psi_u(t,\chi^*)=-\frac1{2\pi i}\int_{(2)}
 \frac{L'}{L}(s+it,\chi^*)X^s\widehat u(s)\,ds.
\]

Here the integral is taken upwards on the vertical line with real
part two. Insert Mellin inversion into the sum and use the absolutely
convergent Dirichlet series
\begin{equation*}
 -\frac{L'}{L}(z,\chi^*)
 =\sum_{n\ge1}\frac{\Lambda(n)\chi^*(n)}{n^z},
 \qquad \Re z>1.
\end{equation*}
The von Mangoldt function is $\log p$ at $p^j$, for a prime $p$ and
an integer $j\ge1$, and is zero otherwise. Its appearance is exactly
what makes the logarithmic derivative available. Analytic continuation,
the functional equation and the local zero count used below are
standard inputs about Dirichlet $L$-functions; they do not follow from
elementary character orthogonality.
The step that shifts the contour uses the residue theorem from complex analysis.
One first integrates around a rectangle with the old and new vertical
lines as opposite sides. The difference of those line integrals is
$2\pi i$ times the sum of the residues at poles inside the rectangle,
up to its horizontal sides. A residue is the coefficient of
$(s-s_0)^{-1}$ in the local expansion at a pole $s_0$.
At a zero of multiplicity $m$ the logarithmic derivative has residue
$m$, so the minus sign before it makes the zero contribution negative.
Taking the rectangle heights to infinity, with the decay stated below,
gives the zero sum. This use of the residue theorem converts analytic
information about zeros into an estimate for the original finite sum.

Shift to \(\Re s=-1/2\). The nontrivial zeros contribute
\(-X^{\rho-it}\widehat u(\rho-it)\); the possible trivial zero at
\(s=-it\) contributes \(O_u(1)\). The functional equation gives
\[
 |L'/L(-1/2+i(t+v),\chi^*)|
 \ll\log\bigl(2r(2+|t+v|)\bigr).
\]
Rapid Mellin decay controls the new vertical integral and makes the
horizontal integrals vanish along heights avoiding zeros. Hence
\begin{equation*}
 \Psi_u(t,\chi^*)=-\sum_\rho X^{\rho-it}\widehat u(\rho-it)
 +O_u\bigl(\log(2r(2+|t|))\bigr).
\end{equation*}
Uniformly for \(0\leq\sigma\leq1\),
\(|\widehat u(\sigma+iv)|\ll_u(1+|v|)^{-3}\).
The local zero count, with multiplicity, is
\(\#\{\rho:|\Im\rho-y|\leq1\}\ll\log(2r(2+|y|))\)
\cite{IK2004}. Thus the zero sum is absolutely convergent and
\[
 \sum_\rho|X^{\rho-it}\widehat u(\rho-it)|
 \ll_u X^\Theta\sum_{j\in\mathbb Z}
 \frac{\log(2r(2+|t|+|j|))}{(1+|j|)^3}
 \ll_u X^\Theta\log(2r(2+|t|)).
\]

For this summation, divide the zero ordinates into intervals of length
one at distance about $|j|$ from $t$. Each interval contains at most a
logarithmic number of zeros, while the Mellin weight is
$O((1+|j|)^{-3})$. Moreover
\begin{equation*}
 \log(2r(2+|t|+|j|))
 \le \log(2r(2+|t|))+\log(1+|j|).
\end{equation*}
Both resulting series converge. The half plane free of zeros enters
through $|X^{\rho-it}|=X^{\Re\rho}\le X^\Theta$.
For the induction step, the support condition restricts the powers of
any one prime to a fixed multiplicative interval. Their number is
therefore bounded independently of the large scale.
If \(\chi\bmod q\) is induced by \(\chi^*\), the difference
consists of prime powers \(p^j\) with \(p\mid q\), \(p\nmid r\).
For \(\operatorname{supp}u\subset[c_1,c_2]\), there are at most
\(1+\log(c_2/c_1)/\log p\) such exponents per prime. Therefore
\[
 |\Psi_u(t,\chi)-\Psi_u(t,\chi^*)|
 \ll_u\sum_{p\mid q}(\log p+1)\ll_u\log(2q),
\]
which proves the assertion.

\end{proof}

For the complete values normalized by Gauss sums used below, the inducing
formula for Gauss sums gives $|\tau_q(\overline\chi)|/\sqrt q\leq1$.
Their magnitude therefore obeys the same bound. Product
thresholds are chosen below these complete levels, so they also obey
the bound without any upper bound on an individual product summand.

\section{Estimates for the density of zeros}
\label{dens:section}

For $Q,T\geq1$, put
\begin{equation*}\label{dens:family-count}
 {\mathscr{Z}}_*(\sigma;Q,T)=\sum_{q\leq Q}\sum_{\chi\bmod q}^{*}
 \#\{\rho:L(\rho,\chi)=0,\ \Re\rho\geq\sigma,
                                  \ |\Im\rho|\leq T\},
\end{equation*}
where only nontrivial zeros are counted, with analytic multiplicity.
The primitive character of conductor one is included; its $L$-function
is $\zeta$, and the pole at one is not counted.

Heath--Brown's Theorems~1--3 \cite{HeathBrown1979Density} give,
respectively,
\begin{align}
 {\mathscr{Z}}_*(\sigma;Q,T)&\ll_\varepsilon
 (Q^2T)^{2(1-\sigma)+\varepsilon},
 &&11/14\leq\sigma\leq1,\label{dens:density-one}\\
 {\mathscr{Z}}_*(\sigma;Q,T)&\ll_\varepsilon
 (Q^2T^{6/5})^{(20/9)(1-\sigma)+\varepsilon},
 &&1/2\leq\sigma\leq1,\label{dens:density-two}\\
 {\mathscr{Z}}_*(\sigma;Q,T)&\ll_\varepsilon
 (Q^2T^2)^{2(1-\sigma)+\varepsilon},
 &&129/167\leq\sigma\leq1.\label{dens:density-three}
\end{align}
In \eqref{dens:density-two} we use the $A(\sigma)=20/9$ branch,
namely equation~(3) of Theorem~2. For each fixed $\varepsilon>0$,
these estimates are uniform in $Q,T\geq1$ and in the displayed
ranges of real parts, including $Q=1$. Bounded smaller conductor or height
parameters may be increased to $2$. The loss is an arbitrary small
power, not a fixed power of a logarithm.

\begin{proposition}[Ingham's zeta density estimate]
\label{dens:lower-zeta-density}
For each fixed $\varepsilon>0$, uniformly for $T\geq2$ and
$1/2\leq\sigma\leq11/14$,
\begin{equation*}
 {\mathscr{Z}}_\zeta(\sigma,T)
 \ll_\varepsilon T^{3(1-\sigma)/(2-\sigma)+\varepsilon}.
\end{equation*}
Here ${\mathscr{Z}}_\zeta(\sigma,T)$ counts the nontrivial zeta zeros
with $\Re\rho\geq\sigma$ and $|\Im\rho|\leq T$, with multiplicity.
\end{proposition}
This is the form of Ingham's estimate allowing an arbitrarily small power loss in
\cite[Theorem~9.19(B), pp.~236--237]{Titchmarsh1986Zeta}.
The convention using positive ordinates and strict inequalities for real parts in the reference
gives the stated convention by conjugation and a shift of the real part of
$1/\log T$. The total zero count covers the endpoint and a
$1/\log T$ neighbourhood of $1/2$; the fixed logarithmic loss is
absorbed into $T^\varepsilon$.

\section{Major arcs and the complete term from the principal character}
\label{sec:major-principal}

Near a rational number with a small denominator, character
orthogonality replaces additive oscillation by multiplicative
characters. The principal character supplies the expected main term;
zeros of the other $L$-functions measure the error. The major arcs are
chosen large enough for the remaining minor arcs to admit stronger
bounds, but their size is limited by the available estimates for the density of zeros. The estimates for a family of characters and those for the
character of conductor one have different roles in this section.

\subsection{Major arc notation}

Put $\Theta_0=\max(\Theta,7/8)$. The major arcs use the full
hypothesis on the primitive family with $\Theta<1$. The principal estimate
below uses its consequence for $\zeta$ and assumes
$\Theta<11/12$.

Fix a real function \(w\in C_c^\infty((0,1))\) such that
\(0\leq w\leq1\) and \(w=1\) on \([1/8,7/8]\), and write
\[
 S_w(\alpha)=\sum_p(\log p)w(p/X)e(p\alpha),\qquad
 T_w(\alpha)=\sum_n w(n/X)e(n\alpha).
\]
For \(P=X^\kappa\), define the arcs on \(\mathbb{R}/\mathbb{Z}\) by
\begin{align*}
 \mathfrak M(P)&=\bigcup_{q\leq P}\,
 \bigcup_{\substack{a\bmod q\\(a,q)=1}}
 \left\{\alpha:\left|\alpha-\frac aq\right|
                 \leq\frac{P}{qX}\right\},\\
 \mathfrak m(P)&= (\mathbb{R}/\mathbb{Z}) \setminus\mathfrak M(P).
\end{align*}
For fixed \(\kappa<1/2\), the individual major arcs are disjoint
when \(X\) is sufficiently large.  Later \(Q\) denotes a dyadic
denominator size, not the reciprocal width of a major arc.

\subsection{The asymptotic on the major arcs}

For $R\geq1$ and $H\geq2$, put
\begin{equation}\label{eq:major-family-count}
 {\mathscr{Z}}_{\mathrm{fam}}(\sigma;R,H)
 =\sum_{r\leq R}\sum_{\chi\bmod r}^{*}{\mathscr{Z}}(\sigma,H;\chi),
\end{equation}
where ${\mathscr{Z}}(\sigma,H;\chi)$ counts the nontrivial zeros, with
multiplicity, in $\Re\rho\geq\sigma$, $|\Im\rho|\leq H$.
The star restricts to primitive characters and includes the
primitive character of conductor one. For the major arcs we use
\eqref{dens:density-two}, the $A(\sigma)=20/9$ branch of
Heath--Brown's Theorem~2 \cite{HeathBrown1979Density}:
\begin{equation}\label{eq:major-hb}
 {\mathscr{Z}}_{\mathrm{fam}}(1-h;R,H)
 \ll_{\varepsilon}(R^2H^{6/5})^{(20/9)h+\varepsilon}
 \qquad (0\leq h\leq1/2).
\end{equation}
This bound is unconditional and uniform in the displayed
range. In particular, for fixed $\delta>0$ and $\nu>0$, taking
$\varepsilon=\nu\delta$ in \eqref{eq:major-hb} gives
\begin{equation}\label{eq:major-hb-away-one}
 {\mathscr{Z}}_{\mathrm{fam}}(1-h;R,H)
 \ll_{\delta,\nu}(R^2H^{6/5})^{(20/9+\nu)h}
 \qquad (\delta\leq h\leq1/2).
\end{equation}

For the two arithmetic signs, put
\[
 I_w^+(n)=\sum_{m\in\mathbb Z}w(m/X)w((n-m)/X),\qquad
 I_w^-(h)=\sum_{m\in\mathbb Z}w(m/X)w((m+h)/X).
\]

\begin{proposition}\label{prop:major-range}
Assume $(Q_\Theta)$ for a fixed $1/2\leq\Theta<1$.
For every fixed $0<\kappa<9/20$, there is $c>0$ such that,
uniformly for even $n\in[X/2,X]$,
\begin{equation*}
 \int_{\mathfrak M(P)}S_w(\alpha)^2e(-n\alpha)\,d\alpha
 =\mathfrak S(n)I_w^+(n)+O_{w,\Theta,\kappa}(X^{1-c}),
\end{equation*}
and, uniformly for even $h$ with $2\leq h\leq X/2$,
\begin{equation*}
 \int_{\mathfrak M(P)}|S_w(\alpha)|^2e(-h\alpha)\,d\alpha
 =\mathfrak S(h)I_w^-(h)+O_{w,\Theta,\kappa}(X^{1-c}),
\end{equation*}
where, for a nonzero even integer $r$,
\[
 \mathfrak S(r)=
 2\prod_{p>2}\left(1-\frac1{(p-1)^2}\right)
   \prod_{\substack{p\mid r\\p>2}}\frac{p-1}{p-2}.
\]
On the stated ranges,
$I_w^+(n),I_w^-(h)\geq X/4+O(1)$ and
$\mathfrak S(n),\mathfrak S(h)\gg1$. In particular, both
coefficients on the major arcs have real part $\gg_w X$ for sufficiently
large $X$.
\end{proposition}

\begin{proof}
Write $L=\log X$, $\Theta_0=\max(\Theta,7/8)$ and
$\delta=1-\Theta_0>0$, and put
\[
 T_w(\xi)=\sum_{m\in\mathbb Z}w(m/X)e(m\xi).
\]
For a primitive character $\chi\bmod r$, define
\begin{align*}
 \widetilde S_w(\xi,\chi)
   &=\sum_{n\geq1}\Lambda(n)\chi(n)w(n/X)e(n\xi),\\
 W_{X,\xi}(s)&=\int_0^\infty w(u/X)e(\xi u)u^{s-1}\,du,
\end{align*}
and put
\[
 E_\chi(\xi)=\widetilde S_w(\xi,\chi)
       -\mathbf1_{r=1}W_{X,\xi}(1).
\]
We first estimate the sum of their second moments on
$|\xi|\leq P/(rX)$.

Fix a sufficiently small $\eta>0$, to be specified below, and
set $H_r=X^\eta P/r$.  Mellin inversion gives
\[
 \widetilde S_w(\xi,\chi)=\frac1{2\pi i}
  \int_{(2)}-\frac{L'}{L}(s,\chi)W_{X,\xi}(s)\,ds.
\]
Moving the contour to the left gives the residue
$W_{X,\xi}(1)$ for $r=1$, and the residue
$-W_{X,\xi}(\rho)$ at each zero.  The Mellin transform decays
rapidly on vertical lines, so the contour can be taken through
heights avoiding zeros and the resulting zero sum converges.
The trivial zeros contribute $O_w(1)$ uniformly on the arcs:
the possible zero at $s=0$ for an even nonprincipal primitive
character contributes $W_{X,\xi}(0)=O_w(1)$, and the remaining
terms at negative integers form an absolutely convergent geometric
sum, since the support of $w(u/X)$ has $u\asymp X$.
For $r=1$ the first trivial zero is at $s=-2$.

If $|\gamma|$ exceeds a sufficiently large fixed multiple of
$1+X|\xi|$, repeated nonstationary integration in
\[
 W_{X,\xi}(\sigma+i\gamma)
 =X^{\sigma+i\gamma}
    \int w(y)e(\xi Xy)y^{\sigma+i\gamma-1}\,dy
\]
gives
\begin{equation*}
 W_{X,\xi}(\sigma+i\gamma)
   \ll_{A,w}X^\sigma(1+|\gamma|)^{-A}
 \qquad (0\leq\sigma\leq1).
\end{equation*}
The phase has first derivative
$2\pi\xi X+\gamma/y$ of magnitude $\gg|\gamma|$ on the
fixed support; its higher derivatives are $O_j(|\gamma|)$.
This proves the asserted uniformity also when $\xi$ varies
through the entire arc.  The local zero count for a primitive
character is $O(\log(2r(|t|+2)))$ per unit interval about $t$.
Consequently the zeros with $|\gamma|>H_r$ contribute
\[
 \ll_{A,w}XH_r^{1-A}\log(2rH_r)\ll_{M,w,\eta}X^{-M},
\]
for any fixed $M$, on choosing $A$ large enough after $\eta$.
Here $H_r\geq X^\eta$ and $X|\xi|\leq P/r$.
Thus, uniformly in $r\leq P$ and the stated $\xi$ range,
\begin{equation}\label{eq:major-smooth-explicit}
 E_\chi(\xi)
  =-\sum_{\substack{\rho\text{ nontrivial}\\|\Im\rho|\leq H_r}}
       W_{X,\xi}(\rho)+O_w(1)+O_{M,w,\eta}(X^{-M}).
\end{equation}
All sums continue to count multiplicity.

Let $\mathcal Z$ be any subcollection of the truncated zeros
of this one character, and suppose $\Re\rho\leq b$ throughout
$\mathcal Z$.  Plancherel on the full real $\xi$ axis gives
\begin{align}
 \int_{\mathbb R}
  \left|\sum_{\rho\in\mathcal Z}W_{X,\xi}(\rho)\right|^2d\xi
 &=\int_0^\infty |w(u/X)|^2
       \left|\sum_{\rho\in\mathcal Z}u^{\rho-1}\right|^2du
       \notag\\[-2pt]
 &\ll_w X^{2b-1}|\mathcal Z|\log(2rH_r).
       \label{eq:major-strip-plancherel}
\end{align}
To verify the inequality, the kernel for a pair
$\rho=\sigma+i\gamma$, $\rho'=\sigma'+i\gamma'$ is
\[
 X^{\sigma+\sigma'-1+i(\gamma-\gamma')}
   \int |w(y)|^2y^{\sigma+\sigma'-2+i(\gamma-\gamma')}\,dy
 \ll_{A,w}
   \frac{X^{\sigma+\sigma'-1}}
        {(1+|\gamma-\gamma'|)^A}.
\]
Integration by parts in $\log y$ proves this uniformly for
$0\leq\sigma,\sigma'\leq1$.  The local zero count implies
that each row of the decaying kernel has sum
$O(\log(2rH_r))$, proving
\eqref{eq:major-strip-plancherel}.  The extension to the full
real axis is applied to this finite zero sum; the bounded
error in the explicit formula is integrated on the original arc.

Keep all zeros with $\Re\rho<1/2$ in one group and partition
$[1/2,\Theta_0]$ into intervals of length at most $1/L$.
For an interval with left endpoint $\lambda$, its upper
endpoint is at most $\lambda+1/L$, so
$X^{2b-1}\ll X^{2\lambda-1}$ in
\eqref{eq:major-strip-plancherel}.  Cauchy's inequality over
these $O(L)$ groups, followed by dyadic division
$R\leq r<2R$ of the conductors, therefore yields
\begin{align}
 \mathcal E_0
 &:=\sum_{r\leq P}\sum_{\chi\bmod r}^{*}
       \int_{|\xi|\leq P/(rX)}|E_\chi(\xi)|^2\,d\xi
       \notag\\
 &\ll_{w,\eta}L^C\left(
      X^\eta P^2+
      \max_{\substack{1\leq R\leq P\\
                       \delta\leq h\leq1/2}}
       X^{1-2h}{\mathscr{Z}}_{\mathrm{fam}}
        \left(1-h;2R,2X^\eta P/R\right)\right).
       \label{eq:major-density-functional}
\end{align}
Indeed the total number of zeros in a conductor block, with
height at most $2X^\eta P/R$, is
$O(X^\eta PRL)\ll X^\eta P^2L$.  This handles the group
below $1/2$, whose factor $X^{2b-1}$ is at most one.
For the other groups use $h=1-\lambda\geq\delta$ and the
displayed family count.  The bounded error in
\eqref{eq:major-smooth-explicit} contributes only
$O(P^2/X)$ to $\mathcal E_0$.

Choose $\nu>0$ so that $A_0=20/9+\nu$ satisfies
$\kappa A_0<1$, and then choose $\eta>0$ sufficiently small
that
\[
 d:=2-A_0(2\kappa+6\eta/5)>0,
 \qquad 2\kappa+\eta<1.
\]
The substitution for conductor and height in
\eqref{eq:major-hb-away-one} is
\[
 (2R)^2(2X^\eta P/R)^{6/5}
       \ll R^{4/5}P^{6/5}X^{6\eta/5}
       \leq P^2X^{6\eta/5}.
\]
Thus every term in the maximum in
\eqref{eq:major-density-functional} is
\[
 \ll X^{1-2h}(P^2X^{6\eta/5})^{A_0h}
       =X^{1-dh}\leq X^{1-d\delta}.
\]
Since $X^\eta P^2=X^{2\kappa+\eta}$, we conclude that
$\mathcal E_0\ll_{w,\Theta,\kappa}X^{1-c_0}$ for some
$c_0>0$, after absorbing the fixed logarithmic factors.

We next pass from primitive characters to the additive
polynomial.  On the arc about $a/q$, let
\[
 V_w(a/q+\xi)=\frac{\mu(q)}{\varphi(q)}T_w(\xi).
\]
If $\chi\bmod q$ is induced by a primitive character
$\chi^*\bmod r$, the two weighted von Mangoldt sums differ
by $O_w(\log(2q))$.  Their difference is supported on
$p^j\asymp X$ with $p\mid q$; for each $p$ there are
$O_w(1+1/\log p)$ admissible exponents, which proves this
bound.  Character expansion on the units modulo $q$ gives
coefficients $\tau_q(\overline\chi)/\varphi(q)$.  In
particular, the part containing the primitive errors is
\[
 A_q(a,\xi)=\frac1{\varphi(q)}
     \sum_{\chi\bmod q}\tau_q(\overline\chi)\chi(a)
                 E_{\chi^*}(\xi),
\]
and orthogonality over the reduced numerators gives exactly
\begin{equation*}
 \sum_{a\bmod q}^{*}|A_q(a,\xi)|^2
   =\frac1{\varphi(q)}\sum_{\chi\bmod q}
        |\tau_q(\overline\chi)|^2|E_{\chi^*}(\xi)|^2.
\end{equation*}
The nonzero induced Gauss sums have
\[
 q=rs,\quad (r,s)=1,\quad s\text{ squarefree},
 \qquad |\tau_q(\overline\chi)|^2=r.
\]
Enlarging their integration intervals to
$|\xi|\leq P/(rX)$, the total weight of the characters induced by
one primitive character is at most
\begin{equation}\label{eq:major-descendant-weight}
 \sum_{\substack{q=rs\leq P\\(r,s)=1}}
       \frac{r\mu(s)^2}{\varphi(q)}
 \leq\frac r{\varphi(r)}\sum_{s\leq P/r}\frac1{\varphi(s)}
 \ll L^2.
\end{equation}
Here $\sum_{s\leq y}1/\varphi(s)\ll\log(2y)$ follows from
$n/\varphi(n)=\sum_{d\mid n}\mu(d)^2/\varphi(d)$ and the
convergence of $\sum_d\mu(d)^2/(d\varphi(d))$; the same
identity gives $r/\varphi(r)\ll\log(2r)$.

For clarity, the errors in this reduction have the following
total size.  The sum over $(n,q)>1$ is
$O_w(\log(2q))$, uniformly in $a$ and $\xi$, and each
difference from the Euler factors introduced by induction obeys the same bound.
The latter may be squared after character orthogonality,
using
$\sum_{\chi\bmod q}|\tau_q(\chi)|^2=\varphi(q)^2$.
Their integrated contribution is $O_w(P^2L^2/X)$, since
the measure of the major arcs is $O(P^2/X)$.
Removing the prime powers changes the polynomial by
$O_w(X^{1/2}L^2)$, and hence contributes $O_w(P^2L^4)$.
Finally Poisson summation gives
\[
 T_w(\xi)-W_{X,\xi}(1)\ll_{M,w}X^{-M}
 \qquad (|\xi|\leq P/X),
\]
because $P/X=o(1)$ and all the nonzero dual frequencies are
bounded away from zero.  The only pole term in the character
expansion is therefore exactly the model $V_w$, up to this
negligible error.  Combining orthogonality,
\eqref{eq:major-descendant-weight} and these estimates gives
\begin{equation}\label{eq:major-error-density}
 \int_{\mathfrak M(P)}|S_w(\alpha)-V_w(\alpha)|^2\,d\alpha
  \ll_w L^C\mathcal E_0+P^2L^C
  \ll_{w,\Theta,\kappa}X^{1-c_1}
\end{equation}
for some $c_1>0$.

It remains to evaluate the two quadratic models. Parseval and
$\sum_{q\leq P}1/\varphi(q)\ll L$ give
\[
 \int_{\mathfrak M(P)}|V_w(\alpha)|^2\,d\alpha
 =\sum_{q\leq P}\frac{\mu(q)^2}{\varphi(q)}
   \int_{|\xi|\leq P/(qX)}|T_w(\xi)|^2\,d\xi
 \ll_w XL.
\]
Put $E=S_w-V_w$. The identities
\[
 S_w^2-V_w^2=2V_wE+E^2,\qquad
 |S_w|^2-|V_w|^2=V_w\overline E+\overline{V_w}E+|E|^2
\]
and Cauchy's inequality imply that, for every integer $r$, either
quadratic replacement costs at most
\begin{equation}\label{eq:major-two-sign-replacement}
 2\|V_w\|_2\|E\|_2+\|E\|_2^2
 \ll_{w,\Theta,\kappa}X^{1-c_2},
\end{equation}
for some $c_2>0$, by \eqref{eq:major-error-density}. Here the norms
are taken on $\mathfrak M(P)$; the bound is independent of $r$.

Write $c_q(r)=\sum_{a\bmod q}^{*}e(ar/q)$ for the Ramanujan sum.
It is real and satisfies $c_q(-r)=c_q(r)$. Since
$\mu(q)/\varphi(q)$ is real, the two model coefficients are
\begin{equation}\label{eq:major-two-models}
 \sum_{q\leq P}\frac{\mu(q)^2c_q(r)}{\varphi(q)^2}
  \int_{|\xi|\leq P/(qX)}G_\pm(\xi)e(-r\xi)\,d\xi,
 \qquad
 G_+(\xi)=T_w(\xi)^2,\quad G_-(\xi)=|T_w(\xi)|^2.
\end{equation}
The difference model therefore retains the conjugation in its
archimedean factor.

For every fixed $B>1$, Poisson summation and smoothness give
\[
 \int_{P/(qX)<|\xi|\leq1/2}|T_w(\xi)|^2\,d\xi
 \ll_{B,w}X(q/P)^B\qquad(q\leq P).
\]
For every $\varepsilon>0$, uniformly in integers $0<|r|\leq X$,
\begin{align}
 \sum_{q>P}\frac{|c_q(r)|}{\varphi(q)^2}
 &\ll_\varepsilon X^\varepsilon/P,
 \label{eq:major-ramanujan-tail-both}\\
 \sum_{q\leq P}\frac{|c_q(r)|}{\varphi(q)^2}(q/P)^B
 &\ll_{B,\varepsilon}X^\varepsilon/P.\notag
\end{align}
Indeed, $|c_q(r)|\leq(q,r)\leq\sum_{d\mid(q,r)}d$ and
$\varphi(q)^{-2}\ll_\eta q^{-2+\eta}$. Writing $q=dk$ gives
\[
 \sum_{q>P}\frac{(q,r)}{q^{2-\eta}}
 \ll_\eta \tau(|r|)P^{-1+\eta},\qquad
 \sum_{q\leq P}(q,r)q^{B-2+\eta}
 \ll_{B,\eta}\tau(|r|)P^{B-1+\eta}.
\]
For the first estimate, split $d\mid r$ at $P$; for the second,
use
$\sum_{k\leq P/d}k^{B-2+\eta}\ll_{B,\eta}(P/d)^{B-1+\eta}$.
Choose $\eta>0$ sufficiently small in terms of $\varepsilon$ and
use $\tau(|r|)\ll_\varepsilon |r|^{\varepsilon/2}$.
The estimates are uniform for small nonzero $r$ as well as for
$r\asymp X$.

Since $|G_+|=|G_-|=|T_w|^2$, completing either integral in
\eqref{eq:major-two-models} to $[-1/2,1/2]$ costs\linebreak
$O_{w,\varepsilon}(X^{1+\varepsilon}/P)$. Fourier orthogonality
gives the exact identities
\begin{align*}
 \int_{-1/2}^{1/2}T_w(\xi)^2e(-n\xi)\,d\xi
 &=I_w^+(n),\\
 \int_{-1/2}^{1/2}|T_w(\xi)|^2e(-h\xi)\,d\xi
 &=I_w^-(h).
\end{align*}
In the second identity, the exponent is $e((m_2-m_1-h)\xi)$,
so that $m_2=m_1+h$. Both completed integrals are $O_w(X)$.
Equation~\eqref{eq:major-ramanujan-tail-both} thus completes the
$q$ sum with the same error. For $r\ne0$ the completed series is
absolutely convergent and, for even $r$, equals
\[
 \sum_{q\geq1}\frac{\mu(q)^2c_q(r)}{\varphi(q)^2}
 =\prod_p\left(1+\frac{c_p(r)}{(p-1)^2}\right)
 =\mathfrak S(r).
\]
Here $c_p(r)=p-1$ when $p\mid r$ and $c_p(r)=-1$ otherwise.
The zero shift is excluded throughout. Taking $0<\varepsilon<\kappa$
and combining the completion bounds with
\eqref{eq:major-two-sign-replacement} proves both asymptotics,
after decreasing $c$.

For the sum sign, both weights equal one for
$X/8\leq m\leq n-X/8$, an interval of length
$n-X/4\geq X/4$. For the difference sign, use
$X/8\leq m\leq7X/8-h$, whose length is
$3X/4-h\geq X/4$. Counting integers in these intervals proves
both lower bounds. The first product in $\mathfrak S(r)$ is
positive and each factor in the second product is at least one.
\end{proof}

The corresponding continuous correlations are
\[
 J_w^+(n)=X\int_{\mathbb R}w(t)w(n/X-t)\,dt,\qquad
 J_w^-(h)=X\int_{\mathbb R}w(t)w(t+h/X)\,dt.
\]
For every fixed $A>0$, Poisson summation gives
\begin{equation*}\label{eq:discrete-continuous-correlations}
 I_w^+(n)=J_w^+(n)+O_{A,w}(X^{-A}),\qquad
 I_w^-(h)=J_w^-(h)+O_{A,w}(X^{-A}),
\end{equation*}
uniformly on the respective ranges. To see this, the functions
$w(t)w(n/X-t)$ and $w(t)w(t+h/X)$ have compact support and
uniformly bounded derivatives of every fixed order. Their nonzero
Poisson frequencies therefore contribute $O_{{\mathsf{M}},w}(X^{1-{\mathsf{M}}})$,
and ${\mathsf{M}}$ may be chosen larger than $A+1$. Equivalently, if
\[
 W(\xi)=\int_{\mathbb R}w(u/X)e(\xi u)\,du,
\]
Fourier inversion gives
\[
 \int_{\mathbb R}|W(\xi)|^2e(-h\xi)\,d\xi
 =\int_{\mathbb R}w(u/X)w((u+h)/X)\,du
 =J_w^-(h).
\]
Thus the discrete and continuous conventions have the same leading
term, with no additional normalization factor.

\subsection{Recombine the principal character before factorization}

For a fixed modulus, character orthogonality is the finite analogue of
orthogonality of exponential functions. For integers coprime to the
modulus, it says
\begin{equation*}
 \frac1{\varphi(q)}\sum_{\chi\bmod q}\chi(n)\overline{\chi(m)}
 =\mathbf1_{n\equiv m\pmod q}.
\end{equation*}
The principal summand has no multiplicative cancellation. Splitting the
prime polynomial first and treating each principal factor separately
would obscure this complete contribution. We identify and estimate it
before making that splitting.

At each point on the minor arcs choose a Dirichlet approximation
\[
 \alpha=a/q+\xi,\qquad (a,q)=1,\qquad
 q\leq X^{1/2},\qquad |\xi|\leq(qX^{1/2})^{-1},
\]
using a measurable choice where intervals overlap.  The principal
character modulo \(q\) is induced by the character of conductor
one, and its Gauss sum is \(\mu(q)\).  Its complete contribution
to the prime sum is consequently
\begin{equation*}\label{eq:principal-complete-term}
 S_{\rm pr}(a/q+\xi)
 =\frac{\mu(q)}{\varphi(q)}S_w(\xi)+O(\log(2q)).
\end{equation*}
The error accounts for primes dividing \(q\).  If the Gauss
expansion is performed for the von Mangoldt sum, the further
difference contributed by prime powers is \(O_w(X^{1/2}\log^2X)\) pointwise.
On the chosen cells these errors have total fourth integral
\(O_w(X^2(\log X)^8)\), which will be absorbed below.

After subtracting this complete principal term, the Heath--Brown
decomposition involves only characters induced from primitive
nonprincipal characters. The Mellin shifts for the long factors
therefore cross no principal pole.

For the principal character,
\begin{equation*}
 L(s,\chi_{0,q})=\zeta(s)\prod_{p\mid q}(1-p^{-s}),
\end{equation*}
so the pole of the zeta function at one survives. For a nonprincipal
character the $L$-function is holomorphic there, even when the character
is induced. This is the distinction needed when a contour is moved to
the critical line in the argument using fourth moments below.

\subsection{The principal term from individual zero strips}

We estimate the complete principal term by grouping the zeros of
$\zeta$ according to their real parts.  On each group, a pointwise
estimate by stationary phase and a local second moment use the same zero
count.  This preserves the dependence on the real part until the
fourth moment has been formed.

Put
\[
 \widetilde S_w(\xi)=\sum_{n\geq1}\Lambda(n)w(n/X)e(n\xi),
 \qquad
 W_{X,\xi}(s)=\int_0^\infty w(u/X)e(\xi u)u^{s-1}\,du.
\]
For a nontrivial zeta zero write $\rho=b_\rho+i\gamma_\rho$.
All zero sums below count multiplicity.  We use the two unconditional
bounds
\begin{equation}\label{eq:principal-two-densities}
 {\mathscr{Z}}_\zeta(\lambda,H)\ll_\varepsilon
 \begin{cases}
 H^{3(1-\lambda)/(2-\lambda)+\varepsilon},
       &\frac12\leq\lambda\leq\frac{11}{14},\\
 H^{2(1-\lambda)+\varepsilon},
       &\frac{11}{14}\leq\lambda\leq1.
 \end{cases}
\end{equation}
Here $H\geq2$ and ${\mathscr{Z}}_\zeta(\lambda,H)$ counts the zeros with
$b_\rho\geq\lambda$ and $|\gamma_\rho|\leq H$.
The first bound is Ingham's estimate, recorded in
Proposition~\ref{dens:lower-zeta-density}. The second is the
specialization to conductor one of \eqref{dens:density-one}, from
Heath--Brown's Theorem~1 \cite{HeathBrown1979Density}.
Both estimates are uniform in the displayed ranges of real parts.

\begin{lemma}\label{lem:principal-strip-fourth}
Fix $\eta>0$ and let $X^\eta\leq K\leq X^{1/2}$.
Let $\mathcal Z$ be any collection of zeta zeros with
$|\gamma_\rho|\leq C K$ and $b_\rho\leq b$, where $0\leq b\leq1$
and $C$ is fixed.  Set
\[
 Z(\xi)=\sum_{\rho\in\mathcal Z}W_{X,\xi}(\rho),
 \qquad R=|\mathcal Z|.
\]
For any fixed positive shell endpoints one has
\begin{align}
 \sup_{|\xi|\asymp K/X}|Z(\xi)|^2
   &\ll_{w,C}X^{2b}K^{-1}R^2,
       \label{eq:principal-strip-sup}\\
 \int_{|\xi|\asymp K/X}|Z(\xi)|^2\,d\xi
   &\ll_{w,C}X^{2b-1}R\log(2K),
       \label{eq:principal-strip-second}\\
 \int_{|\xi|\asymp K/X}|Z(\xi)|^4\,d\xi
   &\ll_{w,C}X^{4b-1}K^{-1}R^3\log(2K).
       \label{eq:principal-strip-fourth}
\end{align}
The constants may depend on the fixed shell endpoints.
\end{lemma}

\begin{proof}
It suffices to prove the estimates on a positive shell; replacing
$\xi$ by $-\xi$ gives the same argument with conjugate zeros.
Choose a fixed smooth nonnegative function $\Phi$, supported in a
compact subinterval of $(0,\infty)$, which majorizes the rescaled
positive shell.  Write $\xi=Kz/X$ and substitute $u=Xy/z$.  Then
\begin{equation*}
 W_{X,Kz/X}(b_\rho+i\gamma_\rho)
  =X^{b_\rho+i\gamma_\rho}z^{-i\gamma_\rho}
      J_{b_\rho,\gamma_\rho}(z),
\end{equation*}
where
\[
 J_{v,\gamma}(z)
  =z^{-v}\int_0^\infty
          w(y/z)e(Ky)y^{v+i\gamma-1}\,dy.
\]
For every fixed nonnegative integer $j$,
\begin{equation}\label{eq:principal-strip-J}
 \frac{d^j}{dz^j}J_{v,\gamma}(z)\ll_{j,w,C}K^{-1/2}
 \qquad(0\leq v\leq1,\ |\gamma|\leq CK).
\end{equation}
Indeed, all derivatives with respect to $z$ fall on a smooth amplitude supported
in a fixed compact subinterval of $(0,\infty)$.  The phase
$2\pi Ky+\gamma\log y$ is independent of $z$.  If $|\gamma|$
is smaller than a sufficiently small fixed multiple of $K$, its
first derivative has magnitude $\gg K$.  Otherwise its second
derivative has magnitude $\gg K$.  The corresponding estimate based on the first
or second derivative gives \eqref{eq:principal-strip-J}.
In particular each zero contributes $O(X^bK^{-1/2})$, proving
\eqref{eq:principal-strip-sup}.

For two zeros, integration by parts in $z$ against
$z^{-i(\gamma_\rho-\gamma_{\rho'})}$, using
\eqref{eq:principal-strip-J}, gives for every fixed $A>0$
\begin{equation}\label{eq:principal-strip-pair-kernel}
 \left|\int \Phi(\xi X/K)
       W_{X,\xi}(\rho)\overline{W_{X,\xi}(\rho')}\,d\xi\right|
 \ll_{A,w,C}
 \frac{X^{b_\rho+b_{\rho'}-1}}
      {(1+|\gamma_\rho-\gamma_{\rho'}|)^A}.
\end{equation}
The Jacobian $K/X$ cancels the factor $K^{-1}$ from the two
amplitudes.  The standard local zero count, with multiplicity,
is $O(\log(2K))$ per unit ordinate interval within $[-CK,CK]$.
Consequently each row of the decaying kernel in
\eqref{eq:principal-strip-pair-kernel} has sum
$O(\log(2K))$.  Expanding $|Z|^2$ and using
$b_\rho+b_{\rho'}\leq2b$ proves
\eqref{eq:principal-strip-second}.  Multiplying this bound by
\eqref{eq:principal-strip-sup} proves
\eqref{eq:principal-strip-fourth}.
\end{proof}

\begin{proposition}\label{prop:principal-lower-theta}
Assume $(Q_\Theta)$ with $\Theta<11/12$, and let
$4/9<\kappa<9/20$.  On the chosen cells on the minor arcs the complete
principal contribution satisfies, for every $\varepsilon>0$,
\begin{equation}\label{eq:principal-lower-theta-fourth}
 \int_{\mathfrak m(P)}|S_{\rm pr}(\alpha)|^4\,d\alpha
       \ll_{w,\Theta,\kappa,\varepsilon}X^{22/9+\varepsilon}.
\end{equation}
This holds for the principal term of either the prime or the
von Mangoldt expansion.
\end{proposition}

\begin{proof}
Set $\Theta_0=\max(\Theta,7/8)$ and $\eta=1/100$.
For a chosen approximation $\alpha=a/q+\xi$, let
$K=\max(1,X|\xi|)$.  The conditions for the minor arcs and Dirichlet approximation give
\[
  P\ll qK\ll X^{1/2}.
\]
We may enlarge the chosen cells to dyadic denominator blocks
$q\asymp Q$ and additive shells $|\xi|\asymp K/X$, treating
$K\asymp1$ as $|\xi|\ll X^{-1}$.
The sum of fourth powers of the principal coefficients is
\begin{equation}\label{eq:principal-gauss-weight}
 \sum_{q\asymp Q}\sum_{\substack{a\bmod q\\(a,q)=1}}
      \left|\frac{\mu(q)}{\varphi(q)}\right|^4
 =\sum_{q\asymp Q}\frac{\mu(q)^4}{\varphi(q)^3}
 \ll_\varepsilon Q^{-2+\varepsilon}.
\end{equation}
Here and below dyadic boxes are retained only if
$P\ll QK\ll X^{1/2}$; fixed factors do not affect the estimates.

When $K<X^\eta$, the trivial bound
$|\widetilde S_w(\xi)|\ll_w X\log X$ and the shell length
$O(K/X)$ give a contribution
\[
 \ll X^{3+\varepsilon}Q^{-2}K
 \ll X^{3-2\kappa+3\eta+\varepsilon}
 \ll X^{22/9+\varepsilon}.
\]
The last inequality has a fixed power margin, since
$\kappa>4/9$ and $\eta=1/100$.

Suppose now $K\geq X^\eta$.  The smooth explicit formula is
\begin{equation*}
 \widetilde S_w(\xi)
 =W_{X,\xi}(1)-\sum_\rho W_{X,\xi}(\rho)+O_w(X^{-2}).
\end{equation*}
The sum over trivial zeros is included in the error.  On the present
shell the pole term is $O_{A,w}(XK^{-A})$.  If $C$ is sufficiently
large, the zero terms with $|\gamma_\rho|>CK$ have no stationary
point, and repeated integration by parts gives
\[
 W_{X,\xi}(\rho)\ll_{A,w}X^{b_\rho}|\gamma_\rho|^{-A}.
\]
The ordinary zero count bounds their total absolute contribution by
\[
 O_{A,w}\bigl(XK^{1-A}\log(2K)\bigr).
\]
As $K\geq X^\eta$, choosing the fixed integer $A$ sufficiently
large makes both contributions smaller than any prescribed negative
power of $X$.  Thus it remains to estimate the zeros of height at
most $CK$.

Put $L=\log X$.  Treat the zeros with $b_\rho<1/2$ as one
group.  By Lemma~\ref{lem:principal-strip-fourth}, with
$b=1/2$ and $R\ll K\log(2K)$, their fourth integral is
$\ll X^{1+\varepsilon}K^2$.  After multiplying by
\eqref{eq:principal-gauss-weight} and using
$K\ll X^{1/2}/Q$, this is $\ll X^{2+\varepsilon}Q^{-4}$.

Partition the remaining interval $[1/2,\Theta_0]$ into intervals
of length at most $L^{-1}$, putting a partition endpoint at
$11/14$ and at $5/6$.  Let $\lambda$ be the left endpoint of one
such interval, and let $Z_\lambda$ be its zero sum.  Its zeros
satisfy $b_\rho\leq\lambda+L^{-1}$, so
$X^{4b_\rho}\leq e^4X^{4\lambda}$.  Its cardinality is at most
${\mathscr{Z}}_\zeta(\lambda,CK)$.  Equations
\eqref{eq:principal-strip-fourth} and
\eqref{eq:principal-two-densities} give
\begin{equation}\label{eq:principal-strip-box}
 Q^{-2}\int_{|\xi|\asymp K/X}|Z_\lambda(\xi)|^4\,d\xi
 \ll X^\varepsilon Q^{-2}X^{4\lambda-1}
          K^{3g(\lambda)-1},
\end{equation}
where we may take
\[
 g(\lambda)=
 \begin{cases}
 \dfrac{42}{17}(1-\lambda),
          &\dfrac12\leq\lambda<\dfrac{11}{14},\\[1mm]
 2(1-\lambda),&\dfrac{11}{14}\leq\lambda\leq\Theta_0.
 \end{cases}
\]
For the first choice we used
$3/(2-\lambda)\leq42/17$.
There are $O(L)$ groups.  The inequality
$|\sum_{j=1}^J z_j|^4\leq J^3\sum_j|z_j|^4$ shows that their
recombination costs only a fixed power of $L$, already absorbed
in the arbitrary $X^\varepsilon$ loss.

If $1/2\leq\lambda<11/14$, then
$3g(\lambda)-1>0$.  Taking $K\ll X^{1/2}/Q$ in
\eqref{eq:principal-strip-box} gives
\begin{align*}
 Q^{-2}X^{4\lambda-1}K^{3g(\lambda)-1}
 &\ll X^{4\lambda-3/2+(63/17)(1-\lambda)}
                 Q^{-1-3g(\lambda)}\\
 &\leq X^{290/119}Q^{-1-3g(\lambda)}.
\end{align*}
The exponent of $X$ is increasing in $\lambda$, with derivative
$5/17$, and its value at $11/14$ is $290/119<22/9$.
The denominator powers are negative and summable over dyadic $Q$.

For $11/14\leq\lambda\leq5/6$ one has
$3g(\lambda)-1=5-6\lambda\geq0$.  The upper bound on $K$
now gives
\[
 Q^{-2}X^{4\lambda-1}K^{5-6\lambda}
 \ll X^{\lambda+3/2}Q^{6\lambda-7}
 \leq X^{7/3}Q^{6\lambda-7}.
\]
These denominator powers are also summable.  Finally, if
$5/6\leq\lambda\leq\Theta_0$, the exponent $5-6\lambda$
is nonpositive.  Using $K\gg P/Q$ gives instead
\[
 Q^{-2}X^{4\lambda-1}K^{5-6\lambda}
 \ll X^{4\lambda-1+\kappa(5-6\lambda)}Q^{6\lambda-7}.
\]
This remains a valid upper bound when $P/Q<1$.
Since $4-6\kappa>0$, its exponent of $X$ is at most
\[
 4\Theta_0-1+\kappa(5-6\Theta_0)
 \leq\frac83-\frac\kappa2<\frac{22}{9}.
\]
Moreover $6\lambda-7\leq-3/2$, so the dyadic sum over $Q$
converges uniformly.  Summing the $O(L)$ additive shells and all
zero strips proves \eqref{eq:principal-lower-theta-fourth} for
the complete von Mangoldt principal term.

The correction from Euler factors is $O_w(\log(2q))$ on the chosen
cells.  Also
$\widetilde S_w(\xi)-S_w(\xi)=O_w(X^{1/2}\log^2X)$.
Their fourth integrals on a set of measure at most one are
$O_w(X^2(\log X)^8)$, and the principal coefficient has magnitude
at most one.  These errors are absorbed by
\eqref{eq:principal-lower-theta-fourth}, proving the assertion
also for the prime expansion.
\end{proof}

After division by $X^2$ in Bessel's inequality, the principal
contribution is $O_\varepsilon(X^{4/9+\varepsilon})$.
In particular it has a fixed power margin below the desired
square root bound throughout the stated range of $\Theta$.

\section{An estimate for difference energy with a fixed logarithmic loss}
\label{sec:difference-energy-fixed-log}

The next estimates concern a finite Dirichlet polynomial, a sum of
multiplicative characters multiplied by powers $n^{it}$. A large value
at one point need not be rare, but many large values constrain how the
points can add and subtract. The energy counts quadruples whose
characters multiply in the same way and whose ordinates have nearly
the same sum. Both conditions will arise directly when a fourth power
is expanded.

We use Dirichlet convolution, defined by
\begin{equation*}
 (a*b)(\ell)=\sum_{mn=\ell}a_m b_n.
\end{equation*}
It is the coefficient operation corresponding to multiplication of
Dirichlet polynomials. The square norm of this convolution counts
weighted solutions of a multiplicative equation. Keeping that equation
intact, rather than estimating each common divisor separately, will
control the logarithmic loss.

Different common divisors may give the same reduced ratio. Combining
their coefficients before taking the mean square retains the joint
convolution norm needed at the terminal scale. The required
mean value at a fixed modulus is as follows.

Constants in this section depend only on the fixed support ratio $c>1$
and integers $\nu\geq1$, $A\geq0$. Write $d_j$ for the coefficients of
$\zeta(s)^j$. Character sums include all characters of one fixed modulus.

\begin{lemma}
\label{lem:de-fixed-modulus-mean}
Let $q\geq1$, $U\geq1$, $M>0$, and
\[
 D(t,\chi)=\sum_{M\leq n\leq cM} b_n\chi(n)n^{it}.
\]
For any real centre $t_0$,
\begin{equation}\label{eq:de-continuous-mean}
 \sum_{\chi\bmod q}\int_{t_0-U}^{t_0+U}|D(t,\chi)|^2\,dt
 \ll_c \varphi(q)\left(U+\frac Mq\right)
       \sum_{\substack{M\leq n\leq cM\\(n,q)=1}}|b_n|^2.
\end{equation}
Moreover, if $\mathcal J$ is the set of integer bins $[j,j+1)$
intersecting $[t_0-U,t_0+U]$, then
\begin{equation}\label{eq:de-bin-mean}
 \sum_{\chi\bmod q}\sum_{j\in\mathcal J}
       \sup_{t\in[j,j+1]}|D(t,\chi)|^2
 \ll_c \varphi(q)\left(U+\frac Mq\right)
       \sum_{\substack{M\leq n\leq cM\\(n,q)=1}}|b_n|^2.
\end{equation}
Thus the same bound holds for any collection of points which is
separated by at least $1$ for each character.  No logarithmic loss is needed
for the sampling step.
\end{lemma}

\begin{proof}
Twisting $b_n$ by $n^{it_0}$ reduces the integral to centre zero.

The twist has absolute value one and leaves the coefficient square
norm unchanged. Expanding the square then gives one term for every
ordered pair of indices. Character orthogonality removes all pairs
except those congruent modulo the modulus. The Fourier transform of
the Gaussian supplies a nonnegative rapidly decaying weight in the
logarithm of their ratio.
Majorizing its characteristic function by a constant times
$\exp(-t^2/U^2)$ and applying character orthogonality, we bound
the resulting quadratic form in absolute value by
\[
 C\varphi(q)U
 \sum_{\substack{n,m\asymp_c M\\n\equiv m\pmod q\\(nm,q)=1}}
 |b_n b_m|\exp\{-c_1U^2\log^2(n/m)\}.
\]
Since $|\log(n/m)|\gg_c|n-m|/M$, each row of the kernel has sum
$O_c(1+M/(qU))$.

For a fixed index, the congruent indices differ from it by integer
multiples of the modulus. Hence the row sum is bounded by
\begin{equation*}
 \sum_{k\in\mathbb Z}\exp\{-c_2(Uqk/M)^2\}
 \ll 1+\frac{M}{qU}.
\end{equation*}
The diagonal gives the first term. Comparison with a Gaussian integral
gives the second. Multiplying by the outer factor gives precisely the
two scales in the asserted mean square.

Applying $2|b_nb_m|\leq|b_n|^2+|b_m|^2$ proves
\eqref{eq:de-continuous-mean}, also for $M<1$.

Apply Sobolev on $[j-1,j+2]$ to $M^{-it}D(t,\chi)$.
Differentiation inserts the bounded factor $\log(n/M)$.

For completeness, if a continuously differentiable function $f$ is
considered on an interval of length three, the fundamental theorem
of calculus applied to $|f|^2$ gives
\begin{equation*}
 \sup |f|^2\ll \int \bigl(|f|^2+|f f'|\bigr).
\end{equation*}
Cauchy's inequality bounds the last integral by the product of the two
square norms. Here we apply this to the polynomial after removing the
phase of its lower length scale, so that
\begin{equation*}
 \frac{d}{dt}\{M^{-it}D(t,\chi)\}
 =i\sum_n b_n\chi(n)(n/M)^{it}\log(n/M).
\end{equation*}
The new coefficients have the same square norm bound up to a constant
because the support ratio is fixed.
The enlarged intervals have bounded overlap, so
\eqref{eq:de-continuous-mean} at height $U+2$ gives
\eqref{eq:de-bin-mean}. Unit separation gives bounded multiplicity per bin.
\end{proof}

The energy below is ordered: changing the order of a quadruple can
produce a different counted quadruple. For a nonempty separated set of
size $R$, it satisfies the elementary bounds
\begin{equation*}
 R^2\le\mathcal E_1(W)\le3R^3.
\end{equation*}
The lower bound uses the quadruples with first and third points equal
and second and fourth points equal. For the upper bound, the first
three points fix the fourth character and restrict its ordinate to an
interval of length two, containing at most three separated points.
These bounds also explain why division by the square root of the
energy in the proof is legitimate.

\begin{lemma}
\label{lem:difference-energy-fixed-log}
Let $q\geq1$, $T\geq1$, $B=qT$, and ${\mathsf{M}}\geq1$.  Let $W$ be a
finite nonempty set of pairs $w=(t,\chi)$, where $|t|\leq T$ and
$\chi$ is an arbitrary character modulo $q$.  Assume that the
ordinates belonging to any one character are separated by at least $1$.
Write $R=|W|$ and
\begin{equation}\label{eq:de-energy-definition}
 \mathcal E_1(W)=\#\bigl\{(w_1,w_2,w_3,w_4)\in W^4:
 \chi_1\chi_2=\chi_3\chi_4,\quad
 |t_1+t_2-t_3-t_4|\leq1\bigr\}.
\end{equation}
Suppose
\[
 F(t,\chi)=\sum_{{\mathsf{M}}\leq n\leq c{\mathsf{M}}}a_n\chi(n)n^{it},\qquad
 |F(t,\chi)|\geq V>0\quad((t,\chi)\in W),
\]
and put
\[
 A_\infty=\max_n|a_n|,\qquad
 \mathcal C(a)=\sum_{\ell\geq1}
                    \bigl((|a|*|a|)(\ell)\bigr)^2.
\]
Then
\begin{equation}\label{eq:energy-coefficient-norm}
 \mathcal E_1(W)V^4\ll_c
 \mathbf1_{{\mathsf{M}}^2\geq B}A_\infty^4R{\mathsf{M}}^4+BR^2\mathcal C(a).
\end{equation}
The first term is absent when \({\mathsf{M}}<\sqrt B\).
The threshold \(V\) is arbitrary, and the characters may be induced
or principal. A different fixed positive energy tolerance changes
only the implied constant.
\end{lemma}

\begin{proof}
Write $\mathcal E=\mathcal E_1(W)$, and set
\[
 \mathcal R(y,a)=\sum_{(t,\chi)\in W}e^{ity}\chi(a),
 \qquad a\in(\mathbb Z/q\mathbb Z)^\times.
\]
The frequencies of ${\mathsf{M}}^{-it}F(t,\chi)$ lie in $[0,\log c]$.

Indeed, each summand of this normalized polynomial has oscillating
factor $e^{it\log(n/{\mathsf{M}})}$. The normalization changes no absolute value.
A reproducing kernel is a function whose convolution leaves all such
frequencies unchanged. It allows a value at one ordinate to be bounded
by a weighted average near another ordinate. In an energy quadruple,
the fourth ordinate is within a fixed distance of the sum of the first
two minus the third, which is exactly the situation required here.
Choose a Schwartz kernel whose Fourier transform is one there.
Applying its reproducing formula and Cauchy's inequality, we obtain
\[
 |F(t_4,\chi_4)|^2
 \ll_c\int_{\mathbb R}h(s)
  |F(t_1+t_2-t_3+s,\chi_1\chi_2\overline\chi_3)|^2\,ds
\]
for every energy quadruple. Here $h$ is a fixed nonnegative,
rapidly decaying envelope of kernel translates by at most one.
Each triple permits at most three fourth points. Summing over triples,
expanding, and using $|\widehat h|\leq\|h\|_1$ therefore gives
\begin{equation}\label{eq:de-cubic}
 V^2\mathcal E\ll_c
 \sum_{\substack{{\mathsf{M}}\leq n_1,n_2\leq c{\mathsf{M}}\\(n_1n_2,q)=1}}
 |a_{n_1}a_{n_2}|
 \left|\mathcal R\left(\log(n_1/n_2),n_1n_2^{-1}\right)\right|^3.
\end{equation}

Put $d_0={\mathsf{M}}^2/B$. For $d=(n_1,n_2)\leq d_0$ and a fixed residue $a$,
the reduced ratios have cross differences divisible by $q$ and
denominators at most $c{\mathsf{M}}/d$. Their logarithms are therefore
separated by $\gg_c qd^2/{\mathsf{M}}^2$ in a fixed compact interval.
This range is empty if $d_0<1$.

To verify the separation, write the two reduced fractions as $m/n$
and $m'/n'$. Equality of their residues means that
$mn'-m'n$ is divisible by $q$. If the fractions differ, this integer
is nonzero, so
\begin{equation*}
 \left|\frac mn-\frac{m'}{n'}\right|
 =\frac{|mn'-m'n|}{nn'}\ge\frac{q}{nn'}
 \gg_c\frac{qd^2}{{\mathsf{M}}^2}.
\end{equation*}
Their values lie in a fixed compact interval away from zero, where
passing to logarithms changes distances by bounded factors. The cutoff
for the common divisor is chosen so that
\begin{equation*}
 B d_0={\mathsf{M}}^2.
\end{equation*}
This balances the two terms obtained after summing the separated
sampling estimate over the small common divisors.

The frequencies of $\mathcal R(y,a)$ lie in $[-T,T]$.
Its exact Schwartz reproducing formula and H\"older's inequality
show that a separated sum over these logarithms is bounded by
\[
 \ll_c\left(T+\frac{{\mathsf{M}}^2}{qd^2}\right)
 \int_{\mathbb R}(1+y^2)^{-10}|\mathcal R(y,a)|^3\,dy.
\]
Indeed the kernel is $O_A(T(1+T|x-y|)^{-A})$.
Summing it over the separated points gives the asserted bound on
a compact interval. Outside that interval, rapid decay supplies the
weight, so the kernel remains integrated over the whole real line.

Orthogonality and separation give, with $\omega(y)=(1+y^2)^{-10}$,
\begin{equation}\label{eq:de-R-moments}
 \sum_a\int\omega(y)|\mathcal R(y,a)|^2\,dy\ll qR,
 \qquad
 \sum_a\int\omega(y)|\mathcal R(y,a)|^4\,dy\ll q\mathcal E.
\end{equation}
For the second bound, bin pair sums for each product character,
with populations $p_{\psi,j}$. Rapid decay of $\widehat\omega$ gives
\[
 \varphi(q)\sum_\psi\sum_{j,k}
 p_{\psi,j}p_{\psi,k}(1+|j-k|)^{-3}
 \ll q\sum_{\psi,j}p_{\psi,j}^2\leq q\mathcal E.
\]
The first bound follows likewise from unit separation.

The bin estimate uses only the summability of the discrepancy kernel.
For every integer displacement, Cauchy's inequality gives
\begin{equation*}
 \sum_j p_{\psi,j}p_{\psi,j+h}
 \le\sum_jp_{\psi,j}^2.
\end{equation*}
Summing this with weights $(1+|h|)^{-3}$ bounds the double bin sum.
Two pair sums in one unit bin differ by at most one, so their squared
population is already counted in the energy. Similarly, for the second
moment each character has only boundedly many ordinates in a unit
interval. These are the two uses of separation and energy in the
moment bounds.
The intermediate third moment is controlled by the second and fourth:
\begin{equation*}
 \sum_a\int\omega(y)|\mathcal R(y,a)|^3\,dy
 \le\left(\sum_a\int\omega|\mathcal R|^2\right)^{1/2}
      \left(\sum_a\int\omega|\mathcal R|^4\right)^{1/2}
 \ll q\sqrt{R\mathcal E}.
\end{equation*}
The coefficient product is bounded by the square of the coefficient
maximum in this range. Multiplying the separated sampling bound by
this estimate for the third moment gives a contribution proportional to
$B+{\mathsf{M}}^2/d^2$ for each common divisor.

Applying Cauchy's inequality and summing over $d$, we bound the
part with small gcds by
\begin{equation}\label{eq:de-small-gcd}
 \ll_c \mathbf 1_{\{d_0\geq1\}} A_\infty^2
       \left(Bd_0+{\mathsf{M}}^2\sum_{d\geq1}d^{-2}\right)\sqrt{R\mathcal E}
 \ll_c \mathbf 1_{\{{\mathsf{M}}^2\geq B\}}
       A_\infty^2{\mathsf{M}}^2\sqrt{R\mathcal E}.
\end{equation}

For the large gcds, collect all terms with the same reduced ratio
before taking a mean square.  For coprime positive integers $m,n$, put
\[
 b_{mn}=\sum_{\substack{d>d_0,\ {\mathsf{M}}\leq dm,dn\leq c{\mathsf{M}}\\(dmn,q)=1}}
                |a_{dm}a_{dn}|.
\]
Every ratio which occurs has
\[
 m,n\leq H_0:=c\min({\mathsf{M}},B/{\mathsf{M}})\ll_c\sqrt B,
 \qquad c^{-1}\leq m/n\leq c.
\]
The coefficients now add all possible common divisors leading to the
same coprime pair. This order of operations retains the cross terms
between different common divisors when the coefficient is squared.
The bound for the reduced numerator and denominator is at most a
constant times the square root of the scale for character and height: for any
positive length,
\begin{equation*}
 \min({\mathsf{M}},B/{\mathsf{M}})\le\sqrt B.
\end{equation*}
This will make the mean square of the polynomial in reduced ratios depend
on that scale without an additional length term.

Define
\[
 K(u,\psi)=\sum_{(m,n)=1}b_{mn}\psi(m)\overline{\psi(n)}(m/n)^{iu}.
\]
For each character $\psi\bmod q$ and integer $j$, let
\[
 \lambda_{\psi,j}=\#\{(w_1,w_2)\in W^2:
 \chi_1\overline\chi_2=\psi,\ t_1-t_2\in[j,j+1)\}.
\]
Two differences in the same bin yield an energy quadruple, after
permuting the four indices.  Hence
\begin{equation}\label{eq:de-difference-bin-energy}
 \sum_{\psi,j}\lambda_{\psi,j}^2\leq\mathcal E.
\end{equation}
Use $|\mathcal R|^3\leq R|\mathcal R|^2$ in the part with large gcds of
\eqref{eq:de-cubic}, expand the square, and group the difference pairs.
More explicitly, after grouping by the difference character and
ordinate bin, the sum is bounded by
\begin{equation*}
 R\sum_{\psi,j}\lambda_{\psi,j}
                  \sup_{u\in[j,j+1]}|K(u,\psi)|.
\end{equation*}
Cauchy's inequality separates the squared bin populations, controlled
by the energy, from the squared polynomial values.

Cauchy's inequality bounds this contribution by
\begin{equation}\label{eq:joint-cubic-bound}
 R\sqrt{\mathcal E}
 \left(\sum_{\psi\bmod q}\sum_{|j|\leq2T+2}
            \sup_{u\in[j,j+1]}|K(u,\psi)|^2\right)^{1/2}.
\end{equation}

Distinct reduced ratios with the same residue modulo $q$ have
nonzero cross difference divisible by $q$.  Their logarithms are
therefore separated by $\gg_c q/H_0^2$.
The Gaussian mean square proof of
Lemma~\ref{lem:de-fixed-modulus-mean}, now applied to these frequencies,
gives
\[
 \sum_{\psi\bmod q}\int_{-2T-3}^{2T+3}|K(u,\psi)|^2\,du
 \ll_c (B+H_0^2)\sum_{m,n}b_{mn}^2
 \ll_c B\sum_{m,n}b_{mn}^2.
\]
Sobolev gives the same bound for suprema over unit bins, since
$|\log(m/n)|\leq\log c$.

The positive coefficient norm satisfies
\[
 \sum_{(m,n)=1}b_{mn}^2
 \leq\sum_{(m,n)=1}
       \left(\sum_{d\geq1}|a_{dm}a_{dn}|\right)^2
 =\mathcal C(a).
\]
The equality follows by reducing equal fractions, or equivalently
by writing $(dm)(en)=(dn)(em)$.

For a direct parametrization, expand the convolution square norm as
\begin{equation*}
 \mathcal C(a)
 =\sum_{uv=xy}|a_u a_v a_x a_y|.
\end{equation*}
Given a solution, write $u=dm$ and $x=dn$ with $(m,n)=1$.
The equality becomes $mv=ny$, so coprimality forces $v=en$ and
$y=em$ for a unique positive integer $e$. Conversely these four
expressions always solve the equation. Thus the reduced pair and the
two common divisors parametrize exactly the same weighted solutions.

Substituting in \eqref{eq:joint-cubic-bound}, we bound the part with large gcds by
\begin{equation*}
 \ll_c R\sqrt{B\mathcal E\,\mathcal C(a)}.
\end{equation*}
Combining with \eqref{eq:de-small-gcd}, dividing by $\sqrt{\mathcal E}$
and squaring proves \eqref{eq:energy-coefficient-norm}.

The last algebraic step is
\begin{equation*}
 V^2\sqrt{\mathcal E}
 \ll_c \mathbf1_{{\mathsf{M}}^2\ge B}A_\infty^2{\mathsf{M}}^2\sqrt R
          +R\sqrt{B\mathcal C(a)}.
\end{equation*}
Squaring and using $(x+y)^2\le2x^2+2y^2$ produces the two terms in
the conclusion. Only the term with small common divisors uses a coefficient
maximum; the other term retains the full convolution square norm.
For another fixed positive tolerance, use narrower bins of fixed width
in \eqref{eq:de-R-moments} and \eqref{eq:de-difference-bin-energy}.
\end{proof}

\begin{corollary}
Under the hypotheses of Lemma~\ref{lem:difference-energy-fixed-log},
put \(L=\log(2+B{\mathsf{M}})\). If \(|a_n|\leq L^A d_\nu(n)\),
with fixed \(A\geq0\) and \(\nu\geq1\), then
\begin{equation}\label{eq:joint-ratio-norm}
 \mathcal C(a)\ll_{c,A,\nu}{\mathsf{M}}^2L^{4A+4\nu^2-3}.
\end{equation}
Consequently,
\begin{equation*}
 \mathcal E_1(W)\ll_{c,A,\nu}
 \mathbf1_{{\mathsf{M}}^2\geq B}A_\infty^4R{\mathsf{M}}^4V^{-4}
       +B{\mathsf{M}}^2R^2V^{-4}L^{4A+4\nu^2-3}.
\end{equation*}
\end{corollary}

\begin{proof}
For coprime \(m,n\), write
\(b^0_{mn}=\sum_{d\geq1}|a_{dm}a_{dn}|\).
Only \(m,n\ll_c {\mathsf{M}}\) with \(c^{-1}\leq m/n\leq c\) occur.
In a dyadic block \(h\leq\max(m,n)<2h\),
submultiplicativity and the mean of squared divisor functions give
\[
 b^0_{mn}\ll L^{2A}d_\nu(m)d_\nu(n)\frac {\mathsf{M}}h
          \bigl(\log(2+{\mathsf{M}}/h)\bigr)^{\nu^2-1}.
\]
Here we use
\(\sum_{d\leq y}d_\nu(d)^2\ll_\nu
y(\log(2y))^{\nu^2-1}\) for $y\geq1$.
Indeed $d_\nu(n)^2\leq d_{\nu^2}(n)$: at each prime, a pair
of exponent allocations among $\nu$ positions is represented by
at least one $\nu\times\nu$ matrix of nonnegative exponents.
Also, for every fixed integer $k\geq1$,
\[
 \sum_{n\leq y}d_k(n)
 \leq y\left(\sum_{m\leq y}\frac1m\right)^{k-1}
 \ll_k y(\log(2y))^{k-1}.
\]
These elementary inequalities prove the mean square bound used here.
After squaring and summing \(m,n\asymp_c h\), one block
contributes at most \({\mathsf{M}}^2L^{4A+4\nu^2-4}\).

The two factors in this logarithmic exponent can be seen separately.
The square of the sum over common divisors contributes at most
$L^{2\nu^2-2}$, and summing the two squared divisor functions over
variables of size $h$ contributes at most $h^2L^{2\nu^2-2}$.
Consequently one block is bounded by
\begin{equation*}
 L^{4A}\frac{{\mathsf{M}}^2}{h^2}L^{2\nu^2-2}
                h^2L^{2\nu^2-2}
 ={\mathsf{M}}^2L^{4A+4\nu^2-4}.
\end{equation*}
The additional logarithm comes only when the nonempty size blocks are
summed.
There are \(O(L)\) nonempty blocks; bounded scales satisfy
the same estimate. The identity for reduced ratios in the preceding
proof gives
\(\mathcal C(a)=\sum_{(m,n)=1}(b^0_{mn})^2\).
This proves \eqref{eq:joint-ratio-norm}.
Apply \eqref{eq:energy-coefficient-norm}.
\end{proof}

The first term still contains the coefficient maximum. To apply
\(d_\nu(n)\ll_{\nu,\varepsilon}n^\varepsilon\) there, a fixed
power margin is required.

For a Heath--Brown product, the original variable intervals also
restrict the coefficient norm. Expanding the coefficient square sum
leads to the following multiplicative equation and a sharper bound.

\begin{lemma}\label{lem:matrix-count}
Fix an integer $r\geq1$ and a support ratio $c>1$.  Suppose
$M_i\geq1$ and $\prod_{i=1}^rM_i=Y$.  The number of tuples of
positive integers satisfying
\[
 M_i\leq n_i,m_i\leq cM_i\quad(1\leq i\leq r),
 \qquad \prod_{i=1}^rn_i=\prod_{i=1}^rm_i
\]
is
\begin{equation*}
 \ll_{r,c}Y\bigl(1+\log Y\bigr)^{(r-1)^2}.
\end{equation*}
The intervals may be further restricted, and their weights may be
complex numbers of modulus at most one.
\end{lemma}

\begin{proof}
Allocating the prime exponents gives a positive integer matrix
$(a_{ij})$ with row products $n_i$ and column products $m_j$.
It therefore suffices to count these matrices.

For each prime, the total exponent on the row side equals the total
exponent on the column side. Distribute this common total among the
matrix entries, with the prescribed row and column sums, and then
multiply the prime powers assigned to each entry. Such a distribution
exists, for example by filling the first available row and column
until one of their remaining totals is exhausted. Every tuple is
therefore represented by at least one matrix, although a matrix
representation need not be unique. Counting all matrices is an upper
bound, which is the direction needed here.

Put $2^{b_{ij}}\leq a_{ij}<2^{b_{ij}+1}$, with $b_{ij}\geq0$.
The row and column restrictions imply
\[
 \sum_j b_{ij}=\log_2 M_i+O_{r,c}(1),\qquad
 \sum_i b_{ij}=\log_2 M_j+O_{r,c}(1).
\]
All entries are $O_{r,c}(1+\log Y)$.

For example, when there are two rows and columns, choosing the
size exponent in the upper left determines the other three up to bounded
choices: the first row determines the entry in the upper right, the first
column determines the lower left, and the last row determines the
last entry. The last column gives a consistency restriction. With
more rows, precisely the square in the upper left remains freely selectable.
This explains the square in the logarithmic exponent.

The upper left
$(r-1)\times(r-1)$ block determines the remaining entries up to
$O_{r,c}(1)$ choices by the row and column restrictions. Thus there are
$O_{r,c}((1+\log Y)^{(r-1)^2})$ size matrices, each admitting at most
\[
 \prod_{i,j}2^{b_{ij}+1}\ll_{r,c}Y,
\]
integer matrices. Restrictions and bounded weights preserve the bound.
\end{proof}

For $\nu\geq1$, let
\[
 c(n)=\sum_{n_1\cdots n_\nu=n}\prod_{j=1}^\nu |b_j(n_j)|,
 \qquad |b_j(n)|\leq\mathbf1_{[M_j,cM_j]}(n),\qquad
 {\mathsf{M}}=\prod_jM_j.
\]
Here the lower lengths are normalized to $M_j\geq1$, and
$L\geq1$ with $\log(2+{\mathsf{M}})\ll L$; all comparability constants are fixed.
A coefficient sequence is called a \emph{majorant for product intervals}
of order $\nu$ if $|a_n|\leq L^A c(n)$, with fixed support ratios.
A fixed number of bounded support refinements preserves this property.
Intervals containing only the integer one may be kept as atoms; they
do not worsen the bound.  Bounded changes when a lower endpoint is
less than one are absorbed into the support constants.

\begin{corollary}
Suppose $a_n$ has a majorant for product intervals of order $\nu$ and
logarithmic order $A$.  If its length is ${\mathsf{M}}$, then
\begin{equation}\label{eq:box-second-norm}
 \sum_n|a_n|^2\ll {\mathsf{M}} L^{2A+(\nu-1)^2}.
\end{equation}
If ${\mathsf{M}}\ll B=qT$, $T\geq1$, and $W$ is a set of character and height
pairs modulo the original modulus $q$, unit separated for each
character and with $|t|\leq T$, then
\begin{equation}\label{eq:box-second-moment}
 \sum_{(t,\chi)\in W}
 \left|\sum_n a_n\chi(n)n^{it}\right|^2
 \ll B{\mathsf{M}} L^{2A+(\nu-1)^2}.
\end{equation}
In particular this applies to every merged short terminal factor.
\end{corollary}

\begin{proof}
Expanding the positive majorant gives a multiplicative equation with
$\nu$ variables on each side. Lemma~\ref{lem:matrix-count}, with
$r=\nu$ and $Y={\mathsf{M}}$, proves \eqref{eq:box-second-norm}.
Explicitly, the positive coefficient square expands as
\begin{equation*}
 \sum_n c(n)^2
 =\sum_{n_1\cdots n_\nu=m_1\cdots m_\nu}
       \prod_{i=1}^\nu |b_i(n_i)b_i(m_i)|.
\end{equation*}
Each variable retains its own interval. The matrix lemma therefore
bounds this weighted count without replacing those intervals by a
single unrestricted divisor function.

The sampled form of Lemma~\ref{lem:de-fixed-modulus-mean} gives
\eqref{eq:box-second-moment}, since
$\varphi(q)(T+{\mathsf{M}}/q)\leq B+{\mathsf{M}}\ll B$.
Merged terminal factors satisfy ${\mathsf{M}}\ll\sqrt B$.
\end{proof}

\begin{corollary}
\label{cor:box-coefficient-norm}
For a majorant for product intervals of order $\nu$, put
\[
 d(\nu)=(2\nu-1)^2.
\]
Then
\begin{equation*}
 \sum_\ell |(a*a)(\ell)|^2
 \leq \mathcal C(a)
 \leq L^{4A}\sum_\ell(c*c)(\ell)^2
 \ll {\mathsf{M}}^2L^{4A+d(\nu)}.
\end{equation*}
Here and below $L$ may be replaced by any comparable logarithm of
all polynomially bounded scales.
\end{corollary}

\begin{proof}
After taking absolute values and expanding, apply
Lemma~\ref{lem:matrix-count} with $r=2\nu$, $Y={\mathsf{M}}^2$ and the interval
list $(M_1,\ldots,M_\nu,M_1,\ldots,M_\nu)$.

There are twice as many original variables because the polynomial was
squared before its coefficient square norm was taken. The product of
their lower length scales is ${\mathsf{M}}^2$, while the number of free matrix
size entries is $(2\nu-1)^2$. These give, respectively, the polynomial
scale and logarithmic exponent in the bound.
\end{proof}

\begin{corollary}
\label{cor:box-energy}
Under the hypotheses of Lemma~\ref{lem:difference-energy-fixed-log},
suppose additionally that the coefficients admit a majorant for product intervals of order $\nu$.  Then
\begin{equation}\label{eq:box-energy}
 \mathcal E(W)V^4\ll
 \mathbf1_{{\mathsf{M}}^2\geq B}A_\infty^4R{\mathsf{M}}^4
       +B{\mathsf{M}}^2R^2L^{4A+(2\nu-1)^2}.
\end{equation}
If ${\mathsf{M}}^2\ll B$, the sampled fourth moment, including all characters
of the original modulus, satisfies
\begin{equation}\label{eq:box-fourth}
 \sum_{(t,\chi)\in W}|F(t,\chi)|^4
       \ll B{\mathsf{M}}^2L^{4A+(2\nu-1)^2}.
\end{equation}
\end{corollary}

\begin{proof}
The term with small gcds in the proof for difference energy is unchanged.
For the joint coefficients for large gcds $b_{mn}$ of that proof, positivity
gives
\[
 \sum_{(m,n)=1}b_{mn}^2
 \leq L^{4A}\sum_{(m,n)=1}
          \left(\sum_{d\geq1}c(dm)c(dn)\right)^2
 =L^{4A}\sum_\ell(c*c)(\ell)^2.
\]
The identity for reduced ratios and Corollary~\ref{cor:box-coefficient-norm}
in the preceding energy proof give \eqref{eq:box-energy}.
Applying the sampled mean square to $F^2$ gives
\eqref{eq:box-fourth}, since $B+{\mathsf{M}}^2\ll B$.

The identity behind this application is
\begin{equation*}
 |F(t,\chi)|^4=|F(t,\chi)^2|^2,
 \qquad
 F(t,\chi)^2=\sum_\ell(a*a)(\ell)\chi(\ell)\ell^{it}.
\end{equation*}
Its length is comparable to ${\mathsf{M}}^2$. Thus the mean square applies to
the convolution coefficients just estimated, and the assumed
shortness prevents its new length from exceeding the main
scale for character and height. No restriction to primitive characters is
introduced in this calculation.
\end{proof}

The interval bound saves $4\nu-4$ logarithms over the divisor bound.
The coefficient maximum in the first energy term still requires a
fixed power margin.

\section{Fixed logarithmic fourth moments for the original factors}
\label{sec:original-factor-fourth-fixed-log}
\label{subsec:of-HB-bookkeeping}

For short factors the fourth moment includes all characters of the
original modulus. For long ordinary factors the contour shift requires
nonprincipal characters, the complete principal term having been
estimated separately.

\begin{lemma}
\label{lem:of-full-character-L-fourth}
For $q\geq1$ and $T\geq1$,
\begin{equation}\label{eq:of-full-L-fourth}
 \sum_{\substack{\chi\bmod q\\\chi\ne\chi_{0,q}}}
 \int_{-T}^T|L(1/2+it,\chi)|^4\,dt
 \ll qT\bigl(\log(2qT)\bigr)^4.
\end{equation}
The implied constant is absolute.  In particular there is no
factor $q^\varepsilon$, and no primitivity requirement on the
characters in this conclusion.
\end{lemma}

\begin{proof}
We estimate the fourth moment by writing the square of the
$L$-function as a sum of two smoothed Dirichlet series and then taking
a second moment. The functional equation is the analytic tool which
allows this representation on the critical line, where the original
Dirichlet series need not converge absolutely. The weights in the
representation decrease rapidly once their integer argument exceeds
the conductor times the height. Thus the effective polynomial length
is comparable to that same product.

For primitive characters modulo $r>1$, on $U\leq|t|\leq2U$, $U\geq1$,
the smoothed approximate functional equation at $s=1/2+it$ is
\begin{align}
 L(s,\chi)^2
 &=\sum_{n\geq1}\frac{d_2(n)\chi(n)}{n^s}
                      V_{\mathfrak a,t}(n/r) \notag\\
 &\quad+\varepsilon_\chi(t)
       \sum_{n\geq1}\frac{d_2(n)\overline\chi(n)}{n^{1-s}}
                      V_{\mathfrak a,-t}(n/r),
 \qquad |\varepsilon_\chi(t)|=1,                 \label{eq:of-AFE}
\end{align}
where $\mathfrak a\in\{0,1\}$ is the parity of $\chi$ and
\begin{equation}\label{eq:of-AFE-weight}
 V_{\mathfrak a,t}(y)=\frac1{2\pi i}\int_{(2)}
 \frac{e^{z^2}}z
 \frac{\Gamma((1/2+it+z+\mathfrak a)/2)^2}
      {\Gamma((1/2+it+\mathfrak a)/2)^2}
 (\pi y)^{-z}\,dz.
\end{equation}
To obtain this identity, integrate the completed $L(s+z,\chi)^2$
against $e^{z^2}/z$, divided by the completion factor at $s$.
The completed function is entire, so the contour shift crosses only
$z=0$. Apply its functional equation and replace $z$ by $-z$;
expansion on the lines on the right gives \eqref{eq:of-AFE}.

The factor $1/z$ has residue one at zero. Crossing it therefore
recovers the square of the $L$-function that we started with. The
functional equation changes the remaining integral on the left line
into a second integral on a right line, where its Dirichlet series can
again be expanded. The dual character is conjugated, and its multiplier
has absolute value one. Both pieces can consequently be estimated by
the same mean square argument.

Put $v=t/U$, $x=n/(rU)$, and
$W_{\mathfrak a,U}(v,x)=V_{\mathfrak a,Uv}(Ux)$.
On a compact $v$ interval away from zero, Stirling's formula gives
uniform bounds for every fixed number of $v$ and $x\partial_x$
derivatives. Shifting to $\Re z=-1/8$ gives
$W_{\mathfrak a,U}(v,x)=1+O(x^{1/8})$ at zero; positive shifts give
rapid decay at infinity, with the same derivative bounds.
Indeed derivatives in $v$ differentiate ratios of gamma functions
at arguments separated by $z$; after multiplication by $U$ the
resulting logarithmic derivatives are $O_j((1+|z|)^{O_j(1)})$
on $|\Im z|\leq U/4$, and the complementary part is controlled
by $e^{-({\Im z})^2}$.  The compact block at low height $|t|\leq2$
has the same estimates with $U=1$ and all its $v$ in a fixed
compact interval, directly from \eqref{eq:of-AFE-weight}.

Choose a fixed smooth $w_0(x)$ which is one near zero and zero
for $x\geq2$.  The portion $w_0(n/(rU))$ has length $O(rU)$.
For support $1\leq n\leq M$, the inequality
$|\log(n/m)|\geq|n-m|/M$ also permits the continuous Gaussian
argument of Lemma~\ref{lem:de-fixed-modulus-mean}.
Together with
\[
 \sum_{n\leq Y}\frac{d_2(n)^2}{n}
 \leq\sum_{n\leq Y}\frac{d_4(n)}n
 \ll\bigl(1+\log Y\bigr)^4,
\]
this estimate for the mean value bounds the squared $L^2$ norm of this portion over all
characters modulo $r$ and the shell by
$O(\varphi(r)U(\log(2rU))^4)$.

In the mean value bound the length term divided by the modulus is
$O(U)$, so it is no larger than the integration height term. The
coefficient square norm is at most a fourth power of a logarithm:
the inequality between the two divisor functions holds prime by prime,
and
\begin{equation*}
 \sum_{n\le Y}\frac{d_4(n)}n
 \le\left(\sum_{m\le Y}\frac1m\right)^4
 \ll(1+\log Y)^4.
\end{equation*}
The expansion of the fourth power on the right includes every ordered
factorization contributing to the sum on the left. This supplies the
explicit logarithmic exponent in the moment estimate.

For $W-w_0$, first multiply in $v$ by a fixed smooth compact cutoff
equal to one on the height block in use. Partition into $x\asymp2^j$
and rescale to a fixed compact interval in $x/2^j$.
Fourier inversion in $v$ and Mellin
inversion in $x/2^j$ separate the weight, with transform $L^1$ norm
\[
 \ll \begin{cases}2^{j/8},&j\leq0,\\2^{-10j},&j>0.
 \end{cases}
\]
by two integrations by parts in each variable. Retaining a fixed
cutoff in the integer variable gives polynomials of length
$\asymp2^jrU$ with coefficients bounded by $Cd_2(n)/\sqrt n$,
possibly twisted by $n^{iu}$. The separated time phase has modulus one.
Lemma~\ref{lem:de-fixed-modulus-mean} bounds their $L^2$ norm,
including the character sum, by
\[
 \ll\{\varphi(r)U(1+2^j)\}^{1/2}
          (\log(2rU)+|j|)^2.
\]
The weight separation is used inside an integral, so the triangle
inequality needed here is the integral form of Minkowski's inequality:
the norm of an integral is at most the integral of the norms. The
negative size blocks have transform weights that form a convergent
geometric series. For positive blocks the norm may grow like the square
root of their length, but the transform weight decreases much faster.
The sums are therefore convergent without adding a factor equal to
the number of size blocks.

Minkowski and the summable transform norms give
$O((\varphi(r)U)^{1/2}(\log(2rU))^2)$ for the remainder, for each parity.
The dual term has the same bound by character conjugation and
$|\varepsilon_\chi(t)|=1$. Squaring and summing the height shells gives
\begin{equation}\label{eq:of-primitive-L-fourth}
 \sum_{\chi\bmod r}^{*}\int_{-T}^T|L(1/2+it,\chi)|^4\,dt
 \ll\varphi(r)T(\log(2rT))^4\qquad(r>1).
\end{equation}

If a character modulo $q$ is induced by the primitive $\chi^*$
of conductor $r>1$, then on the critical line
\[
 |L(s,\chi)|^4\leq
 |L(s,\chi^*)|^4
       \prod_{\substack{p\mid q\\p\nmid r}}(1+p^{-1/2})^4.
\]
Induction is a bijection between characters of conductor $r$
modulo $q$ and primitive characters modulo $r$.

The removed Euler factors arise because an induced character is zero
at primes dividing the larger modulus. Their individual fourth powers
can be larger than one. Instead of bounding their product uniformly
and losing an unspecified small power of the modulus, sum over all
possible conductors first. The next Euler product performs exactly
that sum, one prime at a time.
Write $H_p=(1+p^{-1/2})^4$.  The exact elementary Euler product
\begin{align}
 \sum_{r\mid q}\varphi(r)
       \prod_{\substack{p\mid q\\p\nmid r}}H_p
 &=\prod_{p^a\Vert q}\left(H_p+\sum_{j=1}^a\varphi(p^j)\right)
   =q\prod_{p^a\Vert q}\left(1+\frac{H_p-1}{p^a}\right)
   \ll q                                                   \label{eq:of-Euler-recombine}
\end{align}
holds because $(H_p-1)/p^a=O(p^{-3/2})$.

For a prime power dividing the modulus, conductor exponent zero gives
the first local term. Positive conductor exponents give
\begin{equation*}
 \sum_{j=1}^a\varphi(p^j)=p^a-1.
\end{equation*}
After the common factor $p^a$ is removed, the excess local factor is
summable over primes. The product is therefore bounded independently
of the modulus, which is what preserves a fixed logarithmic loss.
Combining \eqref{eq:of-primitive-L-fourth} for $r>1$ with
\eqref{eq:of-Euler-recombine}, we obtain \eqref{eq:of-full-L-fourth}.
Here $r=1$ is added only to the positive arithmetic sum.
\end{proof}

\begin{lemma}
\label{lem:original-factor-fourth-fixed-log}
Let $q\geq1$, $T\geq1$, $B=qT$, ${\mathsf{M}}\geq1$, and
$L=\log(2+B{\mathsf{M}})$.  Suppose that $W$ is a set of points $(t,\chi)$
with $|t|\leq T$, separated by one for each character.
The following two estimates are uniform in all these data.
\begin{enumerate}
\item If ${\mathsf{M}}^2\ll_c B$ and
\[
 F(t,\chi)=\sum_{{\mathsf{M}}\leq n\leq c{\mathsf{M}}}a_n\chi(n)n^{it},
 \qquad |a_n|\leq C L^A d_\nu(n),
\]
then, including every character modulo $q$,
\begin{equation}\label{eq:of-short-fourth}
 \sum_{(t,\chi)\in W}|F(t,\chi)|^4
 \ll_{c,A,\nu,C}B{\mathsf{M}}^2 L^{4A+4\nu^2-3}.
\end{equation}
\item Suppose that all characters in $W$ are nonprincipal modulo
$q$.  If $w$ belongs to a fixed bounded family of smooth functions
with compact support in $(0,\infty)$, and $A=0$ or $1$, then for
any ${\mathsf{M}}\geq1$,
\begin{equation*}
 \sum_{(t,\chi)\in W}
 \left|\sum_n w(n/{\mathsf{M}})(\log n)^A\chi(n)n^{it}\right|^4
 \ll_w B{\mathsf{M}}^2L^{4+4A}.
\end{equation*}
There is no relation required between ${\mathsf{M}}$ and $B$ in this case.
\end{enumerate}
The same statements hold for suprema over unit bins.  In particular, each
terminal factor can be charged its own exponent
\begin{equation*}\label{eq:factor-costs}
 c_i=\begin{cases}
 4A_i+4\nu_i^2-3,&\text{for a short factor},\\
 4+4A_i,&\text{for an original ordinary factor}.
 \end{cases}
\end{equation*}
\end{lemma}

\begin{proof}
There are two separate reasons a factor can have a manageable fourth
moment. For a short factor, squaring still leaves a polynomial within
the length range of the elementary mean square. For a long ordinary
factor, its coefficients allow Mellin inversion to express it in terms
of an $L$-function, whose fourth moment was just estimated. The second
argument needs nonprincipal characters to avoid a residue at one.

For short factors, the coefficients $c_\ell$ of $F^2$ satisfy
\eqref{eq:joint-ratio-norm}, so
\[
 \sum_\ell|c_\ell|^2\ll {\mathsf{M}}^2L^{4A+4\nu^2-3}.
\]
Lemma~\ref{lem:de-fixed-modulus-mean}, applied to $F^2$ with
$B+{\mathsf{M}}^2\ll B$, gives \eqref{eq:of-short-fourth} and its supremum version.

For an ordinary factor with $A=0$, Mellin inversion and a contour
shift to real part $1/2$ give
\[
 P_{\mathsf{M}}(t,\chi)=\sum_nw(n/{\mathsf{M}})\chi(n)n^{it}
 =\frac{{\mathsf{M}}^{1/2}}{2\pi}\int_{\mathbb R}
    \widehat w(1/2+iu){\mathsf{M}}^{iu}L(1/2+i(u-t),\chi)\,du.
\]
There is no pole since $\chi$ is nonprincipal.

Before the shift, Mellin inversion gives an integral of
$\widehat w(s){\mathsf{M}}^sL(s-it,\chi)$ on a line to the right of one.
Moving that line to real part one half extracts the factor ${\mathsf{M}}^{1/2}$.
Its fourth power is ${\mathsf{M}}^2$, explaining the length scale in the result.
If the character were principal, the crossed pole would instead give a
main term of size comparable to ${\mathsf{M}}$, whose fourth power could not be
bounded by the same expression for arbitrary length.
The Mellin transform is rapidly decreasing with uniform constants.
Minkowski's inequality and
Lemma~\ref{lem:of-full-character-L-fourth}, applied at height
$T+|u|+1$, imply
\[
 \sum_{\substack{\chi\bmod q\\\chi\ne\chi_{0,q}}}
       \int_{-T}^T|P_{\mathsf{M}}(t,\chi)|^4\,dt
 \ll_w B{\mathsf{M}}^2 L^4.
\]
Apply Sobolev to ${\mathsf{M}}^{-2it}P_{\mathsf{M}}^2$. Its derivative introduces the
ordinary weight $iw(x)\log x$, so Cauchy's inequality and the
continuous fourth moments give the same bound. Bounded overlap of
enlarged bins proves the sampled and supremum versions.
For $A=1$, use $\log n=\log {\mathsf{M}}+\log(n/{\mathsf{M}})$, which costs $L^4$;
normalized differentiation again introduces only bounded logarithms.

For the sampling step, differentiating the normalized square produces
terms with one of the two smooth factors replaced by the same factor
weighted by a bounded logarithm of its scaled variable. These belong
to a fixed bounded smooth family, so the continuous estimate applies
to them too. For the logarithmic factor, the decomposition just given
is a sum of an ordinary polynomial multiplied by a constant logarithm
and another ordinary smooth polynomial. The inequality
$|a+b|^4\le8(|a|^4+|b|^4)$ gives the stated extra fourth power of a
logarithm.
\end{proof}

For a fixed smooth weight for the prime polynomial at scale $X$, choose
$C_0>1$ larger than its support endpoint and put $z=(C_0X)^{1/5}$.
For $n\leq z^5$, Heath--Brown's identity is
\begin{equation*}\label{eq:of-HB-identity}
 \Lambda(n)=\sum_{j=1}^{5}(-1)^{j-1}\binom5j
  \bigl(\mu_{\leq z}^{*j}*1^{*(j-1)}*\log\bigr)(n).
\end{equation*}
To verify it, put $M=\mu_{\leq z}$ and use the convolution
identities $\Lambda=\mu*\log$ and $\mu*1=\delta$.
The difference between the two sides is
$\Lambda*(\delta-1*M)^{*5}$.
Since $\delta-1*M=1*(\mu-M)$ is supported on integers greater
than $z$, this difference vanishes at $n\leq z^5$.

The convolution identity can also be seen by expanding the binomial:
\begin{equation*}
 \delta-(\delta-1*M)^{*5}
 =\sum_{j=1}^5(-1)^{j-1}\binom5j(1*M)^{*j}.
\end{equation*}
Convolve both sides with the von Mangoldt function. The elementary
identity $1*\Lambda=\log$ changes each term on the right into the
term displayed in the Heath--Brown identity. In the omitted
remainder, each of the five factors is supported on integers strictly
greater than the cutoff. Their product is consequently larger than
the fifth power of the cutoff, so it cannot contribute in the stated range. This
finite identity is the reason prime coefficients can be replaced by a
bounded number of simpler convolution factors.

A smooth dyadic partition gives $O((\log X)^9)$ boxes: the product
condition fixes the last of at most ten lengths up to bounded choices.
There are at most five M\"obius atoms $\mu(n)\mathbf1_{n\leq z}$;
the others are ordinary smooth factors, with one logarithmic factor.
The lengths multiply to $\asymp X$. One Mellin integral separates
the outer weight and shifts every factor by the same real number.

For a box with variables $n_1,\ldots,n_k$, the outer cutoff depends
on their product. Its Mellin inversion contains
\begin{equation*}
 X^{iu}(n_1\cdots n_k)^{-iu}
 =X^{iu}\prod_{i=1}^k n_i^{-iu}.
\end{equation*}
The original factor $n_i^{-it}$ therefore becomes
$n_i^{-i(t+u)}$ for every variable. There is one common shift, not a
separate choice of shift for each factor. Dyadic partitions mean that
each variable is restricted to an interval between a length and a fixed
multiple of that length; this is the interval information retained in
the coefficient estimates above.
Atoms supported at $1$ are retained in the order count;
their moment bounds also follow from $|W|\ll qT$.

Assume $X^{\beta_0}\leq B\ll X^{1/2}(\log X)^2$, where
$\beta_0>4/9$ is fixed.
Repeatedly merge the two shortest factors whenever their length
product is less than $c_0\sqrt B$, for a fixed constant $c_0>0$.
A merged factor is a Dirichlet convolution of its original atoms.
Its length is less than $c_0\sqrt B$, with a fixed support ratio,
and its coefficients are bounded by
\begin{equation*}\label{eq:of-coefficient-orders}
 C L^{A_i}d_{\nu_i}(n),\qquad
 \nu_i\leq10,\quad A_i\leq1,\quad
 \sum_i\nu_i\leq10,\quad\sum_i A_i\leq1,
\end{equation*}
where now $L=\log(2+X B)$.
This follows by convolution of the original bounded atoms, with the
sole logarithm factored out. All support ratios remain fixed.

An atom here is one original factor in one variable before merging.
Convolving two bounded atoms replaces their coefficients by a sum over
decompositions of an integer into two factors, bounded by the corresponding
divisor function. Further convolution increases its order by addition.
There are at most ten original atoms and only one logarithmic weight,
so the orders and the total logarithmic degree obey the displayed
bounds irrespective of the order of merging.

Since $z\ll X^{1/5}=o(B^{1/2})$, each terminal factor is short or
an original ordinary factor. Lemma~\ref{lem:original-factor-fourth-fixed-log}
therefore gives, over nonprincipal characters of the original modulus,
\begin{equation*}
 RV_i^4\ll B{\mathsf{M}}_i^2L^{c_i}.
\end{equation*}
For a subset $A$, Lemma~\ref{lem:difference-energy-fixed-log} applies
with $\nu_A=\sum_{i\in A}\nu_i\leq10$ and
$A_A=\sum_{i\in A}A_i\leq1$; its second energy term has logarithmic
exponent $e_A=4A_A+4\nu_A^2-3$.

At termination, ${\mathsf{M}}_i{\mathsf{M}}_j\geq c_0\sqrt B$ for $i\ne j$, whence
$\prod_i{\mathsf{M}}_i\gg_k B^{k/4}$ if $k\geq2$.
Since $\prod_i{\mathsf{M}}_i\asymp X$ and $(9/4)\beta_0>1$, we have $k\leq8$
for sufficiently large $X$.

To verify the product inequality, multiply the inequalities for pair lengths
over all unordered pairs. Each length occurs in exactly $k-1$ pairs,
so
\begin{equation*}
 \left(\prod_i{\mathsf{M}}_i\right)^{k-1}
 \ge(c_0\sqrt B)^{k(k-1)/2},
 \qquad
 \prod_i{\mathsf{M}}_i\ge c_0^{k/2}B^{k/4}.
\end{equation*}
If there were at least nine factors, the lower bound would grow faster
than $X$ because $B\ge X^{\beta_0}$ and $9\beta_0/4>1$.
The contradiction proves the stated bound on their number. A factor
formed by merging is short by the rule of the construction, while an
unmerged long factor cannot be a M\"obius atom because its original
cutoff is smaller than the shortness scale.

\paragraph{The common Mellin shift.}
Choose the outer and dyadic weights in a fixed bounded Gevrey family
of order $1<s\leq2$.  Their Mellin tails allow a common shift
$|u|\leq C_A L^s$ with error $O(X^{-A})$ for every fixed $A$.
The shifted factors have enlarged height
$U\asymp T+L^s$, so the preceding moments apply with $B=qU$.

This height follows from the triangle inequality
$|t+u|\le |t|+|u|$. It is an additive enlargement: the original
height and the logarithmic truncation height need not be multiplied.
The same shift multiplies each original coefficient by a phase of
absolute value one, preserving all the positive majorants and support
intervals used above.
At logarithmic precision the count of product boxes and recombination
of the Mellin integral remain to be included, as in
Section~\ref{sec:endpoint-gauss-bridge}.

\section{Complete values, common products and the circle integral}
\label{sec:ab-complete-bridge}

The energy and weighted mass must refer to the same amplitude level.
Thus the sample set is chosen from the complete prime polynomial,
before a product is selected at each point.
Proposition~\ref{prop:principal-lower-theta} treats the principal
contribution, so all long factors below are nonprincipal.

\subsection{Complete polynomials and energy}

Fix $w\in C_c^\infty((0,1))$ and a function
$u\in C_c^\infty((0,\infty))$ equal to $1$ on $\operatorname{supp}w$.
Write
\[
 \Psi_u(t,\chi)=\sum_n\Lambda(n)u(n/X)\chi(n)n^{-it},\qquad
 \mathcal D_u(t,\chi)=X^{it}\Psi_u(t,\chi).
\]
Thus $|\mathcal D_u|=|\Psi_u|$.
Smoothing functions have fixed derivative bounds. We use
$\widehat f(s)=\int_0^\infty f(y)y^{s-1}\,dy$ and put
$\mathcal H=X^\eta$, choosing $\eta>0$ after the final exponent.
The notation $\preccurlyeq$ permits an arbitrarily small power of $X$.

For a set $W$ of pairs $(t,\chi)$, with characters of one modulus,
define the ordered energy
\begin{gather*}
 \mathcal E_L(W)=\#\{(w_1,w_2,w_3,w_4)\in W^4:
   \chi_1\chi_2=\chi_3\chi_4,\
   |t_1+t_2-t_3-t_4|\le L\},\\*
 \mathcal E(W)=\mathcal E_1(W).
\end{gather*}
A set is unit separated if, for each fixed character, its distinct
ordinates are at distance at least $1$.
Empty value levels are omitted throughout.

The factors use $n^{-it}$; reflection $t\mapsto-t$ preserves
separation, energy and the preceding estimates.

For example, expanding a fourth power of a sum over a set of points
produces the character product
$\chi_1\chi_2\overline{\chi_3\chi_4}$ and the ordinate discrepancy
$t_1+t_2-t_3-t_4$. Orthogonality tests whether the character product is
principal, while a smooth integral in the other variable suppresses
large discrepancies. The energy records the quadruples that remain
close enough to contribute. It must be measured on the same set on
which the polynomial has the chosen lower value bound.

\subsection{Mean values and adaptive sampling}

The sampled mean values will follow from the classical multiplicative
large sieve and fourth moment of Dirichlet $L$-functions in the forms
\begin{align}
 \sum_{r\le R}\frac r{\varphi(r)}
 \sum_{\chi\bmod r}^{*}\int_{-U}^{U}
 \left|\sum_{n\asymp {\mathsf{M}}}a_n\chi(n)n^{-it}\right|^2dt
 &\ll (R^2U+{\mathsf{M}})(\log(2{\mathsf{M}}))^C\sum_n|a_n|^2,
 \label{eq:ab-hybrid-ls}\\
 \sum_{\chi\bmod r}^{*}\int_{-U}^{U}
       |L(1/2+it,\chi)|^4dt
 &\ll_\varepsilon(rU)^{1+\varepsilon},
 \qquad r>1,\quad U\ge1.\notag
\end{align}
The large sieve is used in the form of
\cite{Montgomery1971}. For $r>1$, all primitive characters
modulo $r$ are nonprincipal, so the fourth moment follows from the
estimate with a fixed logarithmic loss in
Lemma~\ref{lem:of-full-character-L-fourth}.
For \eqref{eq:ab-hybrid-ls}, Gallagher's lemma bounds the integral by
\[
 U^2\int_{\mathbb R}
 \left|\sum_{e^y<n\le e^{y+1/U}}a_n\chi(n)\right|^2dy.
\]
The multiplicative large sieve on these intervals of length
$O({\mathsf{M}}/U+1)$ gives $(R^2U+{\mathsf{M}}+U)\sum|a_n|^2$, since each coefficient
occurs on a set of $y$ of measure $1/U$. Absorb $U$ into $R^2U$.

The factor accounting for the short interval and integration can be
tracked directly:
\begin{equation*}
 U^2\left(R^2+\frac {\mathsf{M}}U+1\right)\frac1U
 =R^2U+{\mathsf{M}}+U.
\end{equation*}
The last factor on the left is the measure of the set of interval
starting points for which a fixed coefficient is present. The
multiplicative large sieve used at this step is an external mean value
theorem; this calculation explains how its interval length enters the
hybrid bound.
For a fixed modulus the corresponding mean square is
\begin{equation}
 \sum_{\chi\bmod r}\int_{-U}^{U}
  \left|\sum_{n\asymp {\mathsf{M}}}a_n\chi(n)n^{-it}\right|^2dt
 \ll(rU+{\mathsf{M}})(\log(2{\mathsf{M}}))^C\sum_n|a_n|^2,
 \label{eq:ab-fixed-mean}
\end{equation}
which follows from Lemma~\ref{lem:de-fixed-modulus-mean}.

For adaptively chosen unit separated $t_j$ and continuously
differentiable $F$, Sobolev gives
\begin{equation}
 \sum_j|F(t_j)|^2
 \ll\int_{\bigcup_j[t_j-1/2,t_j+1/2]}
             (|F(t)|^2+|F(t)F'(t)|)\,dt.
 \label{eq:ab-sampling}
\end{equation}
For Dirichlet polynomials, differentiation inserts $\log n$.
Applying Cauchy's inequality with \eqref{eq:ab-fixed-mean} or
\eqref{eq:ab-hybrid-ls}, we obtain the sampled versions with logarithmic
losses. For $|F|^4$, the same argument with H\"older gives sampled
fourth moments.

No point needs to be selected in advance of the polynomial. The
calculus inequality is valid for every point of each small interval,
so it also controls a point chosen after the values have been seen.
Separation is used to keep the enlarged intervals from overlapping
more than a bounded number of times. This distinction allows the
representatives of sets of large values later in the proof to be chosen
adaptively.

\subsection{Smooth Heath--Brown factors and their fourth moments}

\begin{lemma}
\label{lem:ab-smooth-HB}
Let $X^\kappa\le B\le X^{1/2}$, where $\kappa>4/9$ is fixed.
For every Dirichlet character modulo $q$, one has an identity
\begin{equation}
 \Psi_u(t,\chi)=
 \sum_{\nu\ll(\log X)^C}c_\nu
 \frac1{2\pi}\int_{\mathbb R}
       \widehat u(i\zeta)X^{i\zeta}
       F_\nu(t+\zeta,\chi)\,d\zeta,
 \qquad F_\nu=\prod_{i=1}^{k_\nu}S_{i,\nu},
 \label{eq:ab-smooth-HB}
\end{equation}
where $|c_\nu|\ll1$, $k_\nu\le8$ for sufficiently large $X$, and
the factors have common coefficients bounded by divisor functions and support
in $[{\mathsf{M}}_{i,\nu},C_0{\mathsf{M}}_{i,\nu}]$, with
\[
 \prod_i{\mathsf{M}}_{i,\nu}\asymp X,\qquad
 {\mathsf{M}}_{i,\nu}{\mathsf{M}}_{j,\nu}\gg B^{1/2}\quad(i\ne j).
\]
Every factor of length greater than a fixed multiple of $B^{1/2}$
is a smoothly weighted sum with coefficients $1$ or $\log n$.
The integral in \eqref{eq:ab-smooth-HB} can be truncated to
$|\zeta|\le \mathcal H$ with error $O_A(X^{-A})$ for any prescribed $A$,
by choosing its fixed decay order sufficiently large.
\end{lemma}

\begin{proof}
The decomposition and merging in Section~\ref{subsec:of-HB-bookkeeping}
give the stated factors. Mellin inversion gives \eqref{eq:ab-smooth-HB}.
The transform decays faster than any fixed power, whereas each product
is polynomially bounded uniformly in its real shift. This permits
truncation at $X^\eta$.

If the product is bounded by a fixed power $X^{C_1}$ and the Mellin
transform is $O_M((1+|\zeta|)^{-M})$, the discarded part has absolute
value at most a constant times
\begin{equation*}
 X^{C_1}\int_{|\zeta|>X^\eta}(1+|\zeta|)^{-M}\,d\zeta
 \ll_M X^{C_1-\eta(M-1)}.
\end{equation*}
Choose the fixed integer $M$ after the desired error exponent and the
fixed positive height parameter. The error then has the required size,
uniformly for the characters and ordinates being considered.
\end{proof}

\begin{lemma}
Let $T \ge 1$. Fix a modulus $q$, a divisor $r\mid q$ with $r>1$, and let
$\chi_q$ be induced from a primitive nonprincipal character
$\chi^*\bmod r$.
For every factor of Lemma~\ref{lem:ab-smooth-HB},
and every unit separated common set $W$ of such characters and
ordinates $|t|\le T\mathcal H$, one has
\begin{equation*}
 |W|V_i^4\preccurlyeq X^{O(\eta)}B {\mathsf{M}}_i^2
 \quad\hbox{whenever }|S_i(t,\chi_q)|\ge V_i\text{ on }W,
 \qquad B=qT.
\end{equation*}
The constants are uniform under restriction of the sample set and over
fixed bounded families of the original smooth weights.  A long ordinary
factor is always used with its original smooth weight; support pieces
or coefficient changes in other tests do not replace that factor.
A common Mellin shift is included in the evaluation ordinate and hence
in the height of the moment estimate.
\end{lemma}

\begin{proof}
Apply Lemma~\ref{lem:original-factor-fourth-fixed-log} at height
$T\mathcal H$, with parameter $B\mathcal H$. The induced characters
form a nonprincipal subset modulo $q$; merged and M\"obius factors
are short, and the remaining long factors are ordinary.
Fixed logarithmic losses are absorbed into the allowed small power.

On the common set, the lower threshold gives
\begin{equation*}
 |W|V_i^4\le\sum_{(t,\chi_q)\in W}|S_i(t,\chi_q)|^4.
\end{equation*}
The full nonprincipal moment at the original modulus bounds the sum
on the right, because the specified induced characters form a subset
of that family. Its scale for character and height is enlarged by the same
height factor. A fixed power of a logarithm is smaller than any fixed
positive power of $X$ for sufficiently large $X$, giving the stated
form of the estimate.
\end{proof}

\subsection{Positive energy and selection of complete polynomials}

\begin{lemma}
\label{lem:ab-energy-transfer}
Suppose $W$ is unit separated, and every $(t,\chi)\in W$ is
mapped to $(t',\chi)$ with $|t'-t|\le \mathcal H$.
Partition the image, counted with its multiplicities, into $L$
sets $W_c'$ on which the map is injective.
Then
\begin{equation*}
 \mathcal E(W)\ll (1+\mathcal H)L^3\sum_{c=1}^L \mathcal E(W_c').
\end{equation*}
In particular this is an $X^{O(\eta)+O(\varepsilon)}$ loss if
$L\ll X^{O(\eta)+O(\varepsilon)}$.
\end{lemma}

\begin{proof}
Let $r$ be the common modulus. The image of an energy quadruple has
discrepancy at most $1+4\mathcal H$.

Each of the four ordinate shifts changes the discrepancy by at most
the height allowance, and the original discrepancy is at most one.
The character product condition is unchanged. Thus the original
energy is bounded by the energy of the image multiset with a larger
tolerance. The multiset is needed because different original points
may have the same image; their multiplicities cannot be discarded.

For $K\geq1$ put
\[
 \mathcal G_K(W)=
 \sum_{\chi_1\chi_2=\chi_3\chi_4}
 \exp\!\left(-\frac{(t_1+t_2-t_3-t_4)^2}{K^2}\right).
\]
Fourier inversion for the Gaussian and character orthogonality
give
\[
 \mathcal G_K(W)=\frac{K}{2\sqrt\pi\,\varphi(r)}
 \sum_{a\bmod r}^{*}\int_{\mathbb R}e^{-K^2\xi^2/4}
 \left|\sum_{(t,\chi)\in W}\chi(a)e^{it\xi}\right|^4d\xi.
\]
The identity also holds for multisets, so H\"older gives
$\mathcal G_K(\bigsqcup_c W_c')\leq L^3\sum_c\mathcal G_K(W_c')$.

The advantage of the Gaussian is that its Fourier transform is
nonnegative. After orthogonality, the resulting integral is therefore
an integral of a nonnegative fourth power. Pointwise,
\begin{equation*}
 \left|\sum_{c=1}^L Z_c\right|^4
 \le L^3\sum_{c=1}^L|Z_c|^4.
\end{equation*}
Taking the sums for the separate image classes as the complex numbers
in this inequality proves the recombination bound, including the
mixed quadruples that use more than one class.
Also,
\[
 \mathcal E_K(W)\le e\,\mathcal G_K(W),
 \qquad \mathcal G_K(W)\ll K \mathcal E_1(W).
\]
For the second bound, bin pair sums by product character in intervals
of length $1/2$. Their populations satisfy
$\sum_jb_jb_{j+m}\leq\sum_jb_j^2$, and the total squared population
is at most $\mathcal E_1(W)$. The Gaussian sum over $m$ is $O(K)$.
Take $K=1+4\mathcal H$.

The comparison with the larger tolerance costs one factor proportional
to that tolerance, while recombination costs the cube of the number of
classes. Their product gives the asserted loss. The injectivity within
each class ensures that its cardinality is the number of original
points assigned to it, even though the full image is a multiset.
\end{proof}

\subsection{Selection from a complete value}

The complete amplitude is fixed before the product is selected.
Cancellation between different products prevents an upper bound for
an individual product from being inferred from a bound for their sum.
The following selection argument uses only lower thresholds.

\begin{proposition}\label{prop:ab-whole-witness}
Fix an original modulus $q$, a conductor divisor $r>1$, and $B=qT$,
where $T\geq1$ and $X^\kappa\leq B\leq X^{1/2}$ for a fixed
$\kappa>4/9$, as in Lemma~\ref{lem:ab-smooth-HB}.
Let $W$ be a unit separated set of nonprincipal characters modulo $q$
induced from conductor $r$, with $|t|\leq CT\mathcal H$ and
\[
 X^{-A}\leq V\leq|\Psi_u(t,\chi)|<2V.
\]
Suppose a product estimate is deduced from the condition on pair lengths,
the fourth moments of the original factors, and the mean square inequalities and inequalities for difference energy with two terms for all subset products.  Require it
uniformly for every product box, every final separated common set
$W_0$, and every admissible collection of polynomially bounded lower
thresholds, including the exponent perturbations described below.
Each test uses that set's own cardinality and energy and its actual
length and height.  In this formulation with lower thresholds, write it as
\begin{equation}\label{eq:ab-abstract-certificate}
 \left(\prod_iV_i\right)^2\mathcal E(W_0)
 \leq\mathscr A\left(X,B,\prod_iV_i\right)|W_0|.
\end{equation}
Then its consequence for complete values is
\begin{equation}\label{eq:ab-complete-certificate}
 V^2\mathcal E(W)
 \leq X^{O(\eta)+\varepsilon}\mathscr A_\eta(X,B,V)|W|.
\end{equation}
Here $\mathscr A_\eta$ is the supremum over
$O(\eta)+O(\varepsilon)$ perturbations of the length, height and
value exponents in the product estimate. The number of subset tests
is fixed independently of $X$.
\end{proposition}

\begin{proof}
The smooth identity \eqref{eq:ab-smooth-HB}, its bounded total
Mellin mass and its negligible tail give at each point of $W$ a
product $F_\nu$ and a shift $|\zeta|\leq\mathcal H$ such that
\[
 |F_\nu(t+\zeta,\chi)|\geq V L^{-C}/C.
\]
Choose one such product, then partition by the product and by the
dyadic magnitudes of its factors. Each factor has a fixed polynomial
upper bound. The product lower bound gives a fixed polynomial lower
bound for every factor on the selected set. There are therefore
only $L^C$ classes, with common coefficients and common original
intervals on each class.

Shifting a unit separated set by at most $\mathcal H$ creates at
most $O(\mathcal H+1)$ points per unit ordinate interval for each
character, counting multiplicity. Color the images into
$O(\mathcal H+1)$ separated injective classes. All characters and
all energy relations remain modulo $q$. Their new height is
$O(T\mathcal H)$, so the estimates for the mean, energy and fourth moments
use the actual scale $B_*=qC'(T\mathcal H+\mathcal H)
\asymp B\mathcal H$, with a fixed sufficiently large $C'$.
The original bound for pair lengths becomes
${\mathsf{M}}_i{\mathsf{M}}_j\gg B_*^{1/2}X^{-O(\eta)}$.
The stipulated estimate uniform under perturbations includes this change;
a bound requiring the unperturbed inequality for pair lengths at $B_*$
cannot be invoked without that robustness.

Let $U_i$ be the dyadic lower endpoints in one such class. Since
$\prod_iU_i\geq VL^{-C'}/C'$, decrease thresholds as needed until
\[
 \prod_iV_i=VL^{-C'}/C',\qquad 0<V_i\leq U_i.
\]
All lower bounds on values remain true after this operation.
Consequently the deficit corresponding to these product thresholds is
$1-\log_XV+o(1)$. The coefficient bounds imply
$\log_XV_i\leq\log_X{\mathsf{M}}_i+o(1)$, so the limiting individual
deficits are nonnegative. No upper bound on the actual selected
product has been used.

Every subset product is a polynomial on the same original variable
intervals, with a fixed support ratio and a fixed divisor bound.
The inequality for difference energy with two terms applies at every length.
There is thus no additional support partition or test depending on length to choose. Apply \eqref{eq:ab-abstract-certificate} on each
final image class, with that class's own cardinality and energy.
Lemma~\ref{lem:ab-energy-transfer} combines these estimates.
The class count, coloring and tolerance enlargement cost
$X^{O(\eta)+\varepsilon}$, while their cardinalities sum to
$|W|$. Absorb the dyadic threshold adjustment into the permitted exponent
perturbations. This proves \eqref{eq:ab-complete-certificate}.
\end{proof}
\subsection{The exact fourth Mellin kernel}

Put
\[
 K_\lambda(t)=\int_0^\infty
          w(y)e(\lambda y)y^{it}\,\frac{dy}{y}.
\]
Mellin inversion gives the exact formula
\begin{equation}
 \sum_n\Lambda(n)w(n/X)\chi(n)e(n\xi)
 =\frac1{2\pi}\int_{\mathbb R}
          K_{X\xi}(t)\mathcal D_u(t,\chi)\,dt.
 \label{eq:ab-mellin-identity}
\end{equation}
The insertion of $u$ is legitimate because it equals one
on $\operatorname{supp}w$.

To check the phases, insert the definition of the normalized complete
polynomial into the integral. Its summand contributes
\begin{equation*}
 \frac1{2\pi}\int_{\mathbb R}K_{X\xi}(t)(n/X)^{-it}\,dt
 =w(n/X)e(n\xi).
\end{equation*}
This is Mellin inversion at the scaled integer. Multiplying by the
auxiliary weight has no effect when the right side is nonzero. The
normalizing phase in the complete polynomial is what places the
scaled integer, rather than the integer itself, in this identity.

\begin{lemma}
\label{lem:ab-four-kernel}
For $T\ge2$, let $\psi$ be fixed smooth and supported in a
compact subinterval of $(0,\infty)$.
For either sign, all real $t_1,t_2,t_3,t_4$, and every fixed $A>0$,
\begin{align}
 &\left|\int_{\mathbb R}\psi(X|\xi|/T)
 K_{X\xi}(t_1)K_{X\xi}(t_2)
 \overline{K_{X\xi}(t_3)K_{X\xi}(t_4)}\,d\xi\right|\notag\\*
 &\qquad\ll_{A,w,\psi}\frac1{XT}
 (1+|t_1+t_2-t_3-t_4|)^{-A}.
 \label{eq:ab-four-kernel}
\end{align}
The integral is restricted to the chosen sign of $\xi$.
For $T=1$, the same conclusion holds with a smooth cutoff
supported in $|X\xi|\le C$.
One can restrict \eqref{eq:ab-mellin-identity} to $|t|\le CT\mathcal H$
with an arbitrarily small error with a negative power, uniformly on
these shells.
\end{lemma}

\begin{proof}
For $T\ge2$ put $X\xi=\pm Tz$ and change variables $v=zy$.
Then
\[
 K_{\pm Tz}(t)=z^{-it}J_z^\pm(t),\qquad
 J_z^\pm(t)=\int_0^\infty
        w(v/z)e(\pm Tv)v^{it}\,\frac{dv}{v}.
\]
Uniformly on the compact $z$ range, and for every fixed $j$,
\[
 |\partial_z^jJ_z^\pm(t)|\ll_{j,w,\psi}T^{-1/2}.
\]
For the phase $\pm2\pi Tv+t\log v$, use the second derivative bound
when $|t|\asymp T$, and integration by parts otherwise.
Differentiation in $z$ affects only the fixed smooth amplitude.

The product has $z$ derivatives $O(T^{-2})$.

There are four oscillatory integrals, each bounded by the inverse
square root of the shell height. Their product is therefore bounded
by its inverse square. The change of variables for the shell then gives
\begin{equation*}
 T^{-1/2}T^{-1/2}T^{-1/2}T^{-1/2}\frac{T}{X}
 =\frac1{XT}.
\end{equation*}
The four remaining phase factors multiply to a single power of the
new variable, with exponent equal to the energy discrepancy. This is
why the fourth Mellin kernel has exactly the discrepancy that appears
in the definition of energy.
Integrating the remaining phase $z^{-i(t_1+t_2-t_3-t_4)}$ by parts,
with $d\xi=(T/X)\,dz$, proves \eqref{eq:ab-four-kernel}.
For $T=1$, rapid decay of each $K_{X\xi}(t)$ and
$|t_1+t_2-t_3-t_4|\leq\sum|t_i|$ give the same conclusion.

For $|t|>CT\mathcal H$, nonstationary integration gives arbitrary
power decay. Combining this with $|\mathcal D_u|\ll X\log X$ and
choosing the decay order after $\eta$ proves the tail estimate,
including polynomially many modulus and character sums.
\end{proof}

\subsection{Gauss weights and the actual global mass}

For $q\asymp Q$, $r\mid q$, and a primitive character
$\chi^*\bmod r$, write $\chi_q$ for its induction to modulus $q$.
A nonzero induced Gauss sum requires
\[
 q=rs,\qquad s\text{ squarefree},\qquad(r,s)=1,
 \qquad |\tau_q(\chi_q)|^2=r.
\]
Only these characters contribute. In the following sums the principal
character is excluded.

For an amplitude level of the complete polynomial $V=X^v$, choose one
representative $t_{j,\chi_q}$ from each occupied unit interval of
\[
 \{t: |t|\le CT\mathcal H,\ V\le|\Psi_u(t,\chi_q)|<2V\}.
\]
Colour the intervals by parity; within one colour their representatives
are unit separated.

Representatives from two adjacent intervals could be arbitrarily close,
because each can lie near their shared boundary. Intervals of the same
parity have at least one full interval between their interiors, so their
representatives are separated by at least one. No independence
condition is imposed on their selection. The sampled inequalities
already control a point chosen after the polynomial is known.

Define the raw mass, summed over these finitely
many colours, by
\begin{equation}
 \mathcal M_v(Q,T)=
 \sum_{q\asymp Q}\ \sum_{\substack{r\mid q\\r>1}}
 \frac r{\varphi(q)}
 \sum_{\chi^*\bmod r}^{*}
 \sum_j|\Psi_u(t_{j,\chi_q},\chi_q)|^2.
 \label{eq:ab-raw-mass}
\end{equation}
Terms with zero Gauss sum are absent.

\begin{lemma}
\label{lem:ab-raw-mass}
If $Q^2T\ll X$, then
\begin{equation*}
 \mathcal M_v(Q,T)\ll X^{2+O(\eta)+\varepsilon}.
\end{equation*}
No representative is required to be chosen independently
of the polynomial, and no additional factor $Q$ occurs.
\end{lemma}

\begin{proof}
For the complete prime polynomial,
\[
 \Psi_u(t,\chi_q)-\Psi_u(t,\chi^*)=O_u(\log(2q))
\]
uniformly in real $t$.
The difference contains only powers of $p\mid q$, $p\nmid r$,
with weight $O_u(1+\log p)$ per prime. Its squared contribution is
$O(X^\varepsilon Q^2T\mathcal H)$.

The complete polynomial has the same coefficients before and after
passage to a conductor, except at prime powers removed by the larger
modulus. The earlier support argument bounds their total independently
of the ordinate. We may therefore compare the polynomials for all induced characters to a common
primitive polynomial before applying the large sieve. Keeping a
separate large sieve estimate for every induced character would repeat the
same primitive information.

Put $U=CT\mathcal H+1$ and apply \eqref{eq:ab-sampling} at height $U$.
The resulting integral for a primitive character has total weight from induced characters
\[
 \sum_{\substack{q=rs\asymp Q\\\mu(s)^2=1,\ (r,s)=1}}
       \frac r{\varphi(q)}
 \leq\frac r{\varphi(r)}
       \sum_{s\asymp Q/r}\frac1{\varphi(s)}
 \ll\frac r{\varphi(r)}.
\]
The dyadic sum of reciprocal totients is bounded, since
\[
 \frac{n}{\varphi(n)}=\sum_{d\mid n}\frac{\mu^2(d)}{\varphi(d)},
 \qquad \sum_d\frac{\mu^2(d)}{d\varphi(d)}<\infty.
\]
Retain $r/\varphi(r)$ for \eqref{eq:ab-hybrid-ls}.

For a contributing induced character, the coprimality of the conductor and
its complementary squarefree factor gives
\begin{equation*}
 \frac{r}{\varphi(q)}
 =\frac{r}{\varphi(r)\varphi(s)}.
\end{equation*}
Summing the reciprocal totients of the complementary factors costs
only the allowed small power. The remaining conductor weight is
already the one in the hybrid large sieve, so it introduces no
additional number of moduli.
The common coefficients $\Lambda(n)u(n/X)$ have square norm
$O(X\log^2X)$; differentiation of $\mathcal D_u$ inserts bounded
$\log(n/X)$.

The elementary bound for the square norm follows from the number of supported
integers and the bound $\Lambda(n)\le\log n$:
\begin{equation*}
 \sum_n|\Lambda(n)u(n/X)|^2\ll X\log^2X.
\end{equation*}
On the fixed support, the scaled logarithm introduced by
differentiation is bounded. Thus adaptive sampling costs only the same
coefficient square norm, together with logarithmic factors.
The large sieve uses $Q^2U+X$, and $Q^2U$ need not be at most $X$.
Since $U\ll T\mathcal H$, however,
\[
 Q^2U+X\ll\mathcal H(Q^2T+X).
\]
Thus Cauchy's inequality and \eqref{eq:ab-hybrid-ls} give
\[
 \mathcal M_v(Q,T)
 \ll X^{O(\eta)+\varepsilon}(Q^2T+X)X
 \ll X^{2+O(\eta)+\varepsilon}.
\]

\end{proof}

\begin{proposition}
\label{prop:ab-weighted-fourth-bridge}
Suppose that on every unit separated level of the complete polynomial
of magnitude $V$, with fixed $q,r$, one has
\[
 \frac{V^2\mathcal E(W)}{|W|}\le\mathscr A(X,QT,V).
\]
Then the nonprincipal contribution of the dyadic shell
$q\asymp Q$, $|X\xi|\asymp T$ has fourth moment at most
\begin{equation*}
 \frac{X^\varepsilon}{XQT}
 \sum_V\mathscr A(X,QT,V)\mathcal M_v(Q,T)
 +O_A(X^{-A}).
\end{equation*}
For $T=1$ use $|X\xi|\le C$.
The bound is valid for the whole induced character family,
with exactly the weights in \eqref{eq:ab-raw-mass}.
\end{proposition}

\begin{proof}
For $(a,q)=1$, character orthogonality gives
\[
 S_{\Lambda,w}(a/q+\xi)
 =\frac1{\varphi(q)}\sum_{\chi\bmod q}
       \tau_q(\overline\chi)\chi(a)
       \sum_n\Lambda(n)w(n/X)\chi(n)e(n\xi)
   +O((\log X)^C).
\]
Remove the principal summand.
Split the remaining sum by the conductor $r$, the dyadic value
level $V$, and the two colours of unit intervals.
H\"older costs at most $\tau(q)^3(\log X)^C\ll X^\varepsilon$.
Values below $X^{-A}$ and the Mellin tails are negligible by
Lemma~\ref{lem:ab-four-kernel}; the number of remaining value
levels is $O(\log X)$.

Fix one resulting part and insert \eqref{eq:ab-mellin-identity}.
The restrictions are independent of $a,\xi$. Summation over $a$ contributes
$\varphi(q)$ and imposes $\chi_1\chi_2=\chi_3\chi_4$, equivalently
for their inducing characters modulo $r$ by surjectivity of unit group
reduction. Each Gauss factor has modulus $\sqrt r$.

The arithmetic weight can now be counted explicitly. Four factors
from the character expansion give $\varphi(q)^{-4}$; summing the
reduced numerator supplies $\varphi(q)$ by orthogonality. The four
Gauss factors supply $r^2$. The net factor is therefore
\begin{equation*}
 \varphi(q)^{-4}\varphi(q)(\sqrt r)^4
 =\frac{r^2}{\varphi(q)^3}.
\end{equation*}
The equality of character products may be checked at the conductor
because reduction from the units modulo the original modulus onto
those modulo its divisor is surjective.

Use a nonnegative smooth majorant of the $\xi$ shell and
Lemma~\ref{lem:ab-four-kernel}.
Replacing the four ordinates by their chosen representatives of unit intervals changes their sum discrepancy by at most $4$,
and every integration domain for four intervals has measure at most $1$.

Within the chosen level, the product of the four polynomial absolute
values is at most a constant times the fourth power of the level
threshold. The shifts within the intervals change the decaying discrepancy
weight only by a bounded amount. We can consequently bound the
fourfold integral by this value factor times a discrepancy sum over
the representatives, where the estimate for positive energy applies.
Taking a fixed discrepancy exponent $A>1$, the decaying discrepancy
sum is $O(\mathcal E(W))$ by the pair bin
Cauchy argument in Lemma~\ref{lem:ab-energy-transfer}.
Combining these bounds, we find that this part contributes at most
\[
 \frac{X^\varepsilon r^2}{XT\varphi(q)^3}\,V^4\mathcal E(W).
\]
Since $r\le q$ and $q/\varphi(q)\ll X^\varepsilon$, this is at most
\[
 \frac{X^\varepsilon}{XqT}
 \left(\frac{V^2\mathcal E(W)}{|W|}\right)
 \left(\frac r{\varphi(q)}\,|W|V^2\right).
\]
The last parenthesis is comparable to the representatives' raw mass.

Indeed every representative on the level has squared absolute value
between $V^2$ and $4V^2$, so that
\begin{equation*}
 \frac r{\varphi(q)}|W|V^2
 \le \frac r{\varphi(q)}\sum_{(t,\chi_q)\in W}|\Psi_u(t,\chi_q)|^2
 \le 4\frac r{\varphi(q)}|W|V^2.
\end{equation*}
This is the point at which choosing the original set from a level of
the complete polynomial is used. The energy estimate is multiplied by
the actual weighted mass of those same representatives, and this mass
is the quantity already bounded globally.
Sum over $q,r,V$ and the colours, using $q\asymp Q$.
\end{proof}

\section{The actual mass of the prime polynomial and nearby zeros}

Write ${\mathscr{Z}}_{\mathrm{fam}}(\sigma;R,Y)$ for the count for the primitive family in
\eqref{eq:major-family-count}, with multiplicity and $|\Im\rho|\leq Y$.
To bound \eqref{eq:ab-raw-mass} by this count, associate each large
value of the prime polynomial with a nearby zero. The family density theorem
then gives Corollary~\ref{cor:HB-extended-level-mass}.

The quantity being estimated here is a weighted sum of squared values
of the original prime polynomial.  It is useful to distinguish this
quantity from the number of sample points.  A sample of size about
$X^v$ contributes about $X^{2v}$ to the mass, and its weight from Gauss sums
depends on the modulus from which the character is induced.  We shall
first explain why such a sample requires a nearby zero, then sum the
weights of all samples that can be associated with one zero.  This
order keeps the density theorem attached to primitive $L$-functions
while allowing the actual values to be sampled at different moduli of induction and different heights.

\begin{lemma}
\label{lem:smooth-zero-witness}
Fix $v_0>0$ and $C>0$.  Suppose $v_0\leq v\leq1$,
$q(2+|t|)\leq X^C$, and $\chi_q$ is induced from a primitive
nonprincipal character $\chi^*$ of conductor $r>1$.
If $\eta>0$ is sufficiently small in terms of $v_0$, then
\[
 |\Psi_u(t,\chi_q)|\geq X^v
\]
implies, for sufficiently large $X$, the existence of a zero
$\rho=\sigma+i\gamma$ of $L(s,\chi^*)$ satisfying
\[
 \sigma\geq v-2\eta,\qquad |\gamma-t|\leq X^\eta.
\]
If, in addition, all nontrivial zeros of the inducing primitive
characters have real part at most $\Theta$, where $1/2\leq\Theta<1$, then
\begin{equation}\label{eq:ab-actual-value-cap}
 |\Psi_u(t,\chi_q)|
 \ll_{u,\Theta}X^\Theta\log(q(2+|t|)).
\end{equation}
\end{lemma}

\begin{proof}
By Lemma~\ref{lem:ab-raw-mass}, the extra Euler factors cost
$O_u(\log(2q))$, so it suffices to treat $\chi^*$.
Mellin inversion for $-L'/L$ and a contour shift give
\begin{equation}\label{eq:ab-smooth-explicit-formula}
 \Psi_u(t,\chi^*)
 =-\sum_\rho X^{\rho-it}\widehat u(\rho-it)
   +O_u\bigl(\log(r(2+|t|))\bigr);
\end{equation}
see \cite{IK2004}. There is no principal pole, and the trivial zero
terms are bounded since $u$ is supported away from zero.

To see the role of the smoothing, the integrand before the contour
shift is
\[
 -\frac{L'}{L}(s+it,\chi^*)X^s\widehat u(s).
\]
A nontrivial zero $\rho$ of $L(s,\chi^*)$ becomes a pole at
$s=\rho-it$, with residue $-X^{\rho-it}\widehat u(\rho-it)$,
counted with its multiplicity.  This gives exactly the summands in
\eqref{eq:ab-smooth-explicit-formula}.  The smooth explicit formula
and the local zero count used below are standard analytic inputs
from the cited source.  Their use reduces the present argument to
estimating an absolutely convergent sum.

For the decay of the Mellin transform, put $x=e^y$ in its defining
integral.  The function $u(e^y)e^{\sigma y}$ is smooth and has compact
support in a fixed interval.  Repeated integration by parts against
$e^{i y\operatorname{Im}s}$ gives one inverse power of
$|\operatorname{Im}s|$ at each step, uniformly for $0\leq\sigma\leq1$.

The zero sum converges absolutely, since uniformly for
$0\leq\sigma\leq1$ and fixed $A>0$,
\[
 |\widehat u(\sigma+iy)|\ll_{u,A}(1+|y|)^{-A}.
\]
The standard local zero count is
$O(\log(r(2+|t|+|k|)))$ for a unit interval whose ordinates
satisfy $|\gamma-t-k|\leq1$.
Zeros with $\sigma<v-2\eta$ consequently contribute
$O(X^{v-2\eta}\log X)$, and zeros with
$|\gamma-t|>X^\eta$ contribute
\[
 O_{u,A}\bigl(X^{1-\eta(A-1)}\log X\bigr).
\]

Here the local count can be summed over unit intervals.  For example,
the part with $|\gamma-t|>X^\eta$ is bounded by a constant times
\[
 X\sum_{|k|>X^\eta-2}(1+|k|)^{-A}
       \log\bigl(r(2+|t|+|k|)\bigr).
\]
The polynomial bound on $q(2+|t|)$ controls the part of the logarithm
independent of $k$, and the tail of the convergent power sum gives
the stated estimate.  Choose $A$ so that
$1-\eta(A-1)<v_0-\eta$.  The two excluded portions are then
$o(X^v)$, as is the logarithmic error.  If the stated rectangle
contained no zero, these portions would exhaust the zero sum, which
would be incompatible with $|\Psi_u|\geq X^v$.

For sufficiently large $A$, absence of a zero in the stated rectangle
would contradict the lower bound. Under the hypothesis on the half plane,
absolute summation of \eqref{eq:ab-smooth-explicit-formula} gives
$O(X^\Theta\log(r(2+|t|)))$; restore the Euler factors to obtain
\eqref{eq:ab-actual-value-cap}.
\end{proof}

\begin{proposition}
\label{prop:weighted-zero-mass}
Let $Q,T\geq1$ be polynomially bounded in $X$, and let
$\mathcal M_v(Q,T)$ be the actual mass in \eqref{eq:ab-raw-mass},
with representatives of height at most $CT\mathcal H$ and
$V=X^v$, where $v_0\leq v\leq1$ and $v_0>0$ is fixed.
Then
\begin{equation}\label{eq:raw-mass-zero-count}
 \mathcal M_v(Q,T)
 \ll X^{2v+O(\eta)+\varepsilon}
 {\mathscr{Z}}_{\mathrm{fam}}(v-3\eta;2Q,C'T\mathcal H).
\end{equation}
The representatives may depend on the inducing character and the
modulus of the induced character.
\end{proposition}

\begin{proof}
Associate each representative with a zero from
Lemma~\ref{lem:smooth-zero-witness}. Unit separation gives at most
$O(X^\eta)$ representatives per zero and induced character; there are boundedly
many colours. All associated zeros have height at most $C'T\mathcal H$.

For a fixed zero $\rho=\sigma+i\gamma$, a representative associated with it must lie in $[\gamma-X^\eta,\gamma+X^\eta]$.
A set of ordinates separated by at least one has at most
$2X^\eta+2$ points in this interval.  This is a count for each
character and each modulus of induction; no common choice of ordinate
across induced characters is needed.  The fixed number of separated colours
changes only the implied constant.

Write $q=rs$.  Nonvanishing of the induced Gauss sum implies
that $s$ is squarefree and $(r,s)=1$.  Hence
\[
 |\tau_q(\chi_q)|^2=r,\qquad
 \frac r{\varphi(q)}
 =\frac r{\varphi(r)}\frac1{\varphi(s)}.
\]
For a fixed primitive character and one of its zeros, the weight
of all possible induced characters with $q\asymp Q$ is therefore at most
\[
 O(X^\eta)\frac r{\varphi(r)}
 \sum_{s\leq2Q/r}\frac1{\varphi(s)}
 \ll X^{O(\eta)+\varepsilon}.
\]

The estimate uses the elementary bounds
$r/\varphi(r)\ll_\varepsilon r^\varepsilon$ and
$\sum_{s\leq Y}1/\varphi(s)\ll_\varepsilon Y^\varepsilon$,
with the exponent losses chosen smaller than the final allowed loss.
Dropping the squarefree and coprimality restrictions only enlarges
this positive sum.  Thus the potentially many moduli $q=rs$ do not
produce a factor of order $Q$; their reciprocal totient weights are
what makes the passage to a primitive zero count possible.

Multiplying by the bound for squared values $O(X^{2v})$ and summing over
primitive zeros, we obtain the required mass bound. The cutoff
$v-3\eta$ absorbs fixed level constants, giving
\eqref{eq:raw-mass-zero-count} also for ordinates depending on the induced character.
\end{proof}

\begin{corollary}
\label{cor:HB-extended-level-mass}
Assume the hypotheses of Proposition~\ref{prop:weighted-zero-mass}.
Let $QT=X^\beta$ with $\beta\leq1/2$, put $D=1-v$ and
$d=1/2-\beta$, and suppose $129/167+3\eta\leq v\leq1$.
Then the mass of the complete prime polynomial satisfies
\begin{equation*}
 \mathcal M_v(Q,T)
 \ll X^{2-2D(1-2\beta)+O(\eta)+\varepsilon}
 =X^{2-4Dd+O(\eta)+\varepsilon}.
\end{equation*}
In limiting exponents, the saving from the unrestricted $X^2$
bound is $4Dd$ for $D\leq D_0:=38/167$.
\end{corollary}

\begin{proof}
Heath--Brown's Theorem~3 \cite{HeathBrown1979Density}, recorded
in \eqref{dens:density-three}, gives
\[
 {\mathscr{Z}}_{\mathrm{fam}}(\sigma;R,Y)
 \ll_\varepsilon(R^2Y^2)^{2(1-\sigma)+\varepsilon}
 \qquad(129/167\leq\sigma\leq1),
\]
uniformly for $R,Y\geq2$; bounded smaller parameters are covered
by increasing them to $2$. Apply the density estimate at the
actual zero height $Y$.
Substitute this estimate in \eqref{eq:raw-mass-zero-count}.
Since $Q^2(T\mathcal H)^2\ll(QT)^2X^{2\eta}$, we obtain the exponent
\[
 2v+4\beta(1-v)+O(\eta)+\varepsilon
 =2-2D(1-2\beta)+O(\eta)+\varepsilon.
\]

More explicitly, the factor counting zeros contributes
$(QT)^{4(1-v)+O(\eta)+\varepsilon}$, whereas the squared value
contributes $X^{2v}$.  Since $v=1-D$ and
$\beta=1/2-d$, the exponent without the losses is
\[
 2(1-D)+4\beta D
 =2-2D+2D-4Dd=2-4Dd.
\]
Both $D$ and $d$ must therefore be positive to obtain a saving by
this argument.  The first measures how far the amplitude is below
$X$, and the second measures how far $QT$ is below $\sqrt X$.
At the terminal scale $d=0$, this particular density estimate gives
no power saving, which explains the separate endpoint argument below.
The restriction $v-3\eta\geq129/167$ is the precise hypothesis
needed to apply the quoted density theorem.

\end{proof}

Within $O(\eta)$ of $v=129/167$, use Lemma~\ref{lem:ab-raw-mass};
the margin for the unrestricted mass in Section~\ref{sec:ar-subset-deficits} absorbs
this strip.

\begin{proposition}
\label{thm:ab-closed-bridge}
Let $P=X^\kappa$, where $4/9<\kappa<1/2$, and use the major arcs
of Proposition~\ref{prop:major-range}.
On each point outside them choose a Dirichlet approximation
\[
 \alpha=a/q+\xi,\qquad q\leq X^{1/2},\quad(a,q)=1,\quad
 |\xi|\leq(qX^{1/2})^{-1}.
\]
For the dyadic parameters $q\asymp Q$ and
$T\asymp\max(1,X|\xi|)$, one has, up to fixed constants,
\[
 P<qT\ll X^{1/2},\qquad Q^2T\ll X.
\]
For a packet estimate \eqref{eq:ab-abstract-certificate}
proved from the stated inputs, the complete nonprincipal minor
fourth moment satisfies
\begin{equation}\label{eq:ab-closed-fourth}
 \int_{\mathfrak m(P)}|S_{\rm np}(\alpha)|^4\,d\alpha
 \ll X^{O(\eta)+\varepsilon}
 \sum_{\substack{Q,T\ {\rm dyadic}\\ P\ll QT\ll X^{1/2}}}
 \frac1{XQT}\sum_V
       \mathscr A_\eta(X,QT,V)\mathcal M_v(Q,T).
\end{equation}
The mass bounds are Lemma~\ref{lem:ab-raw-mass} and, in its
stated amplitude range, Corollary~\ref{cor:HB-extended-level-mass}.
Under $(Q_\Theta)$, every nonprincipal level obeys
\[
 V\ll X^\Theta\log(q(2+T\mathcal H)).
\]

More explicitly, write $QT=X^\beta$ and suppose a complete level
has energy ratio at most $X^H=X^{2-\omega}$ and actual mass at
most $X^{2-m}$, with the allowed losses by small powers.
Its contribution to the fourth integral is then at most
\begin{equation}\label{eq:ab-block-exponent}
 X^{3-\beta-\omega-m+O(\eta)+\varepsilon}.
\end{equation}
After the positive major term is supplied, its exponent for the exceptional set in Bessel's inequality is $1-\beta-\omega-m$.
\end{proposition}

\begin{proof}
If $q\leq P$, exclusion from the major arcs gives
$X|\xi|>P/q$; if $q>P$, then $qT\geq q>P$ directly.
The Dirichlet approximation also gives $qT\ll X^{1/2}$,
and hence $Q^2T\ll Q(QT)\ll X$.

The two possibilities in the first step express the same geometric
fact.  Either the denominator itself is larger than the cutoff for the major arcs, or a small denominator is accompanied by an approximation
error larger than the width of the major arcs.  In the latter case
$T\geq X|\xi|>P/q$.  On the other hand, Dirichlet's inequality gives
$qX|\xi|\leq\sqrt X$, and $q\leq\sqrt X$ also bounds the alternative
$T\asymp1$.  Hence the combined parameter $qT$, rather than $q$
alone, measures the location of the sample on the minor arcs.

Use the measurable cells of Section~\ref{sec:major-principal} and
enlarge them to dyadic shells. With the energy ratio from
Proposition~\ref{prop:ab-whole-witness}, summing
Proposition~\ref{prop:ab-weighted-fourth-bridge} over these shells gives
\eqref{eq:ab-closed-fourth}. The perturbations in $\mathscr A_\eta$
include height and separation changes. For compact $\beta<1/2$,
choose these losses after fixing the range and the finite subset list.
At the top edge, Lemma~\ref{lem:ar-preliminary} gives the limiting
bound $H\leq2$. The same finite closed inequalities give
$H\leq2+o(1)$ under arbitrarily small height and coefficient
perturbations.
The bound on the actual polynomial is
\eqref{eq:ab-actual-value-cap}.
Finally, \eqref{eq:ab-block-exponent} is the identity
\[
 (XQT)^{-1}X^{2-\omega}X^{2-m}
 =X^{3-\beta-\omega-m}.
\]

There are three separate contributions in this multiplication.  The
factor $(XQT)^{-1}$ is the measure and normalization from the circle
integral.  The energy ratio contributes $X^{2-\omega}$, and the
actual mass of squared values contributes $X^{2-m}$.  Their savings
$\omega$ and $m$ add because the estimates are multiplied.  This
additivity is the reason for retaining the actual mass separately
from the energy argument.  Bessel's inequality subsequently divides
by $X^2$, the square of the major coefficient, giving the stated
exponent for the exceptional set.

\end{proof}

\begin{lemma}
Choose $0\leq w\leq1$ with $w=1$ on $[1/8,7/8]$, as in
Proposition~\ref{prop:major-range}, and suppose that its major
coefficient is $\gg X$ for every even $n\in[X/2,X]$.
If, for some fixed $\vartheta\geq0$, the nonprincipal and principal
minor fourth moments are both
$O(X^{2+\vartheta+\varepsilon})$ for every $\varepsilon>0$,
then the number of exceptional even integers in $[X/2,X]$ is
$O(X^{\vartheta+\varepsilon})$.
If the hypotheses hold on every sufficiently large dyadic scale,
then $E_+(X)\ll X^{\vartheta+\varepsilon}$.
\end{lemma}

\begin{proof}
The character expansion has a logarithmic correction from Euler factors,
and removing prime powers costs $O(X^{1/2+\varepsilon})$ pointwise.
Their fourth integrals are $O(X^{2+\varepsilon})$, so
$|A+B|^4\leq8(|A|^4+|B|^4)$ gives
\[
 \int_{\mathfrak m(P)}|S_w(\alpha)|^4\,d\alpha
 \ll X^{2+\vartheta+\varepsilon}.
\]
For exceptional even $n$, the full coefficient of $S_w^2$ vanishes,
so its minor coefficient has magnitude $\gg X$.

Indeed, let $c_n$ be the $n$th Fourier coefficient of
$\mathbf 1_{\mathfrak m(P)}S_w^2$.  Splitting the full coefficient
into its integrals over the major and minor arcs gives $c_n$ equal to the
negative of the major coefficient for an exceptional $n$.
The functions $e(n\alpha)$ are orthonormal on $[0,1]$, so Bessel's
inequality, restricted to the exceptional indices, states
\[
 \sum_{\substack{n\in[X/2,X]\\n\text{ is exceptional and even}}}|c_n|^2
 \leq \int_0^1
       |\mathbf 1_{\mathfrak m(P)}(\alpha)S_w(\alpha)^2|^2\,d\alpha.
\]
Each summand on the left is at least a fixed multiple of $X^2$.
This is the only place in this deduction where the positivity and
uniform size of the major coefficient enter.

Bessel's inequality gives
\begin{equation*}
 X^2\#\{n\in[X/2,X]:n\text{ is exceptional and even}\}
 \ll\int_{\mathfrak m(P)}|S_w(\alpha)|^4\,d\alpha
 \ll X^{2+\vartheta+\varepsilon}.
\end{equation*}
This proves the dyadic assertion.  Applying it with a smaller
positive exponent loss and summing the dyadic intervals proves
the stated bound for $E_+(X)$.
\end{proof}

\section{Subset deficits and the range of conductors and heights}
\label{sec:ar-subset-deficits}

We use inequalities of limiting exponents in this section.  The analytic
interpretation is given after the numerical argument.  Let
\begin{gather*}
 B=X^\beta,\quad {\mathsf{M}}_i=X^{a_i},\quad V_i=X^{u_i},\quad R=X^r,\\*
 \mathcal E(W)=X^e,\quad x_i=a_i-u_i,
 \qquad \frac49<\beta\leq\frac12.
\end{gather*}
The factors have common coefficients, a fixed divisor bound, and a common
set of sample points satisfying the lower thresholds for values.  The decomposition and the fourth moments of the original factors give
\begin{equation}\label{eq:ar-basic}
 \begin{gathered}
 a_i,x_i\geq0,\quad \sum_i a_i=1,\quad
 a_i+a_j\geq\frac\beta2\ (i\ne j),\\*
 0\leq r\leq\beta,\quad 2r\leq e\leq3r,
 \quad r+4u_i\leq\beta+2a_i.
 \end{gathered}
\end{equation}
There are at most eight factors: if there are at least two, summing the
inequalities for pair lengths gives $k\beta\leq4$.

For a subset $A$ of the factors write
$a_A=\sum_{i\in A}a_i$, $u_A=\sum_{i\in A}u_i$, and
$x_A=\sum_{i\in A}x_i$.  The two analytic estimates used for such subsets
are
\begin{align}
 r+2u_A&\leq a_A+\max(a_A,\beta),\notag\\
 e&\leq\max\{r+4x_A,\ \beta-2a_A+2r+4x_A\}.
 \label{eq:ar-de}
\end{align}
The first is the ordinary mean square estimate.  The second follows from
the estimate for difference energy with two terms
\[
 \mathcal E(W)V_A^4
 \ll_\epsilon X^\epsilon
       \bigl(R{\mathsf{M}}_A^4+BR^2{\mathsf{M}}_A^2\bigr).
\]
Thus it has no additional length or energy hypothesis.  Only subset
products of the original factors are required.

Set
\begin{gather*}
 D=\sum_i x_i,\qquad H=2\sum_i u_i+e-r,\qquad
 \omega=2-H,\qquad z=3r-e,\qquad d=\frac12-\beta,\\
 \mathscr L=\frac{e-r}{4}=\frac D2-\frac\omega4,
 \qquad C_0=\beta-r,\qquad K=\frac{C_0+z}{2},\\*
 \tau=\frac{\beta+r}{2},\qquad y_i=a_i-2x_i.
\end{gather*}
The identities and consequences used below are
\begin{align}
 &D=r+\frac{\omega-z}{2},\qquad
   Y:=\sum_i y_i=1-2D,\qquad z\geq0,\quad C_0\geq0,
   \quad y_i\leq\frac{C_0}{2},\label{eq:ar-identities}\\
 &a_A\geq\beta\ \Longrightarrow\ x_A\geq r/2,\notag\\
 &x_A+\frac12(\tau-a_A)_+\geq\mathscr L,
   \qquad x_A<\mathscr L\ \Longrightarrow\ y_A\leq K.
   \label{eq:ar-threshold}
\end{align}
In particular $a_A\geq\tau$ implies $x_A\geq\mathscr L$.

\begin{lemma}\label{lem:ar-preliminary}
One has $H\leq2$.  If $\beta<1/2$ and $\omega\leq d$, then
$5\leq k\leq8$, $0\leq\omega$, and $z\leq2\omega$.
For $k=5$ every pair and every complementary triple has deficit at least
$\mathscr L$; for $k=6$ every triple has this property.  For $k=8$ one has
the stronger inequality $z\leq\omega$.
\end{lemma}

\begin{proof}
For $k\leq4$, summing the fourth moments and using $e\leq3r$ gives
\[
 H\leq1+\frac{k\beta}{2}+\left(2-\frac k2\right)r
 \leq1+2\beta=2-2d.
\]
For $k=5$ take any pair and its triple complement; for $k=6$ take any
partition into two triples.  At least one group, say $C$, has length
exponent at least $1/2\geq\beta$, and hence $x_C\geq r/2\geq\mathscr L$.
If the other group $A$ had $x_A<\mathscr L$, then
\eqref{eq:ar-threshold} and the at most three original fourth moments in
$C$ would give
\[
 Y\leq K+\frac32C_0=2C_0+\frac z2.
\]
Substitution of \eqref{eq:ar-identities} yields
$\omega\geq2d+z/2$.  This is impossible when $\omega\leq d$ and $d>0$,
and is also impossible whenever $\omega<0$ and $d\geq0$.
Otherwise $x_A\geq\mathscr L$, and
$D\geq\mathscr L+r/2$ gives $z\leq2\omega$.

For $k=7$, every product of four factors has length exponent at least
$\beta$.  If some triple has deficit at least $\mathscr L$, its
complement consisting of four factors again gives $z\leq2\omega$.  If no triple has
this property, average $y_A\leq K$ over all triples to obtain
\[
 6Y\leq7(C_0+z).
\]
The seven original fourth moments also give
\[
 3r+2z\leq7\beta-2+2\omega.
\]
The first displayed inequality is
$5r+z+6\omega\geq6-7\beta$.  Eliminating $r$ gives
$\omega\geq2d+z/4$, again impossible in either of the ranges just
specified.

For $k=8$ all products of four factors are long.  Their complementary
mean square bounds give $D\geq r$, or $z\leq\omega$.
These arguments first rule out $\omega<0$ in every case, including
$d=0$, proving $H\leq2$.  With $d>0$ they give all the stated
consequences in the range $\omega\leq d$.
\end{proof}

\begin{proposition}\label{prop:ar-coarse}
If $4/9<\beta<1/2$ and $\omega\leq d$, then
\begin{equation*}
 D\leq\frac92\omega.
\end{equation*}
For five or six factors the stronger inequality $D\leq5\omega/2$ holds.
\end{proposition}

\begin{proof}
The preliminary lemma gives $0\leq\omega$ and excludes $k\leq4$.
For $k=5$, averaging the pair inequalities gives
$2D/5\geq\mathscr L$, hence $D\leq5\omega/2$.

For $k=6$, join $i,j$ when $x_i+x_j\geq\mathscr L$.
An edge has $y_i+y_j\leq C_0$ and a nonedge has
$y_i+y_j\leq K$.  The six original fourth moments imply
\[
 r+z\leq3\beta-1+\omega\leq\beta-d<\beta,
\]
so $K<C_0$.  If there were no two disjoint edges, all the edges would
form a star or a triangle, with at most five edges.  Summing over the
fifteen pairs would give
\[
 5Y\leq5C_0+10K=10C_0+5z,
\]
or $\omega\geq2d$, a contradiction.  Choose two disjoint edges and
let $v$ be the sum of the deficits at the other two vertices.
Then $v\leq D-2\mathscr L=\omega/2$.
The four triples containing those two vertices give
$D+3v\geq4\mathscr L=2D-\omega$.  Consequently
$D\leq\omega+3v\leq5\omega/2$.

For $k=7$, suppose for a contradiction that $D>9\omega/2$ and order
the deficits as $x_1\leq\cdots\leq x_7$.  The mean square of every
complement consisting of four factors gives, for each triple,
\[
 x_A\leq U:=D-r/2=\mathscr L+\omega/2-z/4.
\]
Consider the partition
\[
 A=\{1,2,3\},\qquad B=\{4,7\},\qquad C=\{5,6\}.
\]
First $x_A\leq3D/7<\mathscr L$.
Also $x_7\geq D/7>\omega/2$, so
$x_C\leq U-x_7<\mathscr L$.
Finally $D\leq4x_5+U$ gives $x_5\geq r/8$, while
\[
 r+2z=D-\omega/2+5z/2>4\omega.
\]
Thus $x_5>\omega/2-z/4$ and
$x_B\leq U-x_5<\mathscr L$.
Applying \eqref{eq:ar-threshold} to these three groups gives $Y\leq3K$,
or
\[
 D+\frac32\omega+\frac32z\geq2-3\beta.
\]
The seven original fourth moments give
\[
 3D+\frac72z\leq7\beta-2+\frac72\omega.
\]
Eliminating $D$ yields $16d\leq8\omega+z\leq10\omega\leq10d$,
a contradiction.

For $k=8$, if every triple had deficit below $\mathscr L$, averaging
$y_A\leq K$ over the triples would give
\[
 3Y\leq4(C_0+z),\qquad
 2r+z+3\omega\geq3-4\beta.
\]
The eight original fourth moments give
$2r+z\leq4\beta-1+\omega$.  Together these force
$\omega\geq2d$, a contradiction.  Some triple therefore has deficit
at least $\mathscr L$.  Averaging the five mean square bounds for subsets of four factors in its complement consisting of five factors gives complementary deficit at
least $5r/8$.  Hence
\[
 D\geq\mathscr L+5r/8,
 \qquad r+2z\leq4\omega,
 \qquad D\leq\frac92\omega-\frac52z\leq\frac92\omega.
\]
\end{proof}

\begin{proposition}
Suppose $4/9<\beta<1/2$, $D>0$, and
\begin{equation}\label{eq:ar-density-bad}
 \omega\leq(1-4D)d.
\end{equation}
Then
\begin{equation*}
 D\leq\frac52\omega.
\end{equation*}
For seven or eight factors one has $D\leq2\omega$.
\end{proposition}

\begin{proof}
The preceding results give $0\leq\omega\leq d$ and
$D\leq9\omega/2$.  In particular \eqref{eq:ar-density-bad} implies
$D<1/4$ and
\[
 D\leq\frac{(9/2)d}{1+18d},\qquad
 r\leq D+\omega/2
 \leq\frac{5d}{1+18d}\leq\frac5{36}<\frac\beta2,
\]
where $d<1/18$ has been used.  Every triple has length exponent at
least $3\beta/4\geq\tau$, and therefore every triple has deficit
at least $\mathscr L$.
For $k=5,6$ the result has already been proved.
For $k=8$, averaging over triples gives
$3D/8\geq\mathscr L$, hence $D\leq2\omega$.

For $k=7$, the triple average first gives $D\leq7\omega/2$.
Combining this with \eqref{eq:ar-density-bad} gives
\[
 r\leq D+\omega/2\leq\frac{4d}{1+14d},
 \qquad
 \tau\leq\frac14-\frac d2+\frac{2d}{1+14d}
 \leq\frac{41}{144}<\frac27.
\]
For the last inequality, the function preceding $41/144$ is increasing
on $0\leq d\leq1/18$ and has value $41/144$ at $d=1/18$.
The longest pair has length exponent at least $2/7$, so its deficit
is at least $\mathscr L$.  Averaging the triple inequalities in its
complement consisting of five factors gives complementary deficit at least
$5\mathscr L/3$.  Thus $D\geq8\mathscr L/3$, or $D\leq2\omega$.
\end{proof}

The preceding argument uses the fixed list of nonempty subset products
of at most eight factors.  All their intervals, support ratios and
coefficient orders are retained.  After discarding complete levels
$V<X^{-A}$ for one fixed sufficiently large $A$, put
$\beta_0=\log_X(QT)$ and $D_c=1-\log_XV$.
On a final separated set for a selected product, record its own cardinality
and energy and the factor exponents after adjusting the thresholds.  If $t$ denotes the
total exponent error, then
\[
 |\beta-\beta_0|+|D-D_c|\ll t,\qquad
 a_i+a_j\geq\beta/2-O(t),\qquad
 x_i\geq-O(t),\qquad \sum_i a_i=1+O(t).
\]
Here the analytic scale $X^\beta$ includes the actual shifted height;
the reconstruction denominator remains $XQT$.  All the inequalities
\eqref{eq:ar-basic}--\eqref{eq:ar-de} hold with errors $O(t)$,
uniformly over the fixed coefficient orders and support ratios.
The remaining coordinates are bounded because $D\leq1+A+O(t)$,
$0\leq r\leq\beta+O(t)$ and $2r-O(t)\leq e\leq3r+O(t)$.
Consequently every sequence of admissible records with $t\to0$
and $X\to\infty$ has a convergent subsequence with a fixed number of factors
satisfying the exact displayed inequalities, with
$\beta=\beta_0$ and $D=D_c$ at the limit.  Admissibility means that
all source estimates are valid at that record's $X$; no common
validity threshold as their precisions tend to zero is assumed.

We now combine the local inequalities with the mass of the complete prime polynomial.
Write $\delta=1-\Theta_0$, where
$\Theta_0=\max(\Theta,7/8)$.  The explicit formula gives $D\geq\delta$
for complete values, in limiting exponents.  The raw mass is at most
$X^{2+o(1)}$.  For
$D\leq D_0:=38/167$, with a fixed margin at this endpoint, the
density estimate for the primitive family gives mass at most
$X^{2-4Dd+o(1)}$.
The weighted Mellin reconstruction estimate therefore gives the block exponent for the exceptional set
\begin{equation*}\label{eq:ar-objective}
 \frac12+d-\omega-m,
 \qquad m=0\text{ for the raw mass},\qquad m=4Dd\text{ in the density range}.
\end{equation*}

\begin{proposition}\label{prop:ar-away}
Let $1/12<\delta\leq1/8$ and choose
\[
 \frac{1351}{3006}<\kappa<\frac9{20},\qquad
 d_\kappa=\frac12-\kappa,\qquad
 \frac{(5/2)d_\kappa}{1+10d_\kappa}<\delta.
\]
On every fixed compact interval
$\kappa\leq\beta\leq1/2-\eta_0$, the nonprincipal minor fourth moment
has a fixed power saving over $X^{5/2}$.
\end{proposition}

\begin{proof}
Put $C_\kappa=(9/2)d_\kappa<D_0$ and choose
$C_\kappa<D_b<D_0$.
If $D\geq D_b$, use the raw mass.  A block exponent at least $1/2$
would have $\omega\leq d$, so Proposition~\ref{prop:ar-coarse} would
give $D\leq(9/2)d_\kappa=C_\kappa<D_b$, a contradiction.
If $\delta\leq D\leq D_b$, use the density mass.  A block exponent
at least $1/2$ would satisfy \eqref{eq:ar-density-bad}, and therefore
\[
 D\leq\frac52(1-4D)d
 \leq\frac52(1-4D)d_\kappa,
 \qquad
 D\leq\frac{(5/2)d_\kappa}{1+10d_\kappa}<\delta,
\]
again a contradiction.
To make the uniformity explicit, fix the truncation $A$ above and
use the complete deficit $D_c$ to assign the density label when
$D_c\leq D_b$ and the raw label otherwise.  Each selected product
inherits its complete level's label.  In the limit with zero error retain
both labels at $D=D_b$.  For each factor count and each label, the
exact feasible records with
\[
 \kappa\leq\beta\leq1/2-\eta_0,\qquad
 \delta\leq D\leq1+A
\]
form a compact set, after imposing the corresponding closed condition
$D\leq D_b$ or $D\geq D_b$.  On each such set its labelled block
objective is continuous and, by the two contradictions just proved,
strictly less than $1/2$.  The finite union therefore has a fixed
positive gap; no continuity across the label switch is needed.
A sequence of admissible records with finite error violating a reserved
part of this gap would have an exact limiting record contradicting
that maximum.  Hence a fixed positive error tolerance and a threshold for large $X$ exist.  Choose the source precisions and the height and
losses from the search for zeros only after this tolerance, also reserving the
margin $D_0-D_b$ for the density cutoff, and then take $X$ beyond
all source thresholds.  The selection of complete values and weighted
Mellin reconstruction estimate now give the claimed power saving; their fixed
logarithmic summation costs are absorbed in the reserved gap.
\end{proof}

The existence of the stated $\kappa$ follows from
\[
 \lim_{\kappa\uparrow9/20}
 \frac{(5/2)(1/2-\kappa)}{1+10(1/2-\kappa)}=\frac1{12}.
\]
In particular all fixed $\Theta<11/12$ are admitted.

\section{Equality configurations at the terminal scale}

We first establish an exponent bound and identify all configurations
at which it can be sharp.  Only the mean square, the fourth moments of the original factors and the estimate for difference energy are needed.
For a subset of the factors we retain the original intervals and
apply the estimate with two terms directly to its product polynomial.
Consequently there are finitely many tests and no restriction on
the length of a test.

\begin{lemma}\label{lem:de-terminal-classification}
Suppose that $1\leq k\leq8$ and that real numbers
$a_i,u_i,r,e$ satisfy
\[
 a_i\geq0,\quad \sum_i a_i=1,\quad
 a_i+a_j\geq\frac14\ (i\ne j),\quad
 0\leq r\leq\frac12,\quad 2r\leq e\leq3r.
\]
Put $x_i=a_i-u_i\geq0$, $D=\sum_i x_i$, and suppose that $D>0$.
For every nonempty subset $A\subseteq\{1,\ldots,k\}$, use subscripts
to denote the corresponding sums, and suppose that
\begin{align}
 r+2u_A&\leq a_A+\max(a_A,1/2),\notag
       \\
 e&\leq\max\{r+4x_A,\ 1/2+2r-2a_A+4x_A\}.
       \label{eq:de-terminal-energy}
\end{align}
Suppose also that the fourth moments of the original factors give
\begin{equation*}
 r+4u_i\leq\frac12+2a_i\qquad(1\leq i\leq k).
\end{equation*}
Then
\begin{equation*}
 H:=2\sum_i u_i+e-r\leq2.
\end{equation*}
If $5\leq k\leq8$ and $H=2$, there are only the following possibilities:
\begin{enumerate}
\item $r=D=1/2$, $e=3r$ and $u_i=a_i/2$ for every $i$;
\item $k=5$, $0<r=D<1/2$, $e=3r$, and, after relabelling,
\[
 a_0=u_0=x_0=0,\qquad
 a_i=\frac14,\quad x_i=\frac r4
       \qquad(1\leq i\leq4).
\]
\end{enumerate}
\end{lemma}

\begin{proof}
Set
\[
 h=\frac12-r,\qquad z=3r-e,\qquad
 v_i=u_i-\frac{a_i}{2}.
\]
Thus $h,z\geq0$, and the original fourth moments state that
$v_i\leq h/4$.  We have the identities
\begin{equation*}
 H=2-2D+2r-z
   =2-2h+2\sum_i v_i-z.
\end{equation*}
The subset mean value gives
\begin{equation*}
 v_A\leq\frac h2\quad(a_A\leq1/2),\qquad
 x_A\geq\frac r2\quad(a_A\geq1/2).
\end{equation*}
In the same notation, difference energy is
\begin{equation}\label{eq:de-terminal-v-energy}
 v_A\leq
 \max\left\{\frac{a_A-r}{2}+\frac z4,
                  \frac{h+z}{4}\right\}.
\end{equation}
All these inequalities are closed, including at equal lengths.

If $r=1/2$, the fourth moments give $u_i\leq a_i/2$, so $H\leq2$.
Equality requires $e=3r$ and equality in each fourth moment;
it gives the first configuration in the statement.
We may therefore assume that $h>0$.
For $k\leq4$, the fourth moments and $e\leq3r$ immediately give
\[
 H\leq1+\frac{k h}{2}+2r
      =2+\left(\frac k2-2\right)h\leq2.
\]

We record a useful consequence of a partition $A\cup C$ with
$a_A\leq1/2$, $a_C\geq1/2$, and at most three factors in $C$.
If the first term in \eqref{eq:de-terminal-energy} is the larger one,
combine it with $x_C\geq r/2$ to obtain
$H\leq2-z/2$.  If the second term is the larger one, then
\[
 2u_A\leq\frac14+r+a_A-\frac e2,
 \qquad
 2u_C\leq a_C+\frac34-\frac{3r}{2}.
\]
Here the second bound follows by summing at most three original
fourth moments; adding an unused nonnegative allowance $h/2$ if
necessary covers fewer factors.  Their sum again gives
\begin{equation}\label{eq:de-terminal-partition}
 H\leq2-\frac z2.
\end{equation}
For $k=5$, take a shortest pair as $A$; for $k=6$, take a shortest
triple.  These choices prove the baseline for both values of $k$.
In either equality case we have $z=0$, $D=r>0$, and
\begin{equation}\label{eq:de-terminal-v-total}
 \sum_i v_i=h.
\end{equation}

Suppose first that $k=5$ and $H=2$.
The average deficit of a pair is $2r/5$, so there is a pair $A$
with $x_A<r/2$.  Its mean value implies $a_A<1/2$, and its
first term for difference energy cannot support $e=3r$.
Thus $v_A\leq h/4$.  The other three fourth moments and
\eqref{eq:de-terminal-v-total} force
\[
 v_A=h/4,\qquad v_i=h/4\quad(i\in C).
\]
For any pair of indices in $C$ the value $v_i+v_j=h/2$
exceeds the second term in \eqref{eq:de-terminal-v-energy}.
The first term therefore gives $a_i+a_j\geq1/2$.
Summing the three pair inequalities gives $a_C\geq3/4$.
Since the original condition on pair separation gives $a_A\geq1/4$,
equality holds throughout:
\[
 a_A=1/4,\qquad a_i=1/4\quad(i\in C).
\]
Each of the two values in $A$ is nonnegative, because their sum is
$h/4$ and neither exceeds $h/4$.
If both were positive, pairing either index with any index of $C$
would force the first energy term, and hence $x_i\geq r/4$
for both indices in $A$.  This contradicts
\[
 x_A=a_A/2-v_A=1/8-h/4=r/4.
\]
One of the two values is therefore zero and the other is $h/4$.
Pairing the latter with an index in $C$ gives its length at least
$1/4$.  Its length is consequently $1/4$, while the index with zero value has length zero.  These are precisely the four factors of quarter length and the factor of zero length in the second configuration.

Now suppose that $k=6$ and $H=2$.
For any partition into two triples, call a shorter triple $A$
and its complement $C$.  If the second energy term can support
the value, the argument for \eqref{eq:de-terminal-partition}
forces all three $v_i$, $i\in C$, to equal $h/4$.
Their three pair tests then give $a_C\geq3/4$ just as above.
Pair separation within $A$ gives $a_A\geq3/8$, a contradiction.
Thus the first term supports every such triple test.
The long complementary mean and $D=r$ give
\[
 x_A=x_C=r/2
\]
for every triple partition.  Comparing triples gives $x_i=r/6$
for every $i$.  A pair of maximum length has $a_A\geq1/3$ and
$x_A=r/3$.  Its energy estimate would imply
\[
 3r\leq\max\left\{\frac{7r}{3},
                         \frac{10r}{3}-\frac16\right\}<3r,
\]
because $0<r<1/2$.  No such equality configuration exists.

For $k=7$ we prove the baseline and exclude sparse equality at
the same time.  Assume that $H\geq2$ and write $\omega=2-H\leq0$.
Then
\[
 \sum_i v_i=h+\frac{z-\omega}{2}.
\]
Every triple has length at most $1/2$: its complement consisting of four factors
has length at least $1/2$, by summing its six pair inequalities.
Consequently every triple $D$ satisfies $v_D\leq h/2$.
Fix a triple $A$ and let $C$ be its complement consisting of four factors.
Summing this last inequality over the four triples within $C$ gives
$3v_C\leq2h$.  Hence
\[
 v_A\geq\frac h3+\frac{z-\omega}{2}
       >\frac{h+z}{4},
\]
since $h>0$, $z\geq0$ and $\omega\leq0$.
The second energy term is therefore impossible for every triple.
The first term gives $x_A\geq r/2-z/4$, while the long
complementary mean gives $x_C\geq r/2$.
These inequalities imply $H\leq2-z/2$, hence
$z\leq2\omega\leq0$.
We obtain $z=\omega=0$, $D=r>0$, and $x_A=D/2$ for every triple.
But the average triple deficit is $3D/7$, contradicting $D>0$.
This proves that $H<2$ whenever $k=7$ and $r<1/2$.

Finally let $k=8$.
Summing the twenty-eight pair inequalities forces every pair
length to be $1/4$, and hence every $a_i=1/8$.
Mean value estimates for two complementary products of four factors give $D\geq r$,
so $H\leq2-z\leq2$.
If $H=2$, then $z=0$ and $D=r>0$.
Every deficit of four factors is $r/2$, and comparison gives
$x_i=r/8$ for every $i$.
The energy test for a triple of length $3/8$ and deficit $3r/8$
would give
\[
 3r\leq\max\left\{\frac{5r}{2},
                         \frac{7r}{2}-\frac14\right\}<3r
 \qquad(0<r<1/2).
\]
This excludes the remaining sparse equality case and completes
the proof.
\end{proof}

The normalized inequalities in this lemma follow directly from the
finite analytic estimates: coefficient maxima are bounded by an
arbitrarily small power, whereas the convolution square norms retain
fixed logarithmic powers.  At $B=X^{1/2+o(1)}$ the raw M\"obius factors from the identity with four stages have length at most $X^{1/4+o(1)}$, so their
$B+{\mathsf{M}}_i^2$ corrections give exactly the fourth moment inequalities for the original factors used above.

For the terminal compactness argument with $5\leq k\leq8$,
retain all subset tests.
Their number is at most $2^8-1$, their support ratios are fixed,
and the inequalities for the mean and difference energy are continuous in
the exponent data.  A sequence outside fixed neighborhoods of the
two displayed equality configurations, with $B=X^{1/2+o(1)}$
and $H\to2$, would have a limiting subsequence contradicting the
lemma.  Since the bound on the adjusted product threshold gives $D\geq\delta>0$,
this limit stays in the stated range.  Thus there is a fixed
positive gap on the complement of those neighborhoods.
Choose the width of the terminal strip, and then all coefficient allowances by small powers, inside this gap.
The exact sparse configuration already has one factor of zero length
and four factors of quarter length; no tests using high pure powers or
conditional length boundaries are needed.

\section{A terminal reconstruction estimate with fixed logarithmic loss on the original modulus}
\label{sec:endpoint-gauss-bridge}

The terminal fourth moment requires reconstruction on the original
modulus, with the Gauss weights and additive height enlargement
retained in the mass estimate.

\subsection{Gauss weights and enlarged heights}

There are two different operations here. Character orthogonality reconstructs
an additive sum at a rational point with denominator $q$; this operation
uses characters modulo that same $q$. Estimates for mean values, on the other
hand, are most naturally summed over primitive characters. We keep the
Gauss factor during this second operation, so that passing to the inducing primitive
character does not count all of its induced characters with equal weight.

We shall first bound a sum of squares of complete prime values, then
multiply this mass by an energy bound to obtain a fourth moment. The mass
measures the total size of the selected values. The energy measures how
often four selected points given by a character and an ordinate satisfy the character
relation and an approximate additive relation between their ordinates.
The factorization below keeps these two tasks separate.

Put $L=\log X$.  Choose all cutoffs in logarithmic variables to be
fixed Gevrey functions of order $1<s\leq2$.  Let
\[
 \Psi_u(t,\chi)=\sum_n\Lambda(n)u(n/X)\chi(n)n^{-it},
 \qquad \mathcal D_u(t,\chi)=X^{it}\Psi_u(t,\chi).
\]
The ordinary cutoff $u$ equals one on the support of the additive
cutoff $w$.  The principal character is treated separately throughout.

Multiplication by $X^{it}$ does not change the absolute value. It writes
the oscillating factor in the useful form $(n/X)^{-it}$. Since $n/X$ stays
in a fixed compact interval on the support of $u$, differentiation with
respect to $t$ now introduces a bounded logarithm. This is what allows
sampling at points chosen after seeing the values.

For each nonprincipal character modulo the original modulus $q$, set
\begin{equation}\label{eq:endpoint-G-definition}
 \gamma_q(\chi)=\frac{|\tau_q(\overline\chi)|}{\sqrt q},
 \qquad
 G_q(t,\chi)=\gamma_q(\chi)\mathcal D_u(t,\chi).
\end{equation}
Characters with vanishing Gauss sum contribute nothing.  For the
others, the induced formula for Gauss sums gives $0<\gamma_q(\chi)\leq1$,
and hence $|G_q|\leq|\Psi_u|$.  Thus the selected products for
$G_q$ have common coefficients.

In particular, a large value of $G_q$ is also a large value of the
unweighted polynomial $\Psi_u$. We apply the product decomposition to
that polynomial, whose coefficients are independent of the sampled
character and ordinate. The factor $\gamma_q(\chi)$ is retained when we
return to the mass and to additive reconstruction; it need not be absorbed
into the coefficients of a selected product.

Fix a shell $q\asymp Q$, $|X\xi|\asymp T$, with $T\geq1$;
the notation $T=1$ includes the central interval $|X\xi|\leq C$.
Write $B_0=QT$.  The shells under consideration satisfy
\[
 Q^2T\ll X,\qquad B_0\ll X^{1/2}.
\]
Gevrey Mellin decay permits the original Mellin integral to be
truncated at height $C_0T+H$, where $H=C_A L^s$, with an error
smaller than any prescribed fixed negative power of $X$.
Selecting a product changes each ordinate by at most $H$.
To include these shifts, choose the common enlarged height
\begin{equation}\label{eq:endpoint-enlarged-height}
 U=C_1(T+H),\qquad B_q=qU.
\end{equation}
Thus $B_q\ll B_0L^s\ll X^{1/2}L^2$.

The physical shell is still described by $B_0=QT$: its points have
$|X\xi|$ of size $T$. The analytic estimates for product values use the
larger scale $B_q=qU$, because their ordinates have been shifted. When
$T$ is much smaller than $H$, these two scales differ by a logarithmic
factor. We therefore postpone summing the shells until the dependence on
$T+H$ has been written out.

\subsection{One representative for each character and time interval}

Partition the truncated height interval into unit intervals.  For each
character and occupied interval, choose a point where $|G_q|$ reaches
its maximum, using closed intervals and disregarding their overlaps of measure zero.  Assign the cell to the dyadic level $V$ satisfying
\[
 V\leq \max_{\text{cell}}|G_q|<2V.
\]
The values below $X^{-A}$ contribute a negligible error, and the
polynomial upper bounds leave $O_A(L)$ levels.  Coloring the cells
by parity, we obtain unit separation for each character across all
levels, with one representative per cell.

For clarity, a cell here consists of one fixed character together with a
unit time interval. Two different characters give different cells even
when their intervals coincide. For one character, intervals of the same
parity are separated by an intervening unit interval, so their chosen
points are at least one unit apart. Each cell is assigned to just one
value level. Consequently, joining all levels of a fixed color does not
create additional representatives in that cell. This last observation
will prevent a further factor of $L$ in the mass estimate.

Let $W_{q,V,c}$ be one such level and color, and put $R_{q,V,c}=|W_{q,V,c}|$.
Restrictions may retain entire assigned cells or their intersection with
$|G_q|>Y$; each occupied high cell has maximum above $Y$.
They may depend on $q,\chi$, the shell and time, but not on $a$ or
$\xi$.  On an assigned cell $|G_q|<2V$, and at its representative
$|G_q|\geq V$.

A restriction specifies which points of a cell are kept in the Mellin integral.
Its independence of $a$ allows us to use character orthogonality after
expanding the fourth power. Its independence of $\xi$ allows the same
selection of ordinates to be used throughout the shell. The cell maximum
provides both an upper bound for integration inside the cell and a lower
bound for the discrete mass at its representative.

\subsection{The actual mass, with all conductor weights retained}

Define
\begin{equation*}\label{eq:endpoint-actual-mass}
 \mathcal M_G(Q,U)=
 \sum_{q\asymp Q}\frac q{\varphi(q)}
 \sum_{V,c}\sum_{(t,\chi)\in W_{q,V,c}}|G_q(t,\chi)|^2.
\end{equation*}
By \eqref{eq:endpoint-G-definition}, its summands are exactly
\[
 \frac{|\tau_q(\overline\chi)|^2}{\varphi(q)}
 |\Psi_u(t,\chi)|^2.
\]
\begin{lemma}
\label{lem:endpoint-raw-mass}
Uniformly for the adaptively chosen representatives above,
\begin{equation*}
 \mathcal M_G(Q,U)\ll (Q^2U+X)XL\ll X^2L^{s+1}.
\end{equation*}
The last inequality uses \eqref{eq:endpoint-enlarged-height} and
$Q^2T\ll X$.
\end{lemma}

\begin{proof}
There are only two colors, and within each color the representatives
remain unit separated after all levels are joined.  If a character
modulo $q$ is induced from a primitive nonprincipal character of
conductor $r$, a nonzero Gauss sum requires
\[
 q=rs,\qquad \mu(s)^2=1,\qquad(r,s)=1,
 \qquad |\tau_q(\overline\chi)|^2=r.
\]
Replacing the complete prime polynomial by that of its inducing primitive character
changes it by $O_u(\log(2q))$.  Normalize by $X^{it}$ before applying
Sobolev sampling; its derivative then inserts the bounded factor
$-i\log(n/X)$.

Here is the sampling step explicitly. If $F$ is the normalized smooth
polynomial for one primitive character and the points $t_j$ are separated
by one, the elementary Sobolev estimate in one dimension gives
\begin{equation*}
 \sum_j |F(t_j)|^2
 \ll \int_{-CU-C}^{CU+C}\bigl(|F(t)|^2+|F'(t)|^2\bigr)\,dt.
\end{equation*}
We apply the estimate on bounded intervals around the selected points;
these intervals have bounded overlap. It holds for every choice of the
points, so it also holds for the maxima chosen above. For the derivative,
we use the coefficients multiplied by $-i\log(n/X)$, which have the same
bound for the sum of squares as the original coefficients.

After sampling has replaced the selected points by a full height
integral, the integrand no longer depends on the induction modulus $q$.
For fixed $r$ the weights of induced characters satisfy
\[
 \sum_{\substack{q=rs\asymp Q\\\mu(s)^2=1,\ (r,s)=1}}
 \frac r{\varphi(q)}
 \leq\frac r{\varphi(r)}
       \sum_{Q/r\leq s<2Q/r}\frac1{\varphi(s)}
 \ll\frac r{\varphi(r)}.
\]
The last dyadic sum of reciprocal totients is bounded absolutely:
use $n/\varphi(n)=\sum_{d\mid n}\mu^2(d)/\varphi(d)$ and
$\sum_d\mu^2(d)/(d\varphi(d))<\infty$.

The factor $r/\varphi(q)$ is the reason this sum over induced characters is bounded.
For $q=rs$ with $(r,s)=1$, we have
\begin{equation*}
 \frac{r}{\varphi(q)}=\frac{r}{\varphi(r)}\frac1{\varphi(s)}.
\end{equation*}
Once sampling has produced a common integral, the summation over $s$
therefore contributes only a constant. Replacing the Gauss weight by its
upper bound before this step would replace this useful sum of reciprocal totients by a count of induced characters. The continuous hybrid large sieve is the
analytic mean value theorem used after this reduction.

The continuous hybrid multiplicative large sieve therefore gives
$(Q^2U+X)\sum_n|a_n|^2$ with its usual weight $r/\varphi(r)$.
For the prime coefficients and their normalized derivative,
$\sum|a_n|^2\ll XL$: indeed $\Lambda(n)^2\ll L\Lambda(n)$
on the fixed support, and the elementary Chebyshev bound gives
$\sum_{n\leq CX}\Lambda(n)\ll X$.  This proves the first bound for the main
term.  The corrections from Euler factors contribute $O(Q^2UL^2)$:
the exact identity
$\sum_{\chi\bmod q}|\tau_q(\chi)|^2=\varphi(q)^2$
sums their weights.  This error is absorbed by the displayed bound.
Finally $U\ll TL^s$ gives the second inequality.

Indeed, $T\geq1$ and $H\asymp L^s$ imply $T+H\ll TL^s$, and hence
\begin{equation*}
 (Q^2U+X)XL
 \ll (Q^2TL^s+X)XL
 \ll X^2L^{s+1}.
\end{equation*}
This bound is for all selected levels together. No count of value levels has
been inserted into the sampled square sum.

\end{proof}

\subsection{The weighted fourth power without conductor splitting}

For nonvanishing Gauss sums write
$\epsilon_q(\chi)=\tau_q(\overline\chi)/|\tau_q(\overline\chi)|$.
The complete nonprincipal contribution is, up to the harmless
correction contributed by prime powers with $(n,q)>1$,
\begin{equation}\label{eq:endpoint-G-reconstruction}
 S_{q,\mathrm{np}}(a,\xi)=
 \frac{\sqrt q}{2\pi\varphi(q)}
 \sum_{\substack{\chi\bmod q\\\chi\ne\chi_0}}
 \epsilon_q(\chi)\chi(a)
 \int_{\mathbb R}K_{X\xi}(t)G_q(t,\chi)\,dt.
\end{equation}
All character relations and all energy estimates in this subsection
are taken modulo $q$ itself.

The factor $\epsilon_q(\chi)$ records only a complex phase. Its modulus
is one, whereas the size of the Gauss sum has already been included in
$G_q$. Keeping this distinction makes the expansion of the fourth power
transparent: phases can be discarded in an upper bound, but the sizes
must remain in the mass.

\begin{proposition}
\label{prop:endpoint-G-bridge}
Suppose every retained set $W_{q,V,c}$ satisfies
\[
 \frac{V^2\mathcal E(W_{q,V,c})}{R_{q,V,c}}
 \leq \mathscr A_q(V).
\]
Then the fourth moment of the corresponding part of the shell is
\begin{equation}\label{eq:endpoint-G-fourth}
 \ll\frac{L^3(\log\log(3X))^2}{XB_0}
 \left(\sup_{q,V,c}\mathscr A_q(V)\right)
 \mathcal M_G(Q,U)+O_A(X^{-A}).
\end{equation}
The same assertion holds after retaining any selection of the levels.
\end{proposition}

\begin{proof}
Split \eqref{eq:endpoint-G-reconstruction} into the $O(L)$ dyadic
values and two colors.  H\"older's inequality contributes $O(L^3)$.

More precisely, for $J_0$ parts we use
\begin{equation*}
 \left|\sum_{j=1}^{J_0} z_j\right|^4
 \leq J_0^3\sum_{j=1}^{J_0}|z_j|^4,
 \qquad J_0\ll L.
\end{equation*}
This is the only cost of the union over value levels in the reconstruction step.

For each part, we expand the fourth power and sum over $(a,q)=1$.
Orthogonality supplies $\varphi(q)$ and imposes
$\chi_1\chi_2=\chi_3\chi_4$ modulo $q$.  The Gauss phases have
absolute value one, while the four factors $\sqrt q$ give $q^2$.

The orthogonality identity being used is
\begin{equation*}
 \sum_{a\bmod q}^{*}
 \chi_1(a)\chi_2(a)\overline{\chi_3(a)\chi_4(a)}
 =\begin{cases}
 \varphi(q),&\chi_1\chi_2=\chi_3\chi_4,\\
 0,&\text{otherwise}.
 \end{cases}
\end{equation*}
Thus the four denominators $\varphi(q)$ in reconstruction leave
$\varphi(q)^{-3}$ after the sum over $a$. The remaining integral depends
on the time discrepancy $t_1+t_2-t_3-t_4$, exactly the quantity counted
by the energy.

The fourfold Mellin kernel bound supplies
\[
 \frac{C_A}{XT}
 (1+|t_1+t_2-t_3-t_4|)^{-A}.
\]
Replace each time by its chosen cell representative.  The discrepancy
changes by at most four, each cell integration domain has measure
at most one, and the positive pair bin energy estimate bounds the
resulting decaying discrepancy sum by $O(\mathcal E(W_{q,V,c}))$.

At each of the four integration points the amplitude is at most $2V$.
Replacing the times by representatives therefore leaves a factor
$O(V^4)$ and a sum over representative quadruples. The positive pair bin
estimate compares the decaying discrepancy sum with sharp energy, with a
constant unaffected by the bounded displacement. It applies to the
representatives chosen from the actual values.

This part is consequently at most
\[
 \frac{Cq^2}{XT\varphi(q)^3}V^4\mathcal E(W_{q,V,c}).
\]
To combine this with the actual mass, write
\begin{equation*}
 \frac{q^2 V^4\mathcal E}{XT\varphi(q)^3}
 =\frac{(q/\varphi(q))^2}{XqT}
   \left(\frac{V^2\mathcal E}{R}\right)
   \left(\frac q{\varphi(q)}RV^2\right).
\end{equation*}
The final factor is at most the actual mass of its representatives.

Indeed, every representative in the class has squared magnitude at least
$V^2$, so
\begin{equation*}
 \frac q{\varphi(q)}RV^2
 \leq\frac q{\varphi(q)}
        \sum_{(t,\chi)\in W_{q,V,c}}|G_q(t,\chi)|^2.
\end{equation*}
The middle factor in the preceding factorization is controlled by the
assumed energy estimate. Only the first factor remains outside the mass
sum; it is here that the two powers of $q/\varphi(q)$ are bounded.

Use $q/\varphi(q)\ll\log\log(3X)$, sum the parts, and replace
$qT$ by its comparable shell scale $B_0$.  This proves \eqref{eq:endpoint-G-fourth}.
\end{proof}

\subsection{The part with low amplitude also permits the normalization}

\begin{lemma}
\label{lem:endpoint-G-low-cap}
Retain any measurable restrictions for which $|G_q(t,\chi)|\leq Y$.
In particular this permits the pointwise cutoff on the entire real
time axis, including portions with low values in cells that otherwise have high values.
Then their reconstructed contribution to the shell satisfies
\begin{equation*}
 \sum_{q\asymp Q}\sum_{a\bmod q}^{*}
 \int_{|X\xi|\asymp T}|S_{q,\mathrm{lo}}(a,\xi)|^4\,d\xi
 \ll\frac{B_0^3Y^4}{X}.
\end{equation*}
For $T=1$ use the central interval as before.  In particular the choice
\[
 Y=\sqrt X\left(\frac{X}{B_0^2}\right)^{3/8}
   =X^{7/8}B_0^{-3/4}
\]
gives a bound $O(X^{5/2})$ for this part of one shell.
 If the same bound $Y$ is used for every rational cell with
 $q\leq\sqrt X$ and $1\leq T\ll\sqrt X/q$, their full contribution
 is bounded more precisely by
 \begin{equation}\label{eq:endpoint-G-global-low}
  \int |S_{\mathrm{lo}}(\alpha)|^4\,d\alpha
  \ll \sqrt X\,Y^4L.
 \end{equation}
 Thus $Y=\sqrt X L^K$, with any fixed $K\geq0$, costs only
 $O(X^{5/2}L^{4K+1})$.
\end{lemma}

\begin{proof}
Stationary phase and nonstationary integration by parts give
$\int_{\mathbb R}|K_{X\xi}(t)|\,dt\ll\sqrt T$ uniformly on
the shell.  The formula \eqref{eq:endpoint-G-reconstruction} and
the fact that there are at most $\varphi(q)$ characters give
\[
 \sup_a|S_{q,\mathrm{lo}}(a,\xi)|^2\ll qTY^2.
\]
Character orthogonality, followed by the same Mellin $L^1$ bound,
also gives
\[
 \sum_{a\bmod q}^{*}|S_{q,\mathrm{lo}}(a,\xi)|^2
 \leq \frac{Cq}{\varphi(q)}
 \sum_{\chi\bmod q}
 \left|\int K_{X\xi}(t)G_q(t,\chi)\mathbf1_{\mathrm{lo}}\,dt\right|^2
 \ll qTY^2.
\]
Multiplying these two bounds, integrating over an interval of length
$O(T/X)$ and summing over $q\asymp Q$, we obtain $O(Q^3T^3Y^4/X)$.

The estimate for the fourth power uses the elementary inequality
\begin{equation*}
 \sum_{a\bmod q}^{*}|S_{q,\mathrm{lo}}(a,\xi)|^4
 \leq \left(\sup_a|S_{q,\mathrm{lo}}(a,\xi)|^2\right)
        \sum_{a\bmod q}^{*}|S_{q,\mathrm{lo}}(a,\xi)|^2.
\end{equation*}
The two square bounds give $q^2T^2Y^4$. The width $T/X$ contributes
the third power of $T$, and summing $q^2$ over $q\asymp Q$ contributes
$Q^3$. With the displayed choice of $Y$ depending on the shell, we have
$Y^4=X^{7/2}B_0^{-3}$, so $B_0^3Y^4/X=X^{5/2}$.

 For the final claim, sum geometric heights $T\ll\sqrt X/q$ first:
 \[
  \frac{q^2Y^4}{X}\sum_TT^3
  \ll\frac{\sqrt X\,Y^4}{q}.
 \]
 Summing $q^{-1}$ proves \eqref{eq:endpoint-G-global-low}; disjointness
 of the original rational cells avoids a union factor for the fourth power.

Here the height scales are dyadic, so their cubes form a geometric
progression dominated by the largest height. The subsequent sum over
all integer moduli is different: $\sum_{q\leq\sqrt X}1/q\ll L$.
This is the single logarithm in the global estimate for low values. We sum
integrals on disjoint physical cells, rather than fourth powers of a
sum of overlapping cell contributions.

\end{proof}

\section{An exact identity with four stages and its moments}
\label{sec:explicit-logs}

For the explicit logarithmic exponent in Theorem~\ref{thm:moment-main},
an identity with four stages limits the decomposition to at most eight variables
and $O(L^7)$ product boxes.  Its unmerged M\"obius atoms require the
term in the fourth moment $B+{\mathsf{M}}_i^2$.

\subsection{Fixed parameters}

The reduction in the number of variables saves logarithms when the product
boxes are combined. It also changes the length of a raw M\"obius factor:
the truncation now permits length of order $X^{1/4}$. Near the lower edge
of the terminal strip its square can exceed $B$. We retain the resulting
moment correction throughout the local estimates, and sum its contribution
separately at the end.

The moment and energy bounds are supplied by
Lemmas~\ref{lem:de-fixed-modulus-mean},
\ref{lem:original-factor-fourth-fixed-log} and Corollary~\ref{cor:box-energy}.
Lemma~\ref{lem:de-terminal-classification} applies by the argument
in Section~\ref{sec:c80-coverage}. Its subset estimates hold on the
full original character group.
All local characters are nonprincipal modulo one original modulus $q$;
induced characters are allowed. Energy and sampling estimates include
all characters where required. Put $L=\log X$, $\ell_2=\log\log(10X)$,
$H_{\mathrm M}\asymp L^s$, $1<s\le2$, $U\asymp T+H_{\mathrm M}$, $B=qU$, and
\[
 \mathfrak Z=1+\sqrt X/B,\qquad d_{\mathrm{box}}(n)=(2n-1)^2,\quad c_{\mathrm{box}}(n)=d_{\mathrm{box}}(n)+3\mathbf1_{n=1}.
\]
Let $S$ be any finite nonempty set of points given by a character and an ordinate,
$|t|\le C_UU$, separated by one for each character, and write $R=|S|\ll B$.
The fixed $C_U$ permits adaptive enlargement, uniformly over all subsets.
The energy in \eqref{eq:de-energy-definition} satisfies $\mathcal E(S)\ll R^3$.
For limiting arguments put
$\beta=\log_XB$, $r=\log_XR$, $e=\log_X\mathcal E(S)$ and
$H_{\rm pkt}=\log_X(p^2\mathcal E(S)/R)$.

Thus $\beta,r,e$ record the powers of $X$ in the analytic scale, the
number of sample points, and their energy. The quantity $H_{\rm pkt}$
records the power in the expression to be estimated, before logarithmic
factors are counted. In particular, the target scale $XB^2$ has exponent
$1+2\beta$, which becomes $2$ at $\beta=1/2$. A fixed strict inequality
below $2$ will absorb any fixed logarithmic cost; equality cases need
the more precise estimates below.

Empty sets satisfy all bounds trivially. Constants may depend on fixed
orders, support ratios, $s$ and the allowance in the threshold bound, but not on the set, values or box.

Put $\vartheta=\Theta_0=\max(\Theta,7/8)\leq11/12$ and
$\delta=1-\vartheta>0$. Fix an allowance $C_0$ in the threshold bound large enough for
all adjusted product thresholds. First choose
\[
 0<\eta\leq\min(1/1000,\delta/1000).
\]
Choose the fixed neighborhoods in the equality classification inside
the dense regions and regions with short factors used below. As proved in
Section~\ref{sec:c80-coverage}, the finite subset tests give a
positive gap $\gamma>0$ on the compact complement of these neighborhoods, a continuity width $w_T>0$
and a tolerance for errors in exponents $\rho_T>0$. These are fixed before
the terminal width. Choose $\eta_0>0$ so that
\[
 \eta_0\leq\min(\eta,1/20000),\qquad
 44\eta+\eta_0<\delta/2,\qquad
 \eta_0<\delta/4,\qquad
 \eta_0<\min(\gamma/8,w_T/2).
\]
Choose all arbitrarily small coefficient allowances after these
margins, small enough that the resulting exponent errors are below
$\rho_T/2$; the fixed logarithmic errors and errors from support constants are
below the remaining $\rho_T/2$ once $X$ is sufficiently large. We work on
\begin{equation}\label{c80:eq:strip}
 X^{1/2-\eta_0}\leq B\ll\sqrt X L^s,
 \qquad \sqrt X L^{-C_0}\leq p\ll X^\vartheta L^{C_0}.
\end{equation}
The upper bound concerns the recorded lower threshold $p$. Support constants
are absorbed by the strict margins, and bounded $X$ by the final
constant. The number of subset tests is at most $2^8-1$ and is
independent of every allowance by a small power.

\subsection{The exact identity and factor moments}
Every subset product retains the original variable intervals. For such a
product $F_A$, let $\nu_A$, ${\mathsf{M}}_A$ and $A_A\in\{0,1\}$ denote its order,
length and logarithm indicator, with coefficients $a_A(n)$.
Write $A_{\infty,A}=\max_n|a_A(n)|$.
Corollary~\ref{cor:box-coefficient-norm} gives
\[
 \sum_\ell (|a_A|*|a_A|)(\ell)^2
       \ll {\mathsf{M}}_A^2L^{(2\nu_A-1)^2+4A_A}.
\]
Thus Lemma~\ref{lem:difference-energy-fixed-log} yields
\begin{equation}\label{c80:eq:new-energy}
 Y^4\mathcal E(S)\ll A_{\infty,A}^4{\mathsf{M}}_A^4 R+B {\mathsf{M}}_A^2 L^{d_{\mathrm{box}}(\nu_A)+4A_A}R^2
 \quad (|F_A|\ge Y\text{ on }S).
\end{equation}
The bound remains valid with the first term retained in all length
ranges.

The two terms have distinct uses. The second term involves the square
norm of the convolution coefficients and is the main term in the
logarithmic estimates. The first retains the maximum coefficient and
cannot simply be dropped when a test is long. We shall show that its
length factor is smaller by a fixed power on the sets where it is used.
This permits the elementary divisor bound for the maximum coefficient
to absorb it.
  Applying the sampled mean square to $F_A^2$, we also obtain
\begin{equation}\label{c80:eq:new-fourth}
 \sum_S|F_A|^4\ll (B+{\mathsf{M}}_A^2){\mathsf{M}}_A^2L^{d_{\mathrm{box}}(\nu_A)+4A_A}.
\end{equation}

For the terminal shells, the exact identity of depth four is
\begin{equation*}\label{c80:eq:hb4}
 \Lambda=\sum_{j=1}^4(-1)^{j-1}\binom4j
 \mu_{\le z}^{*j}*1^{*(j-1)}*\log,\qquad z=(C X)^{1/4},
\end{equation*}
on the support of the prime polynomial. Put $M=\mu_{\le z}$, $a=1*M$ and
$\delta_1(n)=\mathbf1_{n=1}$, with $1(n)=1$ and Dirichlet convolution $*$. Since
$\log=1*\Lambda$, the right side is
$\Lambda*\sum_{j=1}^4(-1)^{j-1}\binom4j a^{*j}
=\Lambda*(\delta_1-(\delta_1-a)^{*4})$.
Since $\delta_1-a$ vanishes up to $z$, its fourth convolution is
supported above $z^4$.  Choosing $C$ larger than the fixed support
constant proves the identity on the required support.

For completeness, if $n\leq z$, every divisor of $n$ is at most $z$,
and hence
\begin{equation*}
 a(n)=\sum_{d\mid n}\mu_{\leq z}(d)
     =\sum_{d\mid n}\mu(d)=\delta_1(n).
\end{equation*}
Thus a nonzero term in $(\delta_1-a)^{*4}(n)$ would require a
factorization of $n$ into four integers, each greater than $z$. Such a
term has $n>z^4$. Convolving once more with $\Lambda$ cannot decrease
this support. The binomial expansion is therefore an exact identity on
the support in use, with no remainder term to estimate there.

In the $j$th summand, the $j$ truncated M\"obius functions, $j-1$
constant functions, and one logarithm give $2j$ variables. This explains
both the upper bound of eight and the presence of exactly one logarithmic
atom in each nonzero product.

There are at most eight atoms, at most four of them truncated M\"obius atoms,
and one logarithmic ordinary atom. Terms with the logarithmic atom supported
only at one vanish. The product restriction fixes the last dyadic scale
within boundedly many choices, giving $O(L^7)$ boxes.  On the
original intervals, the common Mellin separation \eqref{eq:ab-smooth-HB}
has weight of fixed $L^1$ norm.

Each variable has $O(L)$ possible dyadic length intervals. Once the first
seven lengths are fixed in a product of eight variables, the restriction
that their product has size $X$ leaves only boundedly many possible
last scales. Products with fewer variables cost no more. The bounded
$L^1$ norm of the common Mellin weight then allows us to select a shifted
product value from this finite sum without another logarithmic count.

Merging the shortest factors first at the actual $B$ gives merged
groups whose squared lengths, and those of their raw subproducts,
are $O(B)$.  For each group put
\[
 c_i=c_{\mathrm{box}}(\nu_i)+4A_i,
 \quad \mathfrak Z_i=\begin{cases}1+{\mathsf{M}}_i^2/B&\text{unmerged M\"obius atom},\\1&\text{otherwise}.\end{cases}
\]
The full sampled moment is
\begin{equation}\label{c80:eq:moment}
 \sum_S |S_i|^4\ll B {\mathsf{M}}_i^2 L^{c_i}\mathfrak Z_i,\qquad \mathfrak Z_i\ll \mathfrak Z.
\end{equation}
For short groups use \eqref{c80:eq:new-fourth}; for an unmerged M\"obius
atom retain $B+{\mathsf{M}}_i^2$, with ${\mathsf{M}}_i\ll X^{1/4}$. For long ordinary atoms apply
Lemma~\ref{lem:original-factor-fourth-fixed-log}, with singleton cost four
and the additional $4A_i$. Sampling permits locations depending on values.
Away shells continue to use the identity with five stages.

The correction for a raw M\"obius atom can be read directly from
\begin{equation*}
 (B+{\mathsf{M}}_i^2){\mathsf{M}}_i^2
   =B{\mathsf{M}}_i^2\left(1+\frac{{\mathsf{M}}_i^2}{B}\right),
 \qquad {\mathsf{M}}_i^2\ll\sqrt X.
\end{equation*}
A merged short group already has ${\mathsf{M}}_i^2\ll B$, so this factor is
bounded there. A long ordinary atom is controlled by the separate
theorem for ordinary factors. These are different reasons for the same
displayed moment shape; in particular, the result for ordinary factors is
not being applied to a raw M\"obius coefficient.

Retain all original atoms supported at the integer one. When the
terminal count is at most four, split merged groups until there are
four groups or all groups are raw. No further splitting is made in
this case.
For five terminal groups choose a group $h$ of maximal order, breaking
ties by shortest length, and split it into its original factors only if
$a_h=\log_X{\mathsf{M}}_h<5\eta$.  For six to eight groups retain the terminal
factors.  These choices preserve the exact product and are made
before the values and sample set are chosen.

The plan consequently depends only on the original intervals and their
orders. When we later pass to a subset of sample points, we keep the
same product and the same plan. What may change is the subset's
cardinality and the numerical values at its points. Uniformity under
this change is the reason for proving the estimates in the next section
on every separated subset.

\section{Direct product estimates on every separated subset}

The subsequent allocation of weights requires bounds for the full product
on every separated subset.  Scalar summation and subset tests give
such bounds.  A bound is \emph{hereditary} if it holds on every
subset, with that subset's own cardinality and the same product threshold.

\subsection{Positive energy and nonnegative weights}
For nonnegative weights $a_z$ define $\mathcal E(a)$ by inserting the product of its four weights in each
sharp energy quadruple. Let $G_K(a)$ be the same weighted sum with ordinate discrepancy $\Delta$
weighted by $\exp(-\Delta^2/K^2)$ instead of $\mathbf1_{|\Delta|\le1}$.
\begin{lemma}\label{c80:lem:weighted}
For $K\ge1$, $\mathcal E(a)\le eG_1(a)$ and $G_K(a)\ll K\mathcal E(a)$. The fourth root of $G_K$ is the $L^4$ norm
of the corresponding exponential sum on a positive measure space, up to its fixed normalization.
If every subset of a fixed separated class obeys $\mathcal E(S)\le D|S|$, then for $0\le a_z\le a_{\max}$,
$m=\sum a_z$,
\[
 \mathcal E(a)\ll D a_{\max}^3m.
\]
\end{lemma}
\begin{proof}
Bin pair sums into intervals of length one half, separately for each product
character, with weighted populations $b_{\psi,j}$. Then
$\sum b_{\psi,j}^2\le\mathcal E(a)$ and, by Cauchy,
$\sum b_{\psi,j}b_{\psi,j+k}\le\sum b_{\psi,j}^2$.
Summing Gaussian decay in $k$ gives $G_K(a)\ll K\mathcal E(a)$;
the reverse comparison is pointwise on $|\Delta|\le1$.
Gaussian Fourier inversion and character orthogonality express $G_K$
as a fourth power integral against a positive measure.

More explicitly, for points $z=(t_z,\chi_z)$ modulo $q$, we have
\begin{equation*}
 G_K(a)=\frac{K}{2\sqrt\pi\,\varphi(q)}
 \sum_{b\bmod q}^{*}\int_{\mathbb R}e^{-K^2x^2/4}
 \left|\sum_z a_z\chi_z(b)e^{it_zx}\right|^4\,dx.
\end{equation*}
Expanding the fourth power, the sum over $b$ imposes equality of the
product characters, and the Gaussian integral supplies
$e^{-\Delta^2/K^2}$. Positivity of this measure lets us use the ordinary
$L^4$ triangle inequality. This is why the Gaussian form is useful even
though the final energy has a sharp cutoff.

For the pair bin comparison, two pair sums in the same bin of length one half
differ by at most one. Pairs in bins $k$ apart carry a Gaussian weight
bounded by a constant times $e^{-c k^2/K^2}$. Cauchy's inequality bounds
their cross population by the same diagonal sum of squares, and
$\sum_k e^{-c k^2/K^2}\ll K$. This accounts for the dependence on $K$.

Writing $a_z=\int_0^{a_{\max}}\mathbf1_{a_z>t}\,dt$ and applying
Minkowski's inequality, the two comparisons and the hereditary
bound, we obtain
\[
 G_1(a)^{1/4}\ll D^{1/4}\int_0^{a_{\max}}|\{z:a_z>t\}|^{1/4}dt
 \le D^{1/4}a_{\max}^{3/4}m^{1/4}.
\]
Taking fourth powers proves the result, as well as the positive energy
triangle and cubic union inequalities.

In the calculation by layers,
$\int_0^{a_{\max}}|\{z:a_z>t\}|\,dt=\sum_z a_z=m$.
H\"older in the $t$ variable therefore gives the factor
$a_{\max}^{3/4}m^{1/4}$. The hypothesis on every subset is used
precisely on the level sets $\{z:a_z>t\}$; a bound for the full class
alone would not justify this step.

For a union of $J_0$ classes, the same argument gives
\begin{equation*}
 \mathcal E\left(\bigcup_{j=1}^{J_0}S_j\right)^{1/4}
 \ll\sum_{j=1}^{J_0}\mathcal E(S_j)^{1/4},
 \qquad
 \mathcal E\left(\bigcup_{j=1}^{J_0}S_j\right)
 \ll J_0^3\sum_{j=1}^{J_0}\mathcal E(S_j).
\end{equation*}
We will use the fourth root form when a geometric sum of amplitudes
gives a better result than counting all classes.

\end{proof}
\subsection{Two scalar summations retaining both energy terms}
\begin{lemma}\label{c80:lem:complement}
Suppose $|F_A S_1S_2S_3|\ge p$ on $S$ and ${\mathsf{M}}_A {\mathsf{M}}_1{\mathsf{M}}_2{\mathsf{M}}_3\asymp X$. Then
\[
 \frac{p^2\mathcal E(S)}R\ll XB^2(\mathfrak Z_1\mathfrak Z_2\mathfrak Z_3)^{1/2}L^{(c_1+c_2+c_3)/2}
 \left(L^{d_A/2}+\frac{A_{\infty,A}^2{\mathsf{M}}_A}{\sqrt{BR}}\right),
 \quad d_A=d_{\mathrm{box}}(\nu_A)+4A_A.
\]
In particular a fixed margin ${\mathsf{M}}_A^2/(BR)\ll X^{-c_*}$ absorbs the second term without a logarithmic loss.
\end{lemma}
\begin{proof}
Partition $|F_A|$ into dyadic amplitudes $Y$, with populations
$r_Y$.  By H\"older's inequality and \eqref{c80:eq:moment}, we have
$r_Y\le MY^{4/3}$, where
\[
 M\ll p^{-4/3}B(X/{\mathsf{M}}_A)^{2/3}L^{\sum c_i/3}(\mathfrak Z_1\mathfrak Z_2\mathfrak Z_3)^{1/3}.
\]

Indeed, on the slice $Y\leq|F_A|<2Y$, the other three factors have
product at least $p/(2Y)$. Thus
\begin{equation*}
 r_Y\left(\frac{p}{2Y}\right)^{4/3}
 \leq \sum_{S_Y}|S_1S_2S_3|^{4/3}
 \leq \prod_{i=1}^3
          \left(\sum_S|S_i|^4\right)^{1/3}.
\end{equation*}
Substituting their moment estimates and ${\mathsf{M}}_1{\mathsf{M}}_2{\mathsf{M}}_3\asymp X/{\mathsf{M}}_A$
gives the stated $M$. The upper bound $r_Y\leq R$ remains available
at the same time.

On each slice trivial energy and \eqref{c80:eq:new-energy} give
$\mathcal E(S_Y)\ll\min(r_Y^3,(A_0r_Y+Dr_Y^2)/Y^4)$,
where $A_0=A_{\infty,A}^4{\mathsf{M}}_A^4$ and $D=B{\mathsf{M}}_A^2L^{d_A}$.
Lemma~\ref{c80:lem:weighted} gives
$\mathcal E(S)^{1/4}\ll\sum_Y\mathcal E(S_Y)^{1/4}$.
Using $\min(u,v+z)\le\min(u,v)+\min(u,z)$, we separate the two
sums.  For the $D$ sum, split at $Y_D=(DR^2/M^3)^{1/8}$, using
$r_Y^{3/4}\le M^{3/4}Y$ below this point and $D^{1/4}R^{1/2}/Y$
above it.  The two geometric sums give $O(D^{1/8}M^{3/8}R^{1/4})$.

The balancing point solves
$M^{3/4}Y=D^{1/4}R^{1/2}/Y$, or $Y^8=DR^2/M^3$.
Below this point the dyadic summands grow like $Y$; above it they
decrease like $Y^{-1}$. Each sum is therefore bounded by a constant
times its value at the transition. This is an amplitude sum, not a
count of the number of dyadic slices, and it introduces no factor $L$.

For the $A_0$ sum, splitting at $(A_0R/M^3)^{1/8}$ gives
$O(A_0^{1/8}M^{3/8}R^{1/8})$.  Taking fourth powers, we obtain
$\mathcal E(S)\ll M^{3/2}(D^{1/2}R+A_0^{1/2}R^{1/2})$, proving the claim.
The divisor bound $A_{\infty,A}^4\ll X^{c_*/2}$ absorbs the second term
using the margin for the whole set.

In fact the unwanted factor satisfies
\begin{equation*}
 \frac{A_{\infty,A}^2{\mathsf{M}}_A}{\sqrt{BR}}
 \ll X^{c_*/4}X^{-c_*/2}=X^{-c_*/4}.
\end{equation*}
The original set's cardinality $R$ occurs here because the scalar sums
have already been combined. Replacing it by the cardinality of each
individual slice would demand a length margin that need not hold.

\end{proof}

\begin{lemma}\label{c80:lem:selector}
Put $f_i=|S_i|/\sqrt{{\mathsf{M}}_i}$, $F=\prod_i f_i\ge v$ on $S$, and $p\asymp\sqrt Xv$. Let a fixed finite
distribution of nonempty subset tests have probabilities $\lambda_A$, inclusion probabilities
$\theta_i\le1/2$, energy exponents $d_A$, and
\[
 r_i=2-4\theta_i,\quad \sigma=\sum r_i\le4,\quad
 D=\sum_A\lambda_A d_A,\quad M=\sum_i r_ic_i/4,\quad
 \mathfrak Z_Q=\prod_i \mathfrak Z_i^{r_i/4}.
\]
Suppose all test lengths obey ${\mathsf{M}}_A^2\ll BRX^{-c_*}$ for a fixed $c_*>0$ on the current set. Then
\begin{equation}\label{c80:eq:selector}
 \mathscr W:=\frac{p^2\mathcal E(S)}{XB^2R}
 \ll L^{D+M}\mathfrak Z_Q(B/R)^{\sigma/4-1}.
\end{equation}
Polynomial factor upper bounds and a polynomial lower bound for $p$ suffice. No common individual factor thresholds are assumed.
\end{lemma}
\begin{proof}
Let $Q=\prod f_i^{r_i}$. When $\sigma>0$ put $q_*=4/\sigma\ge1$ and $Q_0=(B/R)^{\sigma/4}L^M\mathfrak Z_Q$. Strong moment H\"older gives
\[
 \sum_S Q^{q_*}\ll B L^{q_*M}\mathfrak Z_Q^{q_*}=RQ_0^{q_*}.
\]

For each $i$ with $r_i>0$, the H\"older weight is $r_i/\sigma$,
and these weights sum to one. Consequently
\begin{equation*}
 \sum_S\prod_i f_i^{4r_i/\sigma}
 \leq\prod_{i:r_i>0}
          \left(\sum_S f_i^4\right)^{r_i/\sigma}.
\end{equation*}
The moment bound for $S_i$ says
$\sum_S f_i^4\ll BL^{c_i}\mathfrak Z_i$.
Its substitution gives one factor $B$, rather than one factor for each
$i$, and gives exactly the powers of $L$ and $\mathfrak Z_i$ in the
preceding display.

The exact identity $\prod_A(\prod_{i\in A}f_i)^{4\lambda_A}=F^2/Q$ implies that at each point some test satisfies
\[
 (\prod_{i\in A}f_i)^4/L^{d_A}\ge v^2/(QL^D).
\]

To justify the selection, take the weighted geometric mean of the
quantities on the left over all tests. The exponent of $f_i$ is
$4\theta_i=2-r_i$, while the logarithmic exponent is $-D$.
The mean is therefore $F^2/(QL^D)$, which is at least
$v^2/(QL^D)$. At least one member of a finite family is at least its
weighted geometric mean. Choosing the first such member makes a
definite partition without imposing individual thresholds for factor values.

Choose the first such test and partition further by
$2^jQ_0\le Q<2^{j+1}Q_0$.  Denote the populations by $R_{A,j}$,
$R_j=\sum_A R_{A,j}$.  Applying the full energy inequality with two terms
and dividing by ${\mathsf{M}}_A^2$, we obtain
\[
 v^2\mathcal E(S_{A,j})\ll L^DQ_0 2^j(KR_{A,j}+BR_{A,j}^2),\qquad
 K=\max_A A_{\infty,A}^4{\mathsf{M}}_A^2L^{-d_A}.
\]
Since $\sum R_j=R$ and $\sum 2^{q_*j}R_j\ll R$, we sum
geometrically for $j<0$ using $R_j\le R$, and for $j\ge0$ using
$R_j\ll R2^{-q_*j}$, to obtain
\[
 \sum_j2^{j/4}R_j^{1/2}\ll R^{1/2}.
\]

For negative $j$, the factor $2^{j/4}$ is summable. For nonnegative
$j$, the moment bound gives a summand at most
$R^{1/2}2^{j(1/4-q_*/2)}$, again summable because $q_*\geq1$.
The corresponding term linear in the population has exponent
$1/4-q_*/4$; it is strictly negative unless $q_*=1$. This explains
why only that endpoint needs the separate argument with finitely many slices.

They also give $\sum_j2^{j/4}R_j^{1/4}\ll R^{1/4}$ if $q_*>1$. If $q_*=1$, there are $O(L)$
occupied nonnegative slices by polynomial bounds; H\"older instead gives
$\sum_{j\ge0}(2^jR_j)^{1/4}\ll L^{3/4}R^{1/4}$. Gaussian Minkowski and the bounded test count
therefore yield
\[
 v^2\mathcal E(S)\ll L^DQ_0(BR^2+KRL^3).
\]
The margin for the whole set and $A_{\infty,A}^4\ll X^{c_*/2}$
give $KL^3\ll BR$.  Substituting $Q_0$ and $p^2\asymp Xv^2$, we
obtain \eqref{c80:eq:selector}.  When $\sigma=0$, take $Q=Q_0=1$,
$M=0$, $\mathfrak Z_Q=1$, and apply the calculation to the finite
selector partition alone.  The absorption in both cases uses the
margin for the whole set, without requiring one on each scalar slice.

The possible factor $L^3$ occurs only next to $K$. For every test the
length margin gives
\begin{equation*}
 \frac{A_{\infty,A}^4{\mathsf{M}}_A^2L^{-d_A}L^3}{BR}
 \ll X^{-c_*/2}L^3,
\end{equation*}
after the fixed allowance for the divisor bound. This tends to zero. Thus the
endpoint scalar summation adds no logarithmic loss to the main
$BR^2$ term. This is also why the full energy inequality with two terms was
retained until the summation was complete.

\end{proof}

\subsection{Products with at most four groups}

\begin{lemma}\label{lem:c100-four-groups}
Use the decomposition with four stages and the deterministic splitting plan
described above, retaining original atoms supported at the integer one
as atoms throughout the plan.  If at most four groups remain, then for every
nonempty separated subset $S$ on which the full product has magnitude
at least $p$,
\[
 p^2\mathcal E(S)\ll XB^2\mathfrak Z L^{97/2}|S|.
\]
There is no upper bound on $p$ in this assertion.
\end{lemma}

\begin{proof}
Let $b\leq4$ be the number of groups and $R=|S|$.
H\"older's inequality and the individual fourth moments give
\[
 R^{b/2}p^2\ll B^{b/2}X
       L^{\frac12\sum_i c_i}\prod_i\mathfrak Z_i^{1/2}.
\]
Since $\mathcal E(S)\ll R^3$ and $R\ll B$, it follows that
\begin{equation}\label{eq:c100-direct-four}
 \frac{p^2\mathcal E(S)}{XB^2R}
 \ll L^{\frac12\sum_i c_i}
       \prod_i\mathfrak Z_i^{1/2}(R/B)^{2-b/2}
 \ll L^{\frac12\sum_i c_i}\prod_i\mathfrak Z_i^{1/2}.
\end{equation}

If $b<4$, every group is a raw atom, by the splitting rule.
The original identity has $2j$ atoms, precisely $j$ of them M\"obius
atoms, so there are at most two M\"obius atoms in this case.
The moment correction in \eqref{eq:c100-direct-four} is at most
a constant times $\mathfrak Z$.

If $b=4$, at most three groups can be unmerged M\"obius atoms,
since the unique logarithmic ordinary atom belongs to another group.
With at most two such groups the same correction bound is immediate.
If there are three, the remaining group is merged: there were at
least three M\"obius atoms and hence at least six original atoms.
Write ${\mathsf{M}}_h$ for the length of this merged group.
Each raw M\"obius atom has length $O(X^{1/4})$, while the product
of all group lengths is comparable to $X$.  Hence
\[
 {\mathsf{M}}_h\gg X^{1/4}.
\]
A merged group satisfies ${\mathsf{M}}_h^2\ll B$, so $B\gg\sqrt X$.
Consequently $\mathfrak Z=1+\sqrt X/B=O(1)$ and all three
M\"obius corrections are bounded.  Thus in every case
\[
 \prod_i\mathfrak Z_i^{1/2}\ll\mathfrak Z.
\]

It remains to count logarithms.
For four positive integer orders with sum at most eight, the largest
value of $\sum_i c_{\mathrm{box}}(\nu_i)$ occurs at
$(5,1,1,1)$; it is $81+3\cdot4=93$.
This follows directly from the convex increase of $(2\nu-1)^2$,
or by transferring excess order to a largest group.
The unique logarithmic atom adds four to the moment costs, and
H\"older takes one half of their sum.  Therefore
\[
 \frac12\sum_i c_i\leq\frac{93+4}{2}=\frac{97}{2}.
\]
When $b<4$, all orders are one and the cost is at most eight.
Substituting these facts into \eqref{eq:c100-direct-four} proves
the lemma for the full set and every separated subset.
\end{proof}

\subsection{Five groups with a factor of positive length}
When $h$ has maximal order and $a_h\ge5\eta$, retain four individual value bins. In the dense
neighborhood $R,B\ge X^{1/2-\eta}$, if $\nu_h\ge2$, choose the shortest pair among the other four.
Its length is at most $X^{1/2-5\eta/2}$; its complement contains the short group $h$ and at most
two enlarged M\"obius moments. If all five groups are raw, every length exponent is at least
$5\eta$. With at most two M\"obius atoms the shortest pair is sufficient. With four M\"obius
atoms, use the shortest two of them: if $y$ is the ordinary exponent, their sum is at most
$(1-y)/2$. With three M\"obius and two ordinary atoms choose the shortest of each type, with
exponents $x,y$. Then $3x+2y\le1$, so $x+y\le1/2-x/2\le1/2-5\eta/2$. In each case the pair leaves
at most two M\"obius atoms outside it. Its margin for the squared length is at most $X^{-3\eta}$ after
harmless support adjustments.

The purpose of the pair is to leave three factors to which the complement
lemma applies. The first requirement is a fixed saving in its squared
length relative to $BR$. The second is to leave at most two enlarged
M\"obius moments in the complement, so their square root corrections
multiply to at most $\mathfrak Z$. These two requirements explain why
the pair is selected by type as well as by length when all groups are raw.

Lemma~\ref{c80:lem:complement} gives local cost
\[
 \frac{d_{\mathrm{box}}(\nu_A)+\sum_{i\notin A}c_{\mathrm{box}}(\nu_i)}2+2\le35.
\]
For maximal order $4,3,2,1$, the worst costs are respectively
$35,31,35,25/2$: the corresponding excess orders outside $h$ are
$0,1,2,0$. This allows every pair excluding $h$.

The largest costs can be checked directly from the integer orders.
For maximal order four, a pair of order two leaves orders $4,1,1$,
and its cost is $(9+49+4+4)/2+2=35$.
For maximal order three, the largest pair cost is obtained from an
pair of order three leaving orders $3,1,1$; it is
$(25+25+4+4)/2+2=31$.
For maximal order two, a pair of order four may leave orders $2,1,1$,
giving $(49+9+4+4)/2+2=35$.
When all orders are one, the cost is $(9+4+4+4)/2+2=25/2$.
The bounded number of order patterns does not contribute a power of $L$.

The sparse equality neighborhood belongs to the case with short $h$ below;
the compact complement is treated in Section~\ref{sec:c80-coverage}.

\subsection{Subset tests for the remaining products}
For five terminal groups with $a_h<5\eta$, pair separation with $h$ gives $a_j\ge\beta/2-a_h+o(1)$
for every other group. For any chosen leaf $i$, summing this inequality over the other three
leaves yields
\[
 a_i\le1-3\beta/2+2a_h+o(1)\le1/4+3\eta_0/2+10\eta+o(1)<1/4+17\eta.
\]
Thus all four other lengths are at most $X^{1/4+17\eta}$, including M\"obius leaves.

The subtraction in this estimate is useful to spell out. The total
exponent of all groups is $1+o(1)$. Subtracting the heavy exponent
$a_h$ and the lower bounds for the three leaves other than $i$ gives
\begin{equation*}
 a_i\leq1-a_h-3(\beta/2-a_h)+o(1)
      =1-3\beta/2+2a_h+o(1).
\end{equation*}
A selected subset of the heavy atoms has length $O({\mathsf{M}}_h)$, so adjoining
it to one of these leaves gives exponent at most
$1/4+17\eta+5\eta$. The smallness of the whole heavy group controls
every one of these subset lengths.

For $R\le X^{\delta/2}$, the upper bound gives
$p^2\mathcal E/R\ll X^{2-\delta}L^{2C_0}$, a fixed saving over
$XB^2$.  For the remaining sets, expose the $m$ atoms of $h$.
Choose a uniformly random subset $J$ of half the size (sizes
$(m-1)/2,(m+1)/2$ with equal probability if $m$ is odd), and
independently choose one of the four complementary groups uniformly.
Every resulting test has length at most $X^{1/4+22\eta}$, so
\[
 {\mathsf{M}}_A^2/(BR)\ll X^{44\eta+\eta_0-\delta/2}.
\]
The heavy atoms have inclusion probability $1/2$ and the four complementary groups $1/4$.
Lemma~\ref{c80:lem:selector} gives $\sigma=4$ and correction at most $\mathfrak Z$. Writing
$n=\sum n_i$, $S_2=\sum n_i^2$, $n_1=\#\{i:n_i=1\}$, its cost, including the unique logarithm, is
\begin{equation}\label{c80:eq:half-cost}
 m(m+n)+\mathbf1_{m\text{ odd}}+2S_2-2(m+n)+4+3n_1/4\le31.
\end{equation}
Indeed $\mathbb E|J|=m/2$, $\mathbb E|J|^2=(m^2+\mathbf1_{m\text{ odd}})/4$; substitute in $\mathbb E
d_{\mathrm{box}}(|J|+n_i)+\sum c_{\mathrm{box}}(n_i)/4$ and add two. Here $m\ge\max n_i$ and
$m+n\le8$. Convexity concentrates the allowed excess in the interval $[1,m]$; for $m=1,2,3,4$ the
maxima are respectively $11,51/2,117/4,31$. Smaller totals give smaller costs.

For the algebra, let the chosen leaf have order $n_i$. Averaging first
over the subset of half the size and then over the four leaves gives
\begin{align*}
 \mathbb E\,d_{\mathrm{box}}(|J|+n_i)
  &=m^2+\mathbf1_{m\text{ odd}}+mn+S_2-2m-n+1,\\
 \frac14\sum_i c_{\mathrm{box}}(n_i)
  &=S_2-n+1+\frac34n_1.
\end{align*}
Adding these quantities and the logarithmic cost two gives the
displayed polynomial. At the largest allowed total, the extremal leaf
orders are
\begin{equation*}
 \begin{array}{c|c|c}
 m&\text{four leaf orders}&\text{cost}\\ \hline
 1&(1,1,1,1)&11\\
 2&(2,2,1,1)&51/2\\
 3&(2,1,1,1)&117/4\\
 4&(1,1,1,1)&31 .
 \end{array}
\end{equation*}
The quadratic dependence on the orders explains why concentrating the
available excess gives these maxima.

For six groups select a longest $h$, and test the ten pairs outside it with equal probability.
Every such pair has length at most $X^{11/24+3\eta/4}$: its complement contains $h$ of exponent at
least $1/6$, and a triple of length at least $B^{3/4}$ by pair separation. In the dense
neighborhood its ratio of the squared length to $BR$ is at most $X^{-1/12+7\eta/2}$. Residual moment
weights are $1/2$ at $h$ and $1/10$ at each other group; $\sigma=4$. Thus the correction is at
most $\mathfrak Z^{4/5}$ and the local exponent is
\begin{equation}\label{c80:eq:six-cost}
 C_6=\frac1{10}\sum_{|A|=2,\ h\notin A}d_{\mathrm{box}}(\nu_A)
       +\frac12c_{\mathrm{box}}(\nu_h)+\frac1{10}\sum_{i\ne h}c_{\mathrm{box}}(\nu_i)+2
       \le331/10.
\end{equation}
Here and in the next two displays, $c_{\mathrm{box}},d_{\mathrm{box}}$ exclude the unique logarithm; its contribution is the displayed two.

For six groups, each of the five groups outside $h$ occurs in four of
the ten pairs. Its inclusion probability is therefore $2/5$, its
residual exponent is $2-4(2/5)=2/5$, and its moment weight is $1/10$.
The group $h$ has moment weight $1/2$. The sum of residual exponents is
$2+5(2/5)=4$. At most four raw M\"obius atoms can be present, so their
total moment weight is at most $1/2+3/10=4/5$.

The length computation uses pair separation in the complementary
triple. Adding its three pair inequalities gives a total exponent at
least $3\beta/4$. The longest group has exponent at least $1/6$.
Consequently the selected pair has exponent at most
$1-1/6-3(1/2-\eta)/4=11/24+3\eta/4$.
Squaring and dividing by $BR\geq X^{1-2\eta}$ gives
$11/12+3\eta/2-1+2\eta=-1/12+7\eta/2$.

For seven groups select a longest $h$ and test the twenty triples outside it. Their lengths are at
most $X^{27/56+3\eta/4}$ because the complement contains $h$ of exponent at least $1/7$ and a
triple of length at least $B^{3/4}$. The ratio margin is $X^{-1/28+7\eta/2}$. Now $\sigma=2$, with
only $h$ residual, and Lemma~\ref{c80:lem:selector} gives $\mathscr W\ll L^{A_7}\mathfrak
Z_h^{1/2}g^{-1/2}$, where $g=B/R$ and
\[
 A_7=\frac1{20}\sum_{|A|=3,\ h\notin A}d_{\mathrm{box}}(\nu_A)+\tfrac12c_{\mathrm{box}}(\nu_h)+2.
\]
Direct H\"older for fourth moments gives $\mathscr W\ll L^{D_7}\prod \mathfrak Z_i^{1/2}g^{3/2}$ with
$D_7=\sum c_{\mathrm{box}}(\nu_i)/2+2$. Taking powers $3/4$ and $1/4$ cancels $g$. Thus
\begin{equation}\label{c80:eq:seven-cost}
 C_7=(3A_7+D_7)/4\le283/8,
\end{equation}
and the correction is at most $\mathfrak Z^{3/8+4/8}=\mathfrak Z^{7/8}$.

Each of the six groups outside $h$ belongs to half of the twenty
triples, so only $h$ has a residual moment. Its length estimate is
obtained by replacing $1/6$ in the preceding calculation by $1/7$:
$1-1/7-3/8=27/56$. The two available estimates have opposite
powers of $g$. Multiplying the selector estimate to power $3/4$
and the direct estimate to power $1/4$, we have
\begin{equation*}
 g^{-(1/2)(3/4)+(3/2)(1/4)}=1.
\end{equation*}
The selector correction contributes at most
$\mathfrak Z^{(1/2)(3/4)}$, and the direct correction contributes at
most $\mathfrak Z^{(4/2)(1/4)}$. This gives the displayed exponent
$7/8$, while the same interpolation gives $(3A_7+D_7)/4$ for
the logarithmic cost.

For eight groups test all fifty-six triples. Their complement consisting of five groups has length at least
$B^{5/4}$, so ${\mathsf{M}}_A^2/(BR)\ll X^{-1/4+9\eta/2}$. Every inclusion probability is $3/8$, every
residual exponent $1/2$, and $\sigma=4$. All eight orders are one; hence
\[
 C_8=d_{\mathrm{box}}(3)+8c_{\mathrm{box}}(1)/8+2=25+4+2=31,
 \qquad \mathfrak Z_Q\ll \mathfrak Z^{1/2}.
\]

For a complement consisting of five groups, summing all ten pair separation
inequalities counts each group four times. Its total length exponent
is therefore at least $5\beta/4$. This gives
\begin{equation*}
 \log_X\frac{{\mathsf{M}}_A^2}{BR}
 \leq 2-\frac72\beta-r+o(1)
 \leq-\frac14+\frac92\eta+o(1).
\end{equation*}
For the moment correction, each residual exponent $1/2$ has moment
weight $1/8$, and at most four of the atoms are M\"obius atoms.
Their total weight is thus at most $4/8=1/2$.

The unique logarithm always costs two in these selector calculations.
If its inclusion probability is $\theta$, its contribution to the
averaged test cost is $4\theta$, and its residual moment contribution
is $2-4\theta$. Their sum is two.

There are at most two excess orders for six groups, and one for seven.
Concentrating or separating them gives the maxima
$331/10,122/5,51/2$ for $\nu_h=1,2,3$ in \eqref{c80:eq:six-cost},
and $283/8,113/4$ for $\nu_h=1,2$ in \eqref{c80:eq:seven-cost}.

For example, the largest cost for six groups has $\nu_h=1$ and other
orders $(3,1,1,1,1)$. Four pairs have order four and six have order
two, giving
\begin{equation*}
 \frac{4\cdot49+6\cdot9}{10}
 +\frac42+\frac{25+4\cdot4}{10}+2
 =\frac{331}{10}.
\end{equation*}
For $\nu_h=2$, the other orders can be $(2,1,1,1,1)$, giving
$(4\cdot25+6\cdot9)/10+9/2+(9+16)/10+2=122/5$.
For $\nu_h=3$, all five other orders are one and the cost is
$9+25/2+2+2=51/2$.

For seven groups with $\nu_h=1$, the largest pattern has other
orders $(2,1,1,1,1,1)$. Ten triples have order four and ten have
order three. Thus $A_7=37+2+2=41$ and
$D_7=(4+9+5\cdot4)/2+2=37/2$, whence
$(3A_7+D_7)/4=283/8$. With $\nu_h=2$, the other orders are all
one, so $A_7=25+9/2+2=63/2$, while $D_7=37/2$.
This gives the other value $113/4$.

\subsection{Coverage on every subset}\label{sec:c80-coverage}
We now make the passage for finite $X$ and its uniformity explicit.
The original intervals, coefficients and deterministic plan are fixed
before choosing sample points. Restricting the points changes their
cardinality and energy, but not these data. The equality classification
is applied to the terminal groups before any internal atoms are exposed;
the plans for at most four groups and those with short $h$ already have direct estimates and
need no pair separation after splitting.

\paragraph{Records for finite $X$ from the source estimates.}
Let $5\leq k\leq8$, and let $S$ be any nonempty separated set on which
$|\prod_iS_i|\geq p$, with $p$ in \eqref{c80:eq:strip}.
The elementary divisor bounds give, for each fixed $\epsilon>0$,
$A_{\infty,A}\ll_\epsilon X^\epsilon$ and
$|S_i|\ll_\epsilon {\mathsf M}_iX^\epsilon$, uniformly over the
original intervals and the at most $255$ nonempty subset products.
The latter upper bounds and the lower bound for the full product give polynomial
lower bounds for every individual magnitude on $S$.
Thus dyadic binning of those magnitudes partitions $S$ into
$J_0\ll L^k$ nonempty classes $T$. On one such class let $v_i$ be the
dyadic lower endpoints. Since $\prod_i v_i\geq p/2^k$, we may fix
$C\geq2^8$ once and lower one endpoint to obtain $V_i\leq v_i$ with
$\prod_iV_i=p/C$. This changes only lower thresholds, not the points
or polynomials. Define
\begin{gather*}
 a_i=\log_X{\mathsf M}_i,\qquad u_i=\log_XV_i,\qquad
 x_i=a_i-u_i,\qquad D=\sum_i x_i,\\
 r_T=\log_X|T|,\qquad e_T=\log_X\mathcal E(T),\qquad
 H_T=2\sum_i u_i+e_T-r_T.
\end{gather*}
The exact identities, using $\prod_i{\mathsf M}_i\asymp X$, are
\begin{gather*}
 D=\sum_i a_i-\log_X(p/C)
   =1-\log_X(p/C)+O(1/\log X),\\
 H_T=\log_X\bigl((p/C)^2\mathcal E(T)/|T|\bigr).
\end{gather*}
In particular it is the recorded lower threshold, rather than the
possibly larger actual product, that controls $D$.

There are constants $C_*,C_\epsilon,X_\epsilon$, depending only on the
fixed support, smoothing, order and cap data and on $\epsilon$ where
indicated, such that all these records, for $X\geq X_\epsilon$, obey
\[
 \tau_X=C_*\epsilon+C_\epsilon\frac{\log\log(10X)}{\log X}
\]
and the following inequalities. Subscripts $A$ denote sums over $A$,
and $\beta=\log_XB$:
\begin{align*}
 &a_i\geq0,\qquad \left|\sum_i a_i-1\right|\leq\tau_X,
       \qquad a_i+a_j\geq\beta/2-\tau_X\quad(i\ne j),\\
 &x_i\geq-\tau_X,\qquad
       \delta-\tau_X\leq D\leq1/2+\tau_X,\\
 &0\leq r_T\leq\beta+\tau_X,
       \qquad 2r_T\leq e_T\leq3r_T+\tau_X,\\
 &r_T+2u_A\leq a_A+\max(a_A,\beta)+\tau_X
                         \quad(\varnothing\ne A\subseteq\{1,\ldots,k\}),\\
 &e_T\leq\max\{r_T+4x_A,\ \beta+2r_T-2a_A+4x_A\}+\tau_X
                         \quad(\varnothing\ne A\subseteq\{1,\ldots,k\}),\\
 &r_T+4u_i\leq\beta+2a_i+z_i+\tau_X.
\end{align*}
Here $z_i=0$ for a merged group or an ordinary atom; for a raw
M\"obius atom,
\[
 z_i=\max(0,2a_i-\beta),\qquad a_i\leq1/4+\tau_X,
 \qquad 0\leq z_i\leq|\beta-1/2|+2\tau_X.
\]

We give the sources of these finite errors. Length normalization and
terminal pair separation have errors $O(1/\log X)$. The individual
upper bound gives $x_i\geq-\epsilon-O_\epsilon(1/\log X)$.
The cap and exact identity for $D$ give its two bounds. Separation and
diagonal quadruples give
$|T|\ll B$ and $|T|^2\leq\mathcal E(T)\ll|T|^3$.
For every subset product, \eqref{eq:de-bin-mean} and
\eqref{eq:box-second-norm} give directly on $T$
\[
 |T|\prod_{i\in A}V_i^2
 \leq\sum_T|F_A|^2
 \ll(B+{\mathsf M}_A){\mathsf M}_A
           L^{2A_A+(\nu_A-1)^2}.
\]
Taking logarithms gives the mean inequality above.
Likewise \eqref{c80:eq:new-energy}, with $Y=\prod_{i\in A}V_i$,
gives before combining errors
\[
 e_T\leq\max\left\{\begin{aligned}
 &r_T+4x_A+4\log_X A_{\infty,A},\\
 &\beta+2r_T-2a_A+4x_A+
 (d_{\mathrm{box}}(\nu_A)+4A_A)\log_X L
 \end{aligned}\right\}+O(1/\log X).
\]
Both energy terms are retained, at every test length. Finally
\eqref{c80:eq:moment} gives the inequality for fourth moments, since
$\log_X(1+{\mathsf M}_i^2/B)
 =\max(0,2a_i-\beta)+O(1/\log X)$.
These are fresh applications of the estimates for the sampled mean, energy with two terms and
factor moments to $T$; they do not use a product estimate on
$S$. Their forms using suprema over bins permit selection depending on values.
Their constants are uniform over the interval endpoints with fixed
support ratios and bounded normalized smoothing weights. Only the
orders, factor types and positions of subset tests range over finite lists;
we do not assert that all boxes form a finite family as $X$ varies.
The fixed logarithmic costs and support constants therefore fit the
single $\tau_X$ above, independently of $q,U,p,S,T$ and the thresholds.

\paragraph{Equality neighborhoods and stability.}
At $\beta=1/2$ and zero error, these inequalities are precisely the
hypotheses of Lemma~\ref{lem:de-terminal-classification}; the raw
M\"obius corrections vanish. That lemma gives $H\leq2$, with equality
only in the dense configuration or the configuration of five groups with
one factor of zero length. Fix the open neighborhoods
\[
 N_d=\{r>1/2-\eta/2\},\qquad
 N_s=\bigcup_i\{a_i<\eta/4,\ a_j>1/8\ (j\ne i)\}
       \quad(k=5),
\]
and take $N_s=\varnothing$ for $k=6,7,8$.
In the exact feasible set restrict $\delta/2\leq D\leq3/4$.
This set is compact: $a_i\geq0$, $\sum_i a_i=1$, $x_i\geq0$
and $\sum_i x_i=D$ bound $a_i,x_i,u_i$, while
$0\leq r\leq1/2$ and $2r\leq e\leq3r$ bound the remaining
coordinates. The constraints, including their maxima, are continuous.
The closed complement of $N_d\cup N_s$ contains no equality point.
Thus $H\leq2-\gamma$ there for some $\gamma>0$, uniform over
$5\leq k\leq8$ and the finitely many order and type patterns
(empty complements require no bound).

There consequently exist $w_T,\rho_T>0$ such that every finite record
above outside $N_d\cup N_s$ satisfies $H_T\leq2-\gamma/2$ whenever
$|\beta-1/2|\leq w_T$ and $\tau_X\leq\rho_T$.
Otherwise records with widths and errors tending to zero would, after
fixing a pattern and passing to a subsequence, converge to an exact
feasible record outside these open neighborhoods with
$H\geq2-\gamma/2$. Boundedness follows from the relaxed normalization,
$x_i\geq-\tau_X$ and the bounds for $D$; the limiting deficit lies
in $[\delta,1/2]\subset[\delta/2,3/4]$.
The raw M\"obius correction tends to zero by its displayed bound.
This contradicts the exact gap.

For realizable boxes, $N_s$ is contained in the plan with short $h$.
Indeed, call its short group $0$. If one of the other four terminal
groups were merged, consider its last merger. A live ancestor of group
$0$ is unselected and has length $O({\mathsf M}_0)$, since all atomic
lower lengths are at least one. The rule that takes the shortest first makes both
selected inputs $O({\mathsf M}_0)$, so that terminal group has length
$O({\mathsf M}_0^2)$. Its pair separation from group $0$ gives
${\mathsf M}_0^3\gg\sqrt B$, or
$3a_0\geq\beta/2+O(1/\log X)$, incompatible with
$a_0<\eta/4$ in the terminal strip for large $X$.
The four other groups are therefore raw, of order one. If group $0$
has larger order, it is the unique group of maximal order; otherwise it
is the unique shortest group. In either case the prescribed rule
chooses $h=0$, and $a_h<5\eta$.

\paragraph{Joint choice of constants.}
Choose $\eta$, the neighborhoods, and then $\gamma,w_T,\rho_T$ as
above, before $\eta_0$. Choose $\eta_0$ as in the subsection on fixed parameters. The direct arguments have the positive savings in squared length
\[
 3\eta,\quad \delta/2-44\eta-\eta_0,\quad
 1/12-7\eta/2,\quad1/28-7\eta/2,\quad1/4-9\eta/2
\]
for positive five, short five above its cutoff for small cardinalities,
and six, seven and eight groups, respectively.
Fix $c_*>0$ smaller than all five savings. Now choose $\epsilon>0$
so that $C_*\epsilon<\rho_T/2$ and $4\epsilon<c_*/2$.
Only then choose one $X_0$ so that all source estimates apply, the
remaining part of $\tau_X$ is below $\rho_T/2$, and all strict margins
absorb their support errors. The upper extension
$\beta\leq1/2+O(\log L/L)$ is within $w_T$ for large $X$;
the lower extension is within it because $\eta_0<w_T/2$.
The complement lemma's coefficient term is then absorbed by
$X^{2\epsilon-c_*/2}$, and the selector's possible endpoint term
by $X^{4\epsilon-c_*}L^3$. Neither changes the main logarithmic cost.
There is no requirement that $C_\epsilon$ or $X_0$ be uniform as
$\epsilon$ tends to zero. All these choices precede every sample set
and every adjusted lower threshold.

\begin{proposition}\label{prop:c80-terminal}
With these fixed choices there are constants $X_0$ and $C_{\rm fin}$
such that for every $X\geq X_0$, every original modulus and height in
\eqref{c80:eq:strip}, every deterministic plan with four stages, and every
threshold $p$ in that strip, the following holds.
In the plan for five groups in the positive length case, fix also four individual value
bins and the fifth lower threshold, whose product is $p$.
On every nonempty separated subset $\mathcal S$ of that class, or of
$|\prod_iS_i|\geq p$ in any other plan,
\[
 p^2\mathcal E(\mathcal S)
 \leq C_{\rm fin}XB^2\mathfrak Z L^{c_b}|\mathcal S|.
\]
The costs for at most four, positive five, short five, six, seven
and eight groups are, respectively,
\[
 \frac{97}{2},\qquad35,\qquad31,\qquad
 \frac{331}{10},\qquad\frac{283}{8},\qquad31.
\]
The constants depend only on the fixed data, not on the box, subset
or threshold.
\end{proposition}

\begin{proof}
Take an arbitrary allowed $\mathcal S$ and put $R=|\mathcal S|$ anew.
For at most four groups, Lemma~\ref{lem:c100-four-groups} applies
at every cardinality. For short five, if $R\leq X^{\delta/2}$, then
\[
 p^2\mathcal E(\mathcal S)/R\ll p^2R^2
 \ll X^{2-\delta}L^{2C_0},
\]
which has a fixed saving over $XB^2$ because $2\eta_0<\delta$.
For larger $R$, apply the selector argument leading to
\eqref{c80:eq:half-cost}, using this $R$ and its displayed length
saving. For positive five and six to eight groups with
$R\geq X^{1/2-\eta}$, the direct dense arguments apply with this
$R$; the required $B\geq X^{1/2-\eta}$ follows from $\eta_0\leq\eta$.
Their costs are the positive five calculation,
\eqref{c80:eq:six-cost}, \eqref{c80:eq:seven-cost}, and the subsequent
calculation for eight groups. All their correction exponents are at most one.
Within each scalar summation the linear energy term is retained until
the internal slices have been summed; its length saving uses the
cardinality of that invocation's input set, not of every scalar slice.

It remains that $R<X^{1/2-\eta}$ in a positive five plan or a plan for six, seven
or eight groups. Partition this $\mathcal S$ into the classes of factor values $T_j$ constructed above, with $J_0\ll L^8$ and
$\sum_j|T_j|=R$. Recompute $r_{T_j}$ and $e_{T_j}$ on each class.
Since $r_{T_j}\leq\log_XR<1/2-\eta$, every class is outside $N_d$.
For six to eight groups $N_s$ is empty; a box in the positive five plan cannot
lie in $N_s$, since the argument using the history of merging would put its fixed
plan in the case with short $h$. Stability and the exact threshold identity
therefore give
\[
 p^2\mathcal E(T_j)\leq C^2X^{2-\gamma/2}|T_j|.
\]
The cubic union inequality in Lemma~\ref{c80:lem:weighted} now gives
\[
 \frac{p^2\mathcal E(\mathcal S)}{XB^2R}
 \ll L^{24}X^{-\gamma/2+2\eta_0}=O(1),
\]
since $\eta_0<\gamma/8$. This internal partition has been paid for
by a power saving and adds no logarithmic cost to the proposition.
In positive five the four prescribed external bins remain fixed;
the internal partition neither changes their definition nor requires
additional bins in the subsequent allocation.

Every case used the current $R$ with constants fixed in advance.
For a smaller subset the same case distinction is made again, without
changing the intervals, coefficients or plan. A dense parent's sparse
subset uses the argument on the compact complement, and a large short five
set's small subset uses the cap argument. No monotonicity of
$\mathcal E(\mathcal S)/|\mathcal S|$ is asserted or needed.
This proves the estimate on all subsets, including the level sets of
nonnegative weights used in the next section, without using that
section as an input.
\end{proof}

\section{Nonnegative allocation and the circle integral}
\label{sec:c80-global}

The energy of complete values remains to be bounded from the product
estimates.  Capping the available product sizes at the complete
threshold permits a nonnegative distribution among the products.

\subsection{Capped fractional product weights}

At a point with a complete value there may be several useful products,
of very different sizes. We distribute its unit weight among them using
their capped magnitudes. This is a finite allocation at each point: the
weights sum to one, and each is bounded by the product magnitude divided
by the complete threshold. The hereditary estimates allow us to bound
the resulting weighted classes. The amplitude sum will be geometric,
which is the saving over treating all possible product levels as equally
costly classes.

\begin{proposition}\label{c80:prop:fractional}
Let \(W\) be a separated set of complete values with lower
amplitude \(V>0\), and let \(J\) products with common coefficients satisfy
\[
 \sum_{\nu=1}^J\max_{|u|\leq H_{\mathrm M}}
       |F_\nu(t+u,\chi)|\geq cV
 \qquad((t,\chi)\in W),
\]
where \(H_{\mathrm M}\geq1\) and \(c>0\) is fixed.
Assign each product a plan from a fixed finite list.
For a plan \(b\), assume the uniform estimate
\[
 p^2\mathcal E(S)\leq A_b|S|
\]
on every separated subset of the permitted product classes,
for all adjusted product thresholds \(V/(C_1J)\leq p\leq V\),
where \(C_1\) is a sufficiently large fixed constant.
For each product assigned to plan \(b\) and each retained amplitude
level \(p\), let \(K_b\) uniformly bound the number of admissible tuples
of individual bins for factor values used to partition its selected points
into classes on which the hereditary estimate applies. The last factor
threshold, when present, is determined by \(p\) and the tuple.
For an estimate for the full product put \(K_b=1\). Thus \(K_b\) excludes the
counts of products, amplitude levels, and colors for shift separation.
All selected factors have polynomial upper bounds, and the product
threshold has a polynomial lower bound. Then
\begin{equation}\label{c80:eq:witness}
 V^2\mathcal E(W)\ll
 H_{\mathrm M}^4\max_b\{A_bJ^3K_b^3\}|W|.
\end{equation}
\end{proposition}

\begin{proof}
For \(z=(t,\chi)\in W\), put
\(M_\nu(z)=\max_{|u|\leq H_{\mathrm M}}|F_\nu(t+u,\chi)|\)
and choose an attaining shift. Then $\sum_\nu M_\nu(z)\ge cV$; we may take $c\le1$.

Discard the boxes with $M_\nu<cV/(2J)$.  For the others put
$x_\nu=\min(M_\nu,V)$ and $\alpha_\nu=x_\nu/\sum x_\nu$.
The denominator is at least $cV/2$, since either no value is
changed by capping or one term equals $V$.

The discarded terms have total less than $J\cdot cV/(2J)=cV/2$,
so the surviving uncapped terms have sum at least $cV/2$.
If none exceeds $V$, capping changes nothing. If some term exceeds
$V$, its capped value alone equals $V\geq cV/2$. Thus the
denominator stays positive in both cases, and normalizing by it gives
an explicit allocation with no existence argument left to supply.
  Hence
\[
 \sum_\nu\alpha_\nu(z)=1,\qquad \alpha_\nu(z)\le(2/c)x_\nu/V.
\]
Bin \(x_\nu\) by dyadic lower endpoints \(p\), with a saturated
top bin \(p=V\). On each bin \(M_\nu\geq p\) and
\(\alpha_\nu\ll p/V\). Only \(p\gg V/J\) occur.
A sufficiently large fixed $C_1$ covers these thresholds and the adjustments below.

For a nonsaturated bin, $p\leq x_\nu<2p$, so
$\alpha_\nu\leq(4/c)p/V$; the same conclusion holds in the
saturated bin $x_\nu=p=V$. The upper cap prevents levels above $V$,
even when the actual selected product is much larger. The lower cutoff
ensures that every retained threshold lies in the range where the
hereditary product estimate is assumed.

For a plan with a common threshold, bin the first $k_b-1$ actual group values at the shifted point. Each
is polynomially bounded above and below, so there are $K_b\ll L^{k_b-1}$ tuples. If their lower
endpoints are $l_i$, set $l_{k_b}=(p/C)/\prod_{i<k_b}l_i$ with fixed $C\ge2^8$. The actual last
value exceeds $l_{k_b}$, and the product threshold is $p/C\le V/C$.
Replacing $p$ by $p/C$ changes only a fixed constant. Other plans need no such bins.

For the last factor, the product lower bound and the upper endpoints of
the first $k_b-1$ bins give an actual value at least
$p/(2^{k_b-1}\prod_{i<k_b}l_i)$. This exceeds the chosen last
threshold. Each of the other factors has a polynomial lower bound
because the full product is bounded below and every remaining factor
has a polynomial upper bound. Hence each free factor has only $O(L)$
dyadic bins. The positive five plan will need four such free bins;
the plans using the full product need none.

For each character an original unit separated set shifted by at most $H_{\mathrm M}$ has
$O(H_{\mathrm M}+1)$ images in a unit interval, counting multiplicity. Color the images into
$O(H_{\mathrm M}+1)$ separated injective classes. All support intervals and coefficients stay
fixed. Each resulting weighted class $c$ has total mass $m_c$, maximum weight $O(p_c/V)$, and
$\sum_c m_c=|W|$. Lemma~\ref{c80:lem:weighted} gives
\[
 \mathcal E(\alpha_c)\ll A_{b(c)}p_c m_c/V^3.
\]

Here the hereditary constant in the lemma on weighted energy is
$D=A_{b(c)}/p_c^2$, and the maximum weight is $O(p_c/V)$.
Substitution gives
\begin{equation*}
 D\,a_{\max}^3m_c
 \ll \frac{A_{b(c)}}{p_c^2}
          \left(\frac{p_c}{V}\right)^3m_c
 =\frac{A_{b(c)}p_c m_c}{V^3}.
\end{equation*}
Coloring makes the shifted map injective, so the weights and total
masses of its preimages can be carried to these separated classes
without identifying two distinct original points.

Expand each original unit weight using $\sum\alpha_\nu=1$. Original energy quadruples map to
weighted quadruples of discrepancy at most $1+4H_{\mathrm M}$, with the same character relation.
Positivity, Gaussian Minkowski and comparison give
\[
 \mathcal E(W)\ll H_{\mathrm M}\Big(\sum_c\mathcal E(\alpha_c)^{1/4}\Big)^4
 \ll H_{\mathrm M}V^{-3}\Big(\sum_c A_{b(c)}^{1/3}p_c^{1/3}\Big)^3|W|.
\]

The tolerance enlargement contributes the first factor
$H_{\mathrm M}$: four shifts of size at most $H_{\mathrm M}$ change
the discrepancy by at most $4H_{\mathrm M}$. The second inequality
is H\"older in the class index, in the form
\begin{equation*}
 \sum_c A_{b(c)}^{1/4}p_c^{1/4}m_c^{1/4}
 \leq
 \left(\sum_c A_{b(c)}^{1/3}p_c^{1/3}\right)^{3/4}
 \left(\sum_c m_c\right)^{1/4}.
\end{equation*}
The equality $\sum_c m_c=|W|$ follows from the unit total allocation
at each original point. Nonnegativity allows the expansion into
weighted quadruples and the subsequent upper bounds.

For each box, bin tuple and color, the dyadic sum $\sum_p p^{1/3}\ll V^{1/3}$ is geometric,
including the top bin. Pad missing tuples by empty classes.

Writing the lower endpoints as $V,V/2,V/4,\ldots$, up to fixed
factors for threshold adjustment, we have
\begin{equation*}
 \sum_p p^{1/3}
 \ll V^{1/3}\sum_{j\geq0}2^{-j/3}
 \ll V^{1/3}.
\end{equation*}
Thus the number of amplitude bins is never inserted as a factor in
this estimate. The box count and the count of factor bins remain, as does
the number of colors needed for the shifts.
 For a fixed plan the last sum is at
most $A_b^{1/3}JK_bH_{\mathrm M}V^{1/3}$. There are boundedly many plans; cubing their sum costs a
fixed constant and preserves each $A_b$ paired with its own $K_b^3$. This proves \eqref{c80:eq:witness}, with adaptive sampling allowed by the
uniform estimates for separated sets.
The four powers of $H_{\mathrm M}$ in the conclusion now have separate
sources. Three come from cubing the color count in the last H\"older
sum, and one comes from the enlarged energy tolerance. Likewise,
$J^3$ and $K_b^3$ are cubic counts, while the dyadic amplitude sum
costs no further logarithm. Keeping $A_b$ attached to its own $K_b$
avoids combining the worst local cost from one plan with an unrelated
count of factor bins from another.
\end{proof}

\subsection{The physical shell sum}
The exact identity with four stages and its common Mellin integral give
the hypothesis of Proposition~\ref{c80:prop:fractional} for the
complete level normalized by Gauss sums. Indeed
\(\gamma_q(\chi)\leq1\), the Mellin weights have bounded
\(L^1\) mass, and the Gevrey tails are negligible.
Here \(J\ll L^7\) and \(H_{\mathrm M}\asymp L^s\).
The complete levels satisfy \(\sqrt X\ll V\ll X^\vartheta L\)
by Lemma~\ref{lem:endpoint-prime-cap}. Thus all adjusted product
thresholds lie in \eqref{c80:eq:strip} for a fixed sufficiently
large \(C_0\).
The required product bounds are given by
Proposition~\ref{prop:c80-terminal}.

The identity expresses each complete prime value as a linear combination,
with uniformly bounded coefficients, of Mellin integrals of the $J$ products. Taking absolute
values and using the bounded total Mellin mass makes the sum of their
shifted maxima at least a fixed multiple of the complete threshold.
Since $\gamma_q(\chi)\leq1$, a complete value normalized by Gauss sums at
level $V$ satisfies the same requirement for product selection. The negligible tails
are absorbed by fixing the truncation constant before $X$ varies.

Write $K_b\ll L^{\kappa_b}$ for the number of classes of factor values
needed for each product and amplitude level in a given plan. Taking $A_b=XB^2\mathfrak ZL^{c_b}$, the
local costs are
\[
\begin{array}{c|rrrrrr}
\text{plan}&\le4&\text{positive five}&\text{short five}&6&7&8\\ \hline
c_b&97/2&35&31&331/10&283/8&31\\
\kappa_b&0&4&0&0&0&0.
\end{array}
\]
The resulting bound for complete values has exponent
$c_b+21+3\kappa_b+4s$.  Proposition~\ref{prop:endpoint-G-bridge}
contributes $L^3\ell_2^2$, and Lemma~\ref{lem:endpoint-raw-mass}
gives the actual mass bound
\[
 \mathfrak M\ll(Q^2U+X)XL,
 \qquad U\asymp T+H_{\mathrm M},\quad q\asymp Q,\quad B_0=QT.
\]
Here the actual mass of the complete prime polynomial contributes one power of $L$.

To obtain the exponent for complete values, substitute $J\ll L^7$,
$K_b\ll L^{\kappa_b}$ and $H_{\mathrm M}\asymp L^s$ in the
estimate with fractional weights:
\begin{equation*}
 H_{\mathrm M}^4 A_bJ^3K_b^3
 \ll XB^2\mathfrak Z\,
       L^{c_b+4s+3\cdot7+3\kappa_b}.
\end{equation*}
The bins for factor values are already included in $K_b$. Scalar slices
used inside the complement and selector lemmas are not additional
product classes; their energies have been summed within those lemmas.

On physical minor shells,
$Q\ll\sqrt X$, $QT\ll\sqrt X$, $Q^2T\ll X$, $T\ge1$. The $XB^2$ part of the reconstruction estimate at the original modulus has shell factor
\[
 X\frac{Q(T+H_{\mathrm M})^2}{T}\bigl(Q^2(T+H_{\mathrm M})+X\bigr).
\]
For $T\ge H_{\mathrm M}$, this is $O(X^2QT)$; summing $T$ geometrically for each $Q$ gives
$O(X^{5/2})$ and then $O(X^{5/2}L)$ over $Q$. For $T<H_{\mathrm M}$, sum $1/T$ first, then $Q$
geometrically, to get
\[
 XH_{\mathrm M}^2\sum_{Q\ll\sqrt X}Q(Q^2H_{\mathrm M}+X)\ll X^{5/2}H_{\mathrm M}^3.
\]
For the extra part $XB^2(\mathfrak Z-1)=X^{3/2}B$, the shell factor is
\[
 X^{3/2}\frac UT(Q^2U+X).
\]
For $T\ge H_{\mathrm M}$, each shell is $O(X^{5/2})$ and all contribute $O(X^{5/2}L^2)$. For
$T<H_{\mathrm M}$, again sum $1/T$ first and then $Q$, obtaining
\[
 X^{3/2}H_{\mathrm M}\sum_{Q\ll\sqrt X}(Q^2H_{\mathrm M}+X)
 \ll X^{5/2}(H_{\mathrm M}^2+H_{\mathrm M}L).
\]

Here are the summations in more detail. For large $T$, we have
$U\asymp T$ and $Q^2T\ll X$. In the first scale part, the shell
factor is consequently at most a constant times $X^2QT$. At fixed
$Q$, the geometric sum of $QT$ is $O(\sqrt X)$, after which there
are $O(L)$ choices of $Q$. In the correction part the shell bound is
$O(X^{5/2})$, so counting the two dyadic coordinates costs at most
$L^2$.

For small $T$, the sum of $1/T$ over dyadic heights is bounded. In the
first part, the remaining modulus sums are
$\sum_QQ^3\ll X^{3/2}$ and $\sum_QQ\ll X^{1/2}$.
They give respectively $X^{5/2}H_{\mathrm M}^3$ and
$X^{5/2}H_{\mathrm M}^2$. In the correction part,
$\sum_QQ^2\ll X$, whereas $\sum_Q1\ll L$; these give
$X^{5/2}H_{\mathrm M}^2$ and
$X^{5/2}H_{\mathrm M}L$. All sums over $Q$ in this calculation
are over dyadic scales.

Since $s>1$, $L^2,H_{\mathrm M}^2,H_{\mathrm M}L\ll H_{\mathrm M}^3$. Both scale parts therefore
cost at most $X^{5/2}L^{3s}$. This includes the shell count and height extension.

The extra term has been kept because
\begin{equation*}
 XB^2\mathfrak Z=XB^2+X^{3/2}B.
\end{equation*}
It arises from the possible excess of a raw M\"obius length squared
over $B$. The two shell computations show directly that this correction
fits the same total $L^{3s}$ allowance. It is therefore not necessary
to replace every local correction by its worst possible power of $X$.

The total exponent for plan $b$ is thus
\[
 c_b+25+3\kappa_b+7s
 =c_b+32+3\kappa_b+7(s-1).
\]

Each term can be recorded at the point where it is incurred:
\begin{equation*}
 \begin{array}{c|c}
 \text{operation}&\text{power of }L\\ \hline
 \text{local product estimate}&c_b\\
 \text{cubic count of product boxes}&21\\
 \text{cubic count of factor bins}&3\kappa_b\\
 \text{shift separation and tolerance}&4s\\
 \text{reconstruction of the fourth power}&3\\
 \text{square mass of the complete prime polynomial}&1\\
 \text{physical shell sum and height enlargement}&3s
 \end{array}
\end{equation*}
Thus the fixed numerical contribution is $21+3+1=25$ and the
smoothing contribution is $4s+3s=7s$. The factor $\ell_2^2$ from
the two totient ratios is kept separately.

With $s=2$, the total logarithmic costs are
\begin{equation*}
 \frac{175}{2},\qquad86,\qquad70,\qquad
 \frac{721}{10},\qquad\frac{595}{8},\qquad70.
\end{equation*}
Their maximum is $175/2$, from the plan for at most four groups.
The scalar amplitude sums have already been included in the local
lemmas. Hence
\begin{equation*}\label{eq:terminal-fourth-moment}
 \int_{\mathrm{terminal\ minor\ cells}}
             |S_{\rm np,hi}(\alpha)|^4\,d\alpha
 \ll_\vartheta X^{5/2}L^{175/2}\ell_2^2
 \ll_\vartheta X^{5/2}L^{89}.
\end{equation*}
The constants and the strip width may depend on the fixed exponent bound
$\vartheta$, while this logarithmic exponent does not.

\section{Proofs of the main results}
\label{sec:proof-general}

\begin{proof}[Proof of Theorem~\ref{thm:moment-main}]
Put $\Theta_0=\max(\Theta,7/8)$ and $\delta=1-\Theta_0>1/12$.
Choose $\kappa$ as in Proposition~\ref{prop:ar-away}. This is
possible because its scalar threshold tends to $1/12$ as
$\kappa\uparrow9/20$, while its threshold for raw mass tends to
$9/40<38/167$.
Fix the Gevrey smoothing order $s=2$, and choose the terminal
width and all allowances by small powers in the order specified above.

Choose a fixed $\eta_1>0$ with $2\eta_1<\eta_0$. In the
argument with five stages choose a separate auxiliary exponent
$\eta_{\rm a}>0$, in place of the $\eta$ used for $\mathcal H$
there, after $\eta_1$ and the saving in the away range have been fixed.
Take it small enough that every height and exponent perturbation
is smaller than $\eta_1$ and than that saving. For physical shells
below $X^{1/2-2\eta_1}$, the enlarged parameters for five stages
then lie below $X^{1/2-\eta_1}$. Apply
Proposition~\ref{prop:ar-away} on this compact interval.
It gives a total fourth integral $O_\Theta(X^{5/2-c_1})$ for
some fixed $c_1>0$. On all remaining shells the actual enlarged
height $B=qU$ lies in \eqref{c80:eq:strip}; the fixed overlap
absorbs support constants and the extension $B\ll\sqrt X L^2$.
The argument with four stages gives
\[
 \int_{\mathrm{terminal}}|S_{\rm np,hi}(\alpha)|^4\,d\alpha
 \ll_\Theta X^{5/2}L^{175/2}\bigl(\log\log(10X)\bigr)^2
 \ll_\Theta X^{5/2}L^{89}.
\]
By Lemma~\ref{lem:endpoint-G-low-cap}, the contribution from low values
with cutoff $\sqrt X$ has fourth integral $O(X^{5/2}L)$.

Proposition~\ref{prop:principal-lower-theta} gives the complete
principal fourth integral $O_\varepsilon(X^{22/9+\varepsilon})$.
Choose its small exponent below $1/18$, so that it has a fixed
power saving relative to $X^{5/2}$. Prime powers and the
nonunit corrections have fourth integrals $O(X^2L^8)$, and
Gevrey tails have any prescribed fixed power saving. Combining
these finitely many contributions on the measurable rational
cell gives
\begin{equation*}
 \int_{\mathfrak m(X^\kappa)}|S_w(\alpha)|^4\,d\alpha
 \ll_\Theta X^{5/2}L^{89}
 \ll_\Theta X^{5/2}L^{100}.
\end{equation*}

This proves \eqref{eq:moment-main}, with $\kappa$ and every other
parameter fixed before $X$ varies.
\end{proof}

\begin{proof}[Proof of Theorem~\ref{thm:applications}]
Fix the weight $w$ once and for all and take $\kappa=\kappa(\Theta)$
from Theorem~\ref{thm:moment-main}. Write
\[
 M_+(n)=\int_{\mathfrak M(P)}S_w(\alpha)^2e(-n\alpha)\,d\alpha,
 \qquad
 M_-(h)=\int_{\mathfrak M(P)}|S_w(\alpha)|^2e(-h\alpha)\,d\alpha,
\]
and define $m_+(n)$ and $m_-(h)$ by replacing $\mathfrak M(P)$
with $\mathfrak m(P)$. Then
$R_+(n;X)=M_+(n)+m_+(n)$ and
$R_-(h;X)=M_-(h)+m_-(h)$.
Proposition~\ref{prop:major-range} gives a constant $c_w>0$ such
that, for all sufficiently large $X$,
\[
 |M_+(n)|\geq c_wX\quad(X/2\leq n\leq X,\ 2\mid n),\qquad
 |M_-(h)|\geq c_wX\quad(2\leq h\leq X/2,\ 2\mid h).
\]

Let $\mathcal E_+(X)$ be the even $n\in[X/2,X]$ that are not
sums of two primes, and let $\mathcal E_-(X)$ be the even
$h\in[2,X/2]$ for which $R_-(h;X)=0$. For
$n\in\mathcal E_+(X)$, the full coefficient vanishes, so
$m_+(n)=-M_+(n)$. The same identity with the difference sign
holds for $h\in\mathcal E_-(X)$. Bessel's inequality and
Theorem~\ref{thm:moment-main} therefore yield
\begin{align*}
 c_w^2X^2\#\mathcal E_+(X)
 &\leq\sum_{n\in\mathcal E_+(X)}|m_+(n)|^2
 \leq\int_{\mathfrak m(P)}|S_w(\alpha)|^4\,d\alpha,\\
 c_w^2X^2\#\mathcal E_-(X)
 &\leq\sum_{h\in\mathcal E_-(X)}|m_-(h)|^2
 \leq\int_{\mathfrak m(P)}|S_w(\alpha)|^4\,d\alpha.
\end{align*}
Consequently
\[
 \#\mathcal E_+(X)+\#\mathcal E_-(X)
 \ll_\Theta X^{1/2}(\log X)^{100}.
\]

For the sum problem, apply the first estimate at the scales
$X=2^{-j}N$ until a fixed bounded range is reached. Since
\[
 \sum_{\substack{j\geq0\\2^{-j}N\geq2}}(2^{-j}N)^{1/2}
      (\log(2^{-j}N))^{100}
 \ll N^{1/2}(\log N)^{100},
\]
enlarging the constant to cover the remaining integers gives the claimed bound
for $E_+(N)$.

For the difference problem, take $X=2N$. If $h\leq N$ is positive
and even and is not in $\mathcal E_-(2N)$, then
$R_-(h;2N)>0$. The nonnegativity of every summand supplies primes
$p_1<p_2<2N$ with $p_2-p_1=h$. Thus every shift counted by
$E_-(N)$ belongs to $\mathcal E_-(2N)$, and
\[
 E_-(N)\leq\#\mathcal E_-(2N)
 \ll_\Theta N^{1/2}(\log N)^{100}.
\]
This also proves the stated bound with the prime range $p_2<2N$.
\end{proof}

\section*{Acknowledgments}
The authors acknowledge the use of OpenAI's ChatGPT during the preparation
of this manuscript. 
The authors assume responsibility for all content.
Huixi Li is supported by the National Natural Science Foundation of China (Grant No. 12561001). 
Zhixin Liu is supported by  Independent Innovation Foundation of Tianjin University (Grant No. 2025XJ21-0009).


\begin{thebibliography}{10}

\bibitem{BhowmikGrimmelt2026}
G.~Bhowmik and L.~Grimmelt.
\newblock The exceptional set of the {Goldbach} problem.
\newblock {\em Analysis Mathematica}, 2026.
\newblock Published online 23 September 2026.

\bibitem{BhowmikRuzsa2018}
Gautami {Bhowmik} and Imre~Z. {Ruzsa}.
\newblock {Average Goldbach and the Quasi-Riemann Hypothesis}.
\newblock \href{https://arxiv.org/abs/1711.06442}{arXiv:1711.06442}, November
  2017.

\bibitem{CGL2026}
Bin {Chen}, Vishal {Gupta}, and Yung~Chi {Li}.
\newblock {Large Value Estimates for Dirichlet Polynomials with Characters and
  Zero Density of Dirichlet $L$-Functions}.
\newblock \href{https://arxiv.org/abs/2507.08296v2}{arXiv:2507.08296v2}, July
  2026.

\bibitem{ChenPan1980}
Jing~Ren Chen and Cheng~Dong Pan.
\newblock The exceptional set of {G}oldbach numbers. {I}.
\newblock {\em Sci. Sinica}, 23(4):416--430, 1980.

\bibitem{Chen1983}
Jing~Run Chen.
\newblock The exceptional set of {G}oldbach numbers. {II}.
\newblock {\em Sci. Sinica Ser. A}, 26(7):714--731, 1983.

\bibitem{ChenLiu1989}
Jing~Run Chen and Jianmin Liu.
\newblock The exceptional set of {G}oldbach-numbers. {III}.
\newblock {\em Chinese Quart. J. Math.}, 4(1):1--15, 1989.

\bibitem{Chen1973}
Jingrun Chen.
\newblock On the representation of a larger even integer as the sum of a prime
  and the product of at most two primes.
\newblock {\em Sci. Sinica}, 16:157--176, 1973.

\bibitem{Chudakov1937}
N.~G. Chudakov.
\newblock Sur le probl\`eme de {G}oldbach.
\newblock {\em C.R.(Dokl.) Acad. Sci. URSS}, 17:335--338, 1937.

\bibitem{Chudakov1938}
N.~G. Chudakov.
\newblock On the density of the set of even numbers which are not representable
  as a sum of two odd primes.
\newblock {\em Izvestiya Akademii Nauk SSSR. Seriya Matematicheskaya},
  2(1):25--40, 1938.

\bibitem{Estermann1938}
T.~Estermann.
\newblock On {G}oldbach's {P}roblem : {P}roof that {A}lmost all {E}ven
  {P}ositive {I}ntegers are {S}ums of {T}wo {P}rimes.
\newblock {\em Proc. London Math. Soc. (2)}, 44(4):307--314, 1938.

\bibitem{FGIS2022}
J.~B. Friedlander, D.~A. Goldston, H.~Iwaniec, and A.~I. Suriajaya.
\newblock Exceptional zeros and the {G}oldbach problem.
\newblock {\em J. Number Theory}, 233:78--86, 2022.

\bibitem{Gallagher1970Density}
P.~X. Gallagher.
\newblock A large sieve density estimate near {$\sigma =1$}.
\newblock {\em Invent. Math.}, 11:329--339, 1970.

\bibitem{Goldston1992}
D.~A. Goldston.
\newblock On {H}ardy and {L}ittlewood's contribution to the {G}oldbach
  conjecture.
\newblock In {\em Proceedings of the {A}malfi {C}onference on {A}nalytic
  {N}umber {T}heory ({M}aiori, 1989)}, pages 115--155, Salerno, 1992. Univ.
  Salerno.

\bibitem{GoldstonSuriajaya2021}
D.~A. {Goldston} and Ade~Irma {Suriajaya}.
\newblock {Note on the Goldbach Conjecture and Landau-Siegel Zeros}.
\newblock \href{https://arxiv.org/abs/2104.09407}{arXiv:2104.09407}, April
  2021.

\bibitem{GPY2009}
Daniel~A. Goldston, J\'anos Pintz, and Cem~Y. Y\i ld\i~r\i m.
\newblock Primes in tuples. {I}.
\newblock {\em Ann. of Math. (2)}, 170(2):819--862, 2009.

\bibitem{Granville2007}
Andrew Granville.
\newblock Refinements of {G}oldbach's conjecture, and the generalized {R}iemann
  hypothesis.
\newblock {\em Funct. Approx. Comment. Math.}, 37(part 1):159--173, 2007.

\bibitem{GreenTao2006}
Ben Green and Terence Tao.
\newblock Restriction theory of the {S}elberg sieve, with applications.
\newblock {\em J. Th\'eor. Nombres Bordeaux}, 18(1):147--182, 2006.

\bibitem{GreenTao2010}
Benjamin Green and Terence Tao.
\newblock Linear equations in primes.
\newblock {\em Ann. of Math. (2)}, 171(3):1753--1850, 2010.

\bibitem{Grimmelt2022}
Lasse Grimmelt.
\newblock Goldbach numbers in short intervals.
\newblock {\em Ann. Sc. Norm. Super. Pisa Cl. Sci. (5)}, 23(3):1395--1416,
  2022.

\bibitem{GuthMaynard2026}
Larry Guth and James Maynard.
\newblock New large value estimates for {D}irichlet polynomials.
\newblock {\em Ann. of Math. (2)}, 203(2):623--675, 2026.

\bibitem{HardyLittlewood1923}
G.~H. Hardy and J.~E. Littlewood.
\newblock Some problems of `{P}artitio numerorum'; {III}: {O}n the expression
  of a number as a sum of primes.
\newblock {\em Acta Math.}, 44(1):1--70, 1923.

\bibitem{HL1924}
G.~H. Hardy and J.~E. Littlewood.
\newblock Some problems of '{P}artitio numerorum'({V}): {A} further
  contribution to the study of {G}oldbach's problem.
\newblock {\em Proc. London Math. Soc. (2)}, 22:46--56, 1924.

\bibitem{HeathBrown1979Density}
D.~R. Heath-Brown.
\newblock The density of zeros of {D}irichlet's {$L$}-functions.
\newblock {\em Canadian J. Math.}, 31(2):231--240, 1979.

\bibitem{HuangLi2020}
Jing-Jing Huang and Huixi Li.
\newblock On the connection between the {G}oldbach conjecture and the
  {E}lliott-{H}alberstam conjecture.
\newblock In {\em Combinatorial and additive number theory {IV}}, volume 347 of
  {\em Springer Proc. Math. Stat.}, pages 323--346. Springer, Cham, [2021]
  \copyright 2021.

\bibitem{IK2004}
Henryk Iwaniec and Emmanuel Kowalski.
\newblock {\em Analytic number theory}, volume~53 of {\em American Mathematical
  Society Colloquium Publications}.
\newblock American Mathematical Society, Providence, RI, 2004.

\bibitem{KT2005}
Angel~V. Kumchev and D.~I. Tolev.
\newblock An invitation to additive prime number theory.
\newblock {\em Serdica Math. J.}, 31(1-2):1--74, 2005.

\bibitem{LZ2012}
Alessandro Languasco and Alessandro Zaccagnini.
\newblock The number of {G}oldbach representations of an integer.
\newblock {\em Proc. Amer. Math. Soc.}, 140(3):795--804, 2012.

\bibitem{Li1999Goldbach}
Hongze Li.
\newblock The exceptional set of {G}oldbach numbers.
\newblock {\em Quart. J. Math. Oxford Ser. (2)}, 50(200):471--482, 1999.

\bibitem{Li2000Goldbach}
Hongze Li.
\newblock The exceptional set of {G}oldbach numbers. {II}.
\newblock {\em Acta Arith.}, 92(1):71--88, 2000.

\bibitem{Li2023}
Huixi {Li}.
\newblock {Some discussions on the Goldbach conjecture}.
\newblock \href{https://arxiv.org/abs/2306.17769}{arXiv:2306.17769}, June 2023.

\bibitem{Liu2026}
Baiying Liu.
\newblock Improved zero-free half-planes for the quasi-{Riemann} hypothesis.
\newblock \href{https://arxiv.org/abs/2610.12234v2}{arXiv:2610.12234v2},
  October 2026.
\newblock Version 2, 12 October 2026.

\bibitem{Lu2010}
Wen~Chao Lu.
\newblock Exceptional set of {G}oldbach number.
\newblock {\em J. Number Theory}, 130(10):2359--2392, 2010.

\bibitem{Matomaki2008}
Kaisa Matom\"{a}ki.
\newblock On the exceptional set in {G}oldbach's problem in short intervals.
\newblock {\em Monatsh. Math.}, 155(2):167--189, 2008.

\bibitem{MatomakiMerikoski2023}
Kaisa Matom\"aki and Jori Merikoski.
\newblock Siegel zeros, twin primes, {G}oldbach's conjecture, and primes in
  short intervals.
\newblock {\em Int. Math. Res. Not. IMRN}, (23):20337--20384, 2023.

\bibitem{MRT2019}
Kaisa Matom\"aki, Maksym Radziwi\l\l, and Terence Tao.
\newblock Correlations of the von {M}angoldt and higher divisor functions {I}.
  {L}ong shift ranges.
\newblock {\em Proc. Lond. Math. Soc. (3)}, 118(2):284--350, 2019.

\bibitem{Maynard2015}
James Maynard.
\newblock Small gaps between primes.
\newblock {\em Ann. of Math. (2)}, 181(1):383--413, 2015.

\bibitem{MPR2026}
James {Maynard}, Mayank {Pandey}, and Maksym {Radziwi{\l}{\l}}.
\newblock {Exponential sums over primes}.
\newblock \href{https://arxiv.org/abs/2608.14777v1}{arXiv:2608.14777v1}, August
  2026.

\bibitem{Mikawa1991}
Hiroshi Mikawa.
\newblock On prime twins.
\newblock {\em Tsukuba J. Math.}, 15(1):19--29, 1991.

\bibitem{MontgomeryVaughan1975}
H.~L. Montgomery and R.~C. Vaughan.
\newblock The exceptional set in {G}oldbach's problem.
\newblock {\em Acta Arith.}, 27:353--370, 1975.

\bibitem{Montgomery1971}
Hugh~L. Montgomery.
\newblock {\em Topics in multiplicative number theory}, volume 227 of {\em
  Lecture Notes in Mathematics}.
\newblock Springer-Verlag, Berlin--New York, 1971.

\bibitem{OAI2026September30}
{OpenAI}.
\newblock The {Quasi-Riemann Hypothesis}: A zero-free half-plane
  $\mathrm{Re}(s)>7/8$.
\newblock OpenAI Math Release preprint,
  \href{https://github.com/openai/math/blob/main/preprints/The-Quasi-Riemann-Hypothesis-September-30-2026/paper.pdf}{30
  September 2026}, 2026.

\bibitem{PP1993}
A.~Perelli and J.~Pintz.
\newblock On the exceptional set for {G}oldbach's problem in short intervals.
\newblock {\em J. London Math. Soc. (2)}, 47(1):41--49, 1993.

\bibitem{Pintz2009}
J\'anos Pintz.
\newblock Landau's problems on primes.
\newblock {\em J. Th\'eor. Nombres Bordeaux}, 21(2):357--404, 2009.

\bibitem{Pintz2018II}
Janos {Pintz}.
\newblock {A new explicit formula in the additive theory of primes with
  applications II. The exceptional set in Goldbach's problem}.
\newblock \href{https://arxiv.org/abs/1804.09084}{arXiv:1804.09084}, April
  2018.

\bibitem{Pintz2018I}
J\'anos Pintz.
\newblock A new explicit formula in the additive theory of primes with
  applications {I}. {T}he explicit formula for the {G}oldbach problem and the
  generalized twin prime problem.
\newblock {\em Acta Arith.}, 210:53--94, 2023.
\newblock Numbering in the text follows arXiv:1804.05561v1.

\bibitem{PintzAlthofer2026}
J{\'a}nos {Pintz} and Ingo {Alth{\"o}fer}.
\newblock The {Goldbach} exceptional set: an exponent of 7/12 and conditional
  refinements.
\newblock \href{https://arxiv.org/abs/2610.15559v1}{arXiv:2610.15559v1},
  October 2026.

\bibitem{Polymath2014}
D.~H.~J. Polymath.
\newblock Variants of the {S}elberg sieve, and bounded intervals containing
  many primes.
\newblock {\em Res. Math. Sci.}, 1:Art. 12, 83, 2014.

\bibitem{Titchmarsh1986Zeta}
E.~C. Titchmarsh.
\newblock {\em The Theory of the {Riemann} Zeta-Function}.
\newblock Clarendon Press, Oxford, second edition, 1986.
\newblock Revised by D. R. Heath-Brown.

\bibitem{VanDerCorput1937}
J.~G. van~der Corput.
\newblock Sur l'hypoth\`ese de {G}oldbach pour presque tous les nombres pairs.
\newblock {\em Acta Arith}, 2:266--290, 1937.

\bibitem{Vaughan1972}
R.~C. Vaughan.
\newblock On {G}oldbach's problem.
\newblock {\em Acta Arith.}, 22:21--48, 1972.

\bibitem{Vaughan2016}
R.~C. Vaughan.
\newblock Goldbach's conjectures: a historical perspective.
\newblock In {\em Open problems in mathematics}, pages 479--520. Springer,
  [Cham], 2016.

\bibitem{XuSang2026}
Yuanda Xu and Hejian Sang.
\newblock Amplified moments and a refinement of the quasi-{Riemann} hypothesis.
\newblock \href{https://arxiv.org/abs/2610.15291v2}{arXiv:2610.15291v2},
  October 2026.
\newblock Version 2, 15 October 2026.

\bibitem{Zhang2014}
Yitang Zhang.
\newblock Bounded gaps between primes.
\newblock {\em Ann. of Math. (2)}, 179(3):1121--1174, 2014.

\bibitem{Zhao2026}
Genheng {Zhao}.
\newblock {The exceptional set of Goldbach problem and Linnik's constant}.
\newblock \href{https://arxiv.org/abs/2511.05631v2}{arXiv:2511.05631v2},
  January 2026.

\end{thebibliography}
\end{document}